\documentclass[11pt,reqno]{amsart}
\usepackage[T1]{fontenc}
\usepackage{lmodern,amsmath,amssymb,mathtools}
\usepackage[hmargin=1.08in,vmargin=1in]{geometry}
\usepackage{microtype,needspace}
\usepackage[hidelinks]{hyperref}

\newtheorem{theorem}{Theorem}[section]
\newtheorem{proposition}[theorem]{Proposition}
\newtheorem{corollary}[theorem]{Corollary}
\newtheorem{lemma}[theorem]{Lemma}
\theoremstyle{remark}

\numberwithin{equation}{section}

\newcommand{\per}{\operatorname{per}}
\newcommand{\imm}{\operatorname{imm}}
\newcommand{\rk}{\operatorname{rank}}
\newcommand{\tr}{\operatorname{tr}}
\newcommand{\diag}{\operatorname{diag}}
\newcommand{\calI}{\mathcal I}
\newcommand{\calD}{\mathcal D}
\newcommand{\F}{\mathbb F}
\newcommand{\Z}{\mathbb Z}
\newcommand{\N}{\mathbb N}
\newcommand{\cH}{\mathcal H}
\newcommand{\cL}{\mathcal L}
\newcommand{\ind}{\mathbf 1}
\newcommand{\eul}[2]{\genfrac{\langle}{\rangle}{0pt}{}{#1}{#2}}

\title[Closed formulas for permanents]{Closed formulas for the permanents of structured matrices}
\author{M. Fazeel Anwar}
\address{Department of Mathematics, Sukkur IBA University, Sukkur 65200, Pakistan}
\email{fazeel.anwar@iba-suk.edu.pk}
\hypersetup{pdftitle={Closed formulas for the permanents of structured matrices},
            pdfauthor={M. Fazeel Anwar}}

\begin{document}
\raggedbottom

\begin{abstract}
We give finite formulas for permanents of structured matrices and
prove polynomial bounds for their evaluation. Over any field, the
permanent of a matrix of bounded rank-width can be computed from
the matrix alone, answering a question of Lyaudet. In characteristic
\(p>0\), the bound is \(p^{O(k)}n^{O(1)}\), where \(k\) is the
rank-width. In characteristic zero, a determinant formula treats
Cauchy matrices with an added polynomial of fixed degrees in the
row and column nodes. This includes polynomial perturbations of
Hilbert matrices and dense Toeplitz and Hankel families. For
ordered cuts we give a single determinant; in positive
characteristic, rational diagonal sequences give a fixed rational
generating function for the Toeplitz permanents. We also give
formulas for cycle weights and hook immanants. For translation
matrices we obtain signed determinant formulas and classify the
supports for which the permanent is the coefficientwise absolute
value of the determinant. Further formulas cover cyclic block
matrices, residue supports, and bands bounded on one side.
\end{abstract}
\maketitle

\section{Introduction}

Let \(\F\) be a field. We assume characteristic zero except where
otherwise stated. For an \(n\times n\)
matrix \(A=(a_{ij})\) over \(\F\), the permanent of \(A\) is
\[
 \per A=\sum_{\sigma\in S_n}\prod_{i=1}^n a_{i,\sigma(i)},
\]
where \(S_n\) is the symmetric group on \(\{1,\ldots,n\}\). The
permanent of the empty matrix is one.

Barvinok \cite{Barvinok} proved that the permanent of a matrix of
fixed rank can be computed in polynomial time. Flarup and Lyaudet
\cite{FlarupLyaudet} treated matrices of bounded weighted clique-width.

A cut \(I\mid I^c\) is one partition of the indices applied to
both rows and columns. Its width is the sum of the ranks of
\(A[I,I^c]\) and \(A[I^c,I]\). The rank-width is the smallest
possible maximum cut width over binary trees on the indices. In \cite[\S~4]{LyaudetRankWidth}, Lyaudet asked whether the permanent of a matrix of bounded rank-width be computed
in polynomial time?

Theorem~\ref{thm:tree} gives an affirmative answer over any field.
For each fixed bound on the rank-width, it computes the permanent
from the entries of the matrix alone. The proof finds a tree of
bounded width from the cut ranks and evaluates a finite sum by
combining coefficient vectors at the nodes. For rational entries
the algorithm has polynomial bit complexity. In characteristic
\(p>0\), Theorem~\ref{thm:positive} gives the bound
\(p^{O(k)}n^{O(1)}\), including the construction of a suitable tree.
Thus, for fixed \(p\), rank-width \(O(\log(n+1))\) is also allowed.

Cancellation in prime characteristic has been studied in terms of
the rank of the whole matrix. Yu \cite{Yu} proved the bound that
\(\per A=0\) if \(n>(p-1)\rk A\). Glynn \cite{Glynn} gave another
proof and a factorization at the boundary rank. Arvind, Chatterjee,
Datta and Mukhopadhyay \cite[Corollary 4.7]{ACDM} gave a parameterized
algorithm over finite fields for matrices of bounded rank.
Here the rank bounds concern submatrices across a tree, and the
field may be infinite. The proof reduces the endpoint polynomials
modulo the \(p\)th powers of their variables. For rational Toeplitz
diagonal sequences this gives a fixed rational generating function
for the permanents.

The ranks of the off-diagonal submatrices at successive cuts are
used in the theory of quasiseparable matrices
\cite{PernetStorjohann,PernetSignargoutVillard}. Bounding these two
ranks allows the total rank to grow with the order of the matrix. For
example, the matrix \((2^{-|i-j|})_{i,j=1}^n\) has full rank, and both
ranks at each nontrivial cut in the natural ordering are one.

Borchardt's identity \cite{Borchardt} expresses the permanent of a
Cauchy matrix as a quotient of two determinants. Han and
Krattenthaler \cite{HK} treated rectangular Cauchy permanents, and
Khare and Sahi \cite{KS} obtained permanent and immanant identities
for entrywise power series. These formulas do not directly give the
evaluation bound considered here for
\[
 A_{ij}=\frac1{x_i-y_j}+P(x_i,y_j),
\]
where the two degrees of \(P\) are fixed and its coefficients and
the nodes may vary with \(n\). Theorem~\ref{thm:cauchy} gives a finite
formula with polynomially many determinants, each of order at most
\(2n\). Its proof expands a bordered Cauchy permanent according to
the blocks and then uses polynomial interpolation. The resulting
families need not have bounded rank-width.

For circulants with fixed jumps, Minc \cite{Minc1985,Minc1987}
gave recurrences for the permanent. Golin, Leung and Wang
\cite{GolinLeungWang} gave a transfer matrix proof and treated
jumps that vary linearly with the order. Codenotti and Resta
\cite{CodenottiRestaSparse} expressed the permanent of a
three-shift circulant as a linear combination of four
determinants and treated a class with four shifts. For further
formulas for sparse Toeplitz and banded matrices, see
\cite{CodenottiToeplitz,RestaSburlati,Schwartz}. Alekseyev and
Khomovsky \cite{AK1,AK2} give fixed-band recurrences and their
symmetry reductions. The coefficient
signs for three-shift circulants were studied in
\cite{LWW,DePoiSparse}. Rodtes and Anwar \cite{RodtesAnwar} gave
a formula for the permanent of a block matrix in terms of
submatrices of its blocks. For a circulant with four arbitrary
shifts and one weight per shift, the question of a closed formula
that can be evaluated in polynomial time uniformly in the shifts
remains open.

Our main results are the following. Theorems~\ref{thm:tree}
and~\ref{thm:positive} give the tree formula and its refinement in
positive characteristic. Theorem~\ref{thm:cauchy} treats polynomial
perturbations of Cauchy matrices in characteristic zero. These are
the two general constructions used for the permanent. For cuts
along an ordering, Theorem~\ref{thm:path} and
Corollary~\ref{cor:positive-ordered} give a single determinant.
The subsequent results treat prescribed cycle weights and hook
immanants. Theorems~\ref{thm:circulant}, \ref{thm:block},
\ref{thm:toeplitz} and~\ref{thm:hankel} collect the formulas for
the corresponding matrix classes. Their corollaries give the
specializations. In each polynomial-time assertion, the parameters
that must be fixed are stated explicitly.

\section{Preliminaries}\label{sec:preliminaries}

Write \([n]=\{1,\ldots,n\}\) and \(\N=\{0,1,2,\ldots\}\).
For sets of row and column indices \(I,J\), \(A[I,J]\) is the
submatrix in the given order. We use \(I^c\) for the complement
in the relevant index set, and \([w]f\) for the coefficient
of the monomial \(w\) in \(f\). The determinant of the empty
matrix is one. Permuting rows or columns preserves the permanent.

For \(I\subseteq[n]\), use the partition \(I\mid I^c\) for the
row labels and for the column labels. The two crossing blocks are
\(A[I,I^c]\) and \(A[I^c,I]\). In an ordering, we use
\(I_i=\{1,\ldots,i\}\).
In a rooted binary tree on the indices, each leaf is an index
and each node \(I\) is the set of leaves below it. Write
\begin{equation}\label{eq:rank-width}
 w_A(I)=\rk A[I,I^c]+\rk A[I^c,I],
 \qquad
 \operatorname{rw}_{\F}(A)=\min_{\mathcal T}\max_{I\in\mathcal T}w_A(I).
\end{equation}
The minimum is over all such trees on \([n]\). The ranks are taken
over the field \(\F\) of the entries of \(A\). In the terminology
of Kant\'e and Rao \cite[Definitions 4.1 and 4.3]{KanteRao}, \(w_A(I)\) is the
bi-cut-rank and \(\operatorname{rw}_{\F}(A)\) is the
\(\F\)-bi-rank-width of the directed graph with an arc \(i\to j\)
of weight \(a_{ij}\) whenever \(i\ne j\) and \(a_{ij}\ne0\).
For a symmetric matrix the crossing blocks are
transposes, so this rank-width is twice the width defined using one
cut block. In formulas involving cut coordinates, \(r,s\)
denote upper bounds on the two ranks. A cut block of rank zero has
empty coordinates. In the Cauchy formulas, the locally specified
\(r,s\) instead bound the two polynomial degrees.

For each node choose rank factorizations
\begin{equation}\label{eq:tree-factors}
 A[I,I^c]=R_I F_I,\qquad A[I^c,I]^{\mathsf T}=C_I G_I .
\end{equation}
Here \(R_I,C_I\) have full column rank and rows indexed by \(I\).
The matrices \(F_I,G_I\) have full row rank and columns indexed by
\(I^c\). Put \(r_I=\rk R_I\) and \(s_I=\rk C_I\).
Let \(x_I\) and \(y_I\) be columns of \(r_I\) and \(s_I\)
indeterminates, respectively. The row and column
coordinates at \(i,j\in I\) are the linear forms
\[
 R_I(i,:)x_I,\qquad C_I(j,:)y_I .
\]
The variables at the root are empty.

Suppose \(U,V\) are the two children of \(I\). Choose matrices \(B_{UV}\)
and \(B_{VU}\) satisfying
\begin{equation}\label{eq:coupling}
 A[U,V]=R_U B_{UV} C_V^{\mathsf T},\qquad
 A[V,U]=R_V B_{VU} C_U^{\mathsf T}.
\end{equation}
They exist because the columns and rows of each cross block lie
in the respective spaces of the two factors. They are obtained
using left inverses of the two coordinate matrices.
For \(J=U,V\) define \(S_{J,I},T_{J,I}\) by
\[
 F_J[:,I^c]=S_{J,I}F_I,\qquad
 G_J[:,I^c]=T_{J,I}G_I .
\]
The rows of each child factor span the corresponding cut row
space. Restriction to the parent complement puts them in the
row spaces of \(F_I,G_I\). Write \(\rho_I\) for the substitution
\(x_J\mapsto S_{J,I}x_I,\ y_J\mapsto T_{J,I}y_I\), for both children.

The coordinate changes satisfy
\(R_J S_{J,I}=R_I[J,:]\) and \(C_J T_{J,I}=C_I[J,:]\).
Indeed, factor the corresponding submatrix with columns in \(I^c\)
in the two ways and use the full row rank of \(F_I\) or \(G_I\).

\subsection{}

For this subsection let \(\operatorname{char}\F=p>0\). For
nonnegative integers \(r,s\), put
\begin{equation}\label{eq:positive-basis}
 \mathcal B_{p,r,s}=
 \{x^\alpha y^\beta:0\leq\alpha_i,\beta_j<p,\
                          |\alpha|=|\beta|\},\qquad
 D_{p,r,s}=|\mathcal B_{p,r,s}|.
\end{equation}
Empty lists of variables are allowed. Replacing \(\beta_j\) by
\(p-1-\beta_j\) in the equality of the two degrees gives
\begin{equation}\label{eq:positive-dimension}
 D_{p,r,s}=[t^{s(p-1)}](1+t+\cdots+t^{p-1})^{r+s}
          \leq p^{r+s}.
\end{equation}
At a tree node \(I\), calculations will take place in
\begin{equation}\label{eq:positive-quotient}
 \mathcal R_I=
 \F[x_I,y_I]/(x_{I,1}^{p},\ldots,x_{I,r_I}^{p},
                         y_{I,1}^{p},\ldots,y_{I,s_I}^{p}).
\end{equation}
Each class has a unique representative with every exponent less
than \(p\). The monomials in \(\mathcal B_{p,r_I,s_I}\) are its
balanced monomials. For a scalar \(b\) and two variables define
\begin{equation}\label{eq:positive-pair}
 E_{b;x,y}=\sum_{m=0}^{p-1}\frac{b^m}{m!}
                         \partial_x^m\partial_y^m.
\end{equation}
All denominators in this formula are nonzero in \(\F\).
Lemma~\ref{lem:positive-ideal} proves that these operations and the
coordinate substitutions are well defined in the quotients.

\subsection{}

In this and the next subsection the field has characteristic zero,
\(n\geq1\), and \(r,s\geq1\). Let the row nodes
\(x_1,\ldots,x_n\) be distinct, let the column nodes
\(y_1,\ldots,y_n\) be distinct, and assume \(x_i\ne y_j\).
Write \(C=((x_i-y_j)^{-1})\). Borchardt's identity is
\begin{equation}\label{eq:cauchy-borchardt}
 \per C=\frac{\det((x_i-y_j)^{-2})}{\det((x_i-y_j)^{-1})}.
\end{equation}
The denominator is nonzero by the Cauchy determinant formula.
The following notation describes the coefficients used for adding
a polynomial to this kernel.
For $a\in\N^r$ and $b\in\N^s$, write
$|a|=\sum_\alpha a_\alpha$, $a!=\prod_\alpha a_\alpha!$, and
$z^a=\prod_\alpha z_\alpha^{a_\alpha}$, with the same conventions for $b$.
Inequalities between multiindices are componentwise. Put
\begin{equation}\label{eq:cauchy-indexset}
 \calI_n=\{(a,b)\in\N^r\times\N^s:|a|=|b|\leq n\},\qquad
 D_n=|\calI_n|
 =\sum_{h=0}^n\binom{h+r-1}{r-1}\binom{h+s-1}{s-1}.
\end{equation}
In particular, $D_n=O_{r,s}(n^{r+s-1})$.

For an $r\times s$ matrix $L$, define
\begin{equation}\label{eq:cauchy-edef}
 e_L(a,b)=[z^a w^b]
 \prod_{\alpha=1}^r\prod_{\beta=1}^s
 \left(\sum_{m=0}^n\frac{L_{\alpha\beta}^m}{m!}
 z_\alpha^m w_\beta^m\right),\qquad (a,b)\in\calI_n.
\end{equation}
This is a coefficient of a finite product. Equivalently,
\begin{equation}\label{eq:cauchy-etable}
 e_L(a,b)=
 \sum_{\substack{M\in\N^{r\times s}\!,\ M\mathbf1=a\!,\ M^{\mathsf T}\mathbf1=b}}
 \prod_{\alpha,\beta}\frac{L_{\alpha\beta}^{M_{\alpha\beta}}}{M_{\alpha\beta}!}.
\end{equation}
Here and below $0^0=1$. Formal exponentials will be used in proofs
to abbreviate \eqref{eq:cauchy-edef}; no infinite series has to be evaluated.

\subsection{}
Choose auxiliary column nodes $\zeta_1,\ldots,\zeta_r$ and auxiliary
row nodes $\eta_1,\ldots,\eta_s$. We call this choice admissible if
the entries of $(x,\eta)$ are distinct, the entries of $(y,\zeta)$
are distinct, and the two lists are disjoint. Put
\begin{equation}\label{eq:cauchy-uv}
 U_{i\alpha}=\frac1{x_i-\zeta_\alpha},\qquad
 V_{\beta j}=\frac1{\eta_\beta-y_j},\qquad
 W_{\beta\alpha}=\frac1{\eta_\beta-\zeta_\alpha}.
\end{equation}

For $(a,b)\in\calI_n$, form a row list consisting of $(x_i,0)$,
$1\leq i\leq n$, followed, for each $\beta$, by
\[
 (\eta_\beta,0),\ldots,(\eta_\beta,b_\beta-1).
\]
If $b_\beta=0$, this part of the list is empty. Form the column list
in the same way from $(y_j,0)$ and
$(\zeta_\alpha,0),\ldots,(\zeta_\alpha,a_\alpha-1)$.
For $q=1,2$, define a matrix on these lists by
\begin{equation}\label{eq:cauchy-kdef}
 K_q(a,b)_{(u,\ell),(v,m)}
 =(-1)^\ell\frac{(q+\ell+m-1)!}{(q-1)!\,\ell!\,m!}
 (u-v)^{-q-\ell-m}.
\end{equation}
Its order is $n+|a|\leq2n$. Define
\begin{equation}\label{eq:cauchy-rdef}
 R_{a,b}=\frac{\det K_2(a,b)}{\det K_1(a,b)}.
\end{equation}
Lemma~\ref{lem:cauchy-confluent} proves that the denominator is nonzero and
identifies this quotient as a permanent with repeated nodes.

The following finite sums will be used in the main formula:
\begin{align}
 Q_{a,b}
 &=\sum_{\substack{c\leq a,\ d\leq b\\|c|=|d|}}
 e_{-W^{\mathsf T}}(c,d)
 \frac{R_{a-c,b-d}}{(a-c)!\,(b-d)!},\label{eq:cauchy-qcoef}\\
 \calD(x,y;\eta,\zeta;L)
 &=\sum_{(a,b)\in\calI_n}a!b!\,e_L(a,b)Q_{a,b}.
 \label{eq:cauchy-dformula}
\end{align}
Thus $\calD$ is explicitly defined by determinants, factorials and
finite sums. Its definition contains no unevaluated permanent.

\subsection{}

For successive matrices \(W_1,\ldots,W_n\) of the same order,
let \(W\) have blocks indexed by \(0,\ldots,n\), with
\(W_{i-1,i}=W_i\) and all other blocks zero. If \(e_0,e_n\)
are chosen coordinate vectors in blocks \(0,n\), then
\begin{equation}\label{eq:product-border}
 \det\begin{pmatrix}I-W&e_n\\-e_0^{\mathsf T}&0\end{pmatrix}
    = e_0^{\mathsf T}(I-W)^{-1}e_n.
\end{equation}
Indeed \(\det(I-W)=1\), and the Schur complement gives the
identity. The \((0,n)\) block of the inverse is \(W_1\cdots W_n\),
by the finite expansion \(I+W+\cdots+W^n\).

\section{Main results}\label{sec:results}

\subsection{}

Let \(\mathcal T\) be a rooted binary tree on \([n]\), with cut
coordinates as in Section~\ref{sec:preliminaries}.
For a node \(I\) with children \(U,V\), let
\[
 \mathcal D_I=
 \sum_{\alpha,\beta}(B_{UV})_{\alpha\beta}
       \partial_{x_{U,\alpha}}\partial_{y_{V,\beta}}
 +\sum_{\alpha,\beta}(B_{VU})_{\alpha\beta}
       \partial_{x_{V,\alpha}}\partial_{y_{U,\beta}} .
\]
At node \(I\), let \(\mathcal M_I\) be the monomials in \(x_I,y_I\)
of bidegree \((h,h)\), \(0\leq h\leq |I|\).
For \(u\in\mathcal M_U,\ v\in\mathcal M_V,\ w\in\mathcal M_I\), the
entry in characteristic zero is
\begin{equation}\label{eq:tree-entries}
 t_I(u,v;w)=[w]\,\rho_I\!\left(
           \sum_{\ell=0}^{|I|}\frac{\mathcal D_I^\ell}{\ell!}\,uv\right).
\end{equation}
Over any field the expression inside \(\rho_I\) can be formed
using the integer coefficients in \eqref{eq:integer-merge}.
At a leaf \(\{i\}\) put
\[
 l_i(w)=[w]\big(a_{ii}+(R_{\{i\}}(i,:)x_{\{i\}})
                              (C_{\{i\}}(i,:)y_{\{i\}})\big).
\]
The permanent is obtained by summing the contributions from all
nodes, with the constant monomial at the root.

\Needspace{5\baselineskip}
\begin{theorem}[Answer to Lyaudet's question]\label{thm:tree}
Let \(k\) be a fixed nonnegative integer and \(A\) a square matrix
over a field \(\F\). If \(\operatorname{rw}_{\F}(A)\leq k\), then
\(\per A\) can be computed from the entries of \(A\) alone in
polynomial time, with field operations and zero tests counted at
unit cost.

For a tree \(\mathcal T\) as above,
\begin{equation}\label{eq:tree-formula}
 \per A=
 \sum_{\substack{w_I\in\mathcal M_I\ \text{for all nodes }I\\
                 w_{\{1,\ldots,n\}}=1}}
 \prod_{i=1}^n l_i(w_{\{i\}})
 \prod_{\substack{I\ \text{with}\\\text{children }U,V}}
       t_I(w_U,w_V;w_I).
\end{equation}
For fixed \(r,s\), the sum can be evaluated in polynomial time
from \(A\) and \(\mathcal T\). For rational entries both finding
the tree and evaluating the formula have polynomial bit complexity.
\end{theorem}

To evaluate \eqref{eq:tree-formula}, let \(p_I(w)\) be the
coefficient at node \(I\). Start with \(p_{\{i\}}=l_i\) and use
\[
 p_I(w)=\sum_{u\in\mathcal M_U,\ v\in\mathcal M_V}
             p_U(u)p_V(v)t_I(u,v;w).
\]
The answer is \(p_{[n]}(1)\). Each coefficient records partial
matchings inside \(I\), with linear factors for their unmatched
rows and columns. A matching on \(I\) has unique restrictions to
its two children and a unique set of edges between them. This is
the counting identity behind the formula. The polynomial-time
bound comes from the number of coefficients at a node; it does
not require listing all assignments in the displayed sum.
The full argument is given in Section~\ref{sec:proof-tree}.

\subsection{}

Assume \(\operatorname{char}\F=p>0\). Use the cut coordinates of
\eqref{eq:tree-factors} and \eqref{eq:coupling}, and write
\(\mathcal B_I=\mathcal B_{p,r_I,s_I}\). At a node with children
\(U,V\), put
\[
 E_I=
 \prod_{a=1}^{r_U}\prod_{b=1}^{s_V}
 E_{(B_{UV})_{ab};x_{U,a},y_{V,b}}
 \prod_{a=1}^{r_V}\prod_{b=1}^{s_U}
 E_{(B_{VU})_{ab};x_{V,a},y_{U,b}}.
\]
The factors commute. For \(u\in\mathcal B_U\),
\(v\in\mathcal B_V\) and \(w\in\mathcal B_I\), define
\[
 t_I^{(p)}(u,v;w)=[w]\rho_I\bigl(E_I(uv)\bigr)
\]
in the quotient \eqref{eq:positive-quotient}. At a leaf use the
same polynomial \(l_i\) as in Theorem~\ref{thm:tree}, reduced in
this quotient. The next theorem bounds the number of coefficients
independently of the number of leaves below the node.

\begin{theorem}\label{thm:positive}
Let \(A\) have order \(n\) over a field of characteristic \(p>0\).
For every rooted binary tree on its indices,
\begin{equation}\label{eq:positive-tree}
 \per A=
 \sum_{\substack{w_I\in\mathcal B_I\ (I\in\mathcal T)\\w_{[n]}=1}}
 \prod_{i=1}^n l_i(w_{\{i\}})
 \prod_{\substack{I\in\mathcal T\\I\text{ has children }U,V}}
                t_I^{(p)}(w_U,w_V;w_I).
\end{equation}
At a node of width at most \(k\), the formula uses at most
\(p^k\) coefficients. If \(\operatorname{rw}_{\F}(A)\leq k\),
the permanent can be computed from \(A\) alone in
\(p^{O(k)}n^{O(1)}\) field operations and zero tests. The
exponent of \(n\) is independent of \(p\) and \(k\).
For fixed \(p\), this is polynomial when
\(k=O(\log(n+1))\).
\end{theorem}

The field need not be finite. The displayed sum is evaluated by
combining the coefficient vectors at successive nodes, as in
Theorem~\ref{thm:tree}. A tree supplied with width \(k\) can be
used directly. From the matrix alone, it suffices to find a tree
of width at most \(4k\); this does not change the stated bound.

\subsection{}

In characteristic zero, consider
\begin{equation}\label{eq:cauchy-mainmatrix}
 A_{ij}=\frac1{x_i-y_j}+P(x_i,y_j),\qquad
 P(X,Y)=\sum_{a=0}^{r-1}\sum_{b=0}^{s-1}H_{ab}X^aY^b.
\end{equation}
The nodes satisfy the distinctness and disjointness assumptions
in the preliminaries. The positive integers \(r,s\) bound the
degrees in the two variables. We specify the interpolation data
for the formula.
Take the rational numbers
\[
 \xi_\alpha=\alpha\quad(1\leq\alpha\leq r),\qquad
 \upsilon_\beta=r+\beta\quad(1\leq\beta\leq s).
\]
Define the constant matrices
\[
 (B_x)_{a\alpha}=-\xi_\alpha^{-a-1}\quad(0\leq a<r),\qquad
 (B_y)_{\beta b}=\upsilon_\beta^{-b-1}\quad(0\leq b<s).
\]
They are invertible Vandermonde matrices with nonzero column or row
factors. For $t\ne0$, put
\begin{align}
 \zeta_\alpha(t)&=\xi_\alpha/t,&
 \eta_\beta(t)&=\upsilon_\beta/t,\notag\\
 T_r(t)&=\diag(t^{-1},\ldots,t^{-r}),&
 T_s(t)&=\diag(t^{-1},\ldots,t^{-s}),\notag\\
 H(t)&=B_x^{-1}T_r(t)H T_s(t)B_y^{-1}.
 \label{eq:cauchy-transformed}
\end{align}
The $H$ in \eqref{eq:cauchy-transformed} is the coefficient matrix in
\eqref{eq:cauchy-mainmatrix}, in increasing order of the two exponents.

Put
\begin{equation}\label{eq:cauchy-degrees}
 m_0=n(r+s-2),\qquad d_0=n(2r+2s-2),
\end{equation}
and
\begin{equation}\label{eq:cauchy-denpoly}
 q(t)=\prod_{i=1}^n\prod_{\alpha=1}^r(x_it-\xi_\alpha)
       \prod_{\beta=1}^s\prod_{j=1}^n(\upsilon_\beta-y_jt).
\end{equation}
In particular $q(0)\ne0$. Choose $d_0+1$ distinct nonzero rational
numbers $t_0,\ldots,t_{d_0}$ for which the auxiliary nodes in
\eqref{eq:cauchy-transformed} are admissible. They can be chosen as positive
integers of size $O_{r,s}(n)$. Indeed, there are at most
$2n(r+s)$ forbidden nonzero values. Finally put
\begin{equation}\label{eq:cauchy-lagrange}
 \lambda_\ell=[t^{m_0}]
 \prod_{\substack{0\leq j\leq d_0\\j\ne\ell}}
 \frac{t-t_j}{t_\ell-t_j}.
\end{equation}
This is the coefficient of an explicitly given polynomial of degree $d_0$.

The next theorem evaluates the permanent after a polynomial of fixed
degree has been added to a Cauchy kernel.

\begin{theorem}\label{thm:cauchy}
For the matrix $A$ in \eqref{eq:cauchy-mainmatrix},
\begin{equation}\label{eq:cauchy-mainformula}
 \boxed{\displaystyle
 \per A=\frac1{q(0)}\sum_{\ell=0}^{d_0}
 \lambda_\ell t_\ell^{m_0}q(t_\ell)
 \calD\bigl(x,y;\eta(t_\ell),\zeta(t_\ell);H(t_\ell)\bigr).}
\end{equation}
For fixed $r,s$, this formula can be evaluated in
\begin{equation}\label{eq:cauchy-timebound}
 O_{r,s}\bigl(n^{r+s+3}+n^{2r+2s-1}\bigr)
\end{equation}
field operations and zero tests. Every determinant has order at most
$2n$. For rational input, the bit complexity is polynomial in $n$
and the input bit length.
\end{theorem}

Only the bounds on the degrees are fixed here. The nodes and all the
coefficients of $P$ may vary with $n$. Formula \eqref{eq:cauchy-mainformula}
is nonrecursive. Its auxiliary evaluation points are chosen by the
stated finite rule.

\subsection{}

Suppose the two ranks at the initial interval cuts are bounded
by \(r,s\). We use the following generators for the upper and
lower entries.
\begin{equation}\label{eq:generators}
 \begin{split}
 a_{ij}&=p_i^{\mathsf T} B_{i+1}\cdots B_{j-1}q_j,\\
 a_{ji}&=u_i^{\mathsf T} C_{i+1}\cdots C_{j-1}v_j
                 \qquad (i<j).
 \end{split}
\end{equation}
The first line uses vectors and square matrices of order \(r\), and
the second uses order \(s\). An empty product is the identity.

Write \(d_i=a_{ii}\). Put \(N=\lfloor n/2\rfloor\) and let
\[
 \mathcal M=\{x^\alpha y^\beta:|\alpha|=|\beta|\leq N\},
 \qquad D=|\mathcal M|,
\]
where \(x,y\) have \(r,s\) variables. Zero generators may be used on either side. Denote by \(\Pi\) the
projection onto the span of \(\mathcal M\). Define
\[
 \begin{gathered}
 D_i=q_i^{\mathsf T}\partial_x,\quad E_i=v_i^{\mathsf T}\partial_y,
 \qquad L_i=p_i^{\mathsf T}x,\quad M_i=u_i^{\mathsf T}y,\\
 \Phi_i f=f(B_i x,C_i y),\\
 \mathcal U_i f=\Pi\big(
 \Phi_i(D_iE_i f)+L_i\Phi_i(D_i f)+M_i\Phi_i(E_i f)
                  +(d_i+L_iM_i)\Phi_i f\big).
 \end{gathered}
\]
The matrices \(W_i\), of order \(D\), have entries
\begin{equation}\label{eq:wi}
 (W_i)_{f,g}=[g]\,\mathcal U_i f\qquad(f,g\in\mathcal M).
\end{equation}
\begin{equation}\label{eq:dimension}
 D=\sum_{h=0}^N\binom{h+r-1}{r-1}\binom{h+s-1}{s-1}.
\end{equation}

Use the block matrix \(W\) of \eqref{eq:product-border} for the
matrices \(W_i\), and let \(e_0,e_n\) select the constant monomial
at the first and last cuts. The constant-to-constant entry of
the product of the \(W_i\) counts the perfect matchings.

\Needspace{5\baselineskip}
\begin{theorem}\label{thm:path}
The rank bounds imply the representation \eqref{eq:generators}.
With the matrices just defined,
\begin{equation}\label{eq:border}
 \per A=(W_1\cdots W_n)_{1,1}
       =\det\begin{pmatrix}
                I-W & e_n\\
                -e_0^{\mathsf T}&0
             \end{pmatrix}.
\end{equation}
The order of the determinant is \((n+1)D+1\). For fixed \(r,s\),
evaluation from the entries of \(A\) takes
\(O_{r,s}(n^2+n^{r+s+1})\) arithmetic operations, and \(O(n^2)\)
when \(r=s=1\).
\end{theorem}

\subsection{}

Now let \(\operatorname{char}\F=p>0\) in the generator
representation \eqref{eq:generators}. Use the same five terms
defining \(\mathcal U_i\), with reduction modulo the individual
\(p\)th powers in place of \(\Pi\). Denote this operator by
\(\mathcal U_i^{(p)}\), and set
\[
 (W_i^{(p)})_{f,g}=[g]\mathcal U_i^{(p)}f,
             \qquad f,g\in\mathcal B_{p,r,s}.
\]
The constant monomial is the first index. Form \(W^{(p)}\)
from these matrices as in \eqref{eq:product-border}.
The basis now has size independent of \(n\).

\begin{corollary}\label{cor:positive-ordered}
For the generators \eqref{eq:generators} in characteristic \(p\),
\begin{equation}\label{eq:positive-ordered}
 \per A=(W_1^{(p)}\cdots W_n^{(p)})_{1,1}
 =\det\begin{pmatrix}I-W^{(p)}&e_n\\-e_0^{\mathsf T}&0\end{pmatrix}.
\end{equation}
The determinant has order \((n+1)D_{p,r,s}+1\). From the supplied
generators the product can be evaluated in
\(p^{O(r+s)}n\) field operations.
\end{corollary}

\subsection{}

Write \(\mathcal C(\sigma)\) for the set of cycles of \(\sigma\),
including its fixed points. We regard \(\sigma\) as the directed
cycle cover with edges \(i\to\sigma(i)\), each of weight \(a_{i,\sigma(i)}\).
A permutation consisting of one cycle gives a directed Hamiltonian
cycle. Simultaneous permutation of rows and columns relabels the
vertices and preserves all cycle lengths.

Fix \(d\geq1\), and let \(\boldsymbol x=(x_1,\ldots,x_d)\) and
\(\boldsymbol z=(z_1,\ldots,z_d)\) be indeterminates. Put
\begin{equation}\label{eq:cycle-colors}
 F_A(\boldsymbol x,\boldsymbol z)
 =\sum_{\sigma\in S_n}\left(\prod_{i=1}^n a_{i,\sigma(i)}\right)
   \prod_{C\in\mathcal C(\sigma)}
       \left(\sum_{c=1}^d z_c x_c^{|C|}\right).
\end{equation}
The index \(c\) is a color for a cycle. The variable \(z_c\)
counts cycles of color \(c\), and \(x_c\) counts their vertices.

For one color put \(P_A(z)=F_A((1),(z))\) and \(C_k(A)=[z^k]P_A(z)\).
Thus \(C_k(A)\) is the sum of the weights of the permutations
with \(k\) cycles. In particular, \(C_1(A)\) is the weighted
Hamiltonian-cycle sum. The polynomial \(P_A(z)\) is the
\(\alpha\)-permanent with its parameter left indeterminate
\cite{VereJones}. The sign of a permutation with \(k\) cycles
is \((-1)^{n-k}\), so
\begin{equation}\label{eq:cycle-specializations}
 P_A(1)=\per A,\qquad P_A(-1)=(-1)^n\det A.
\end{equation}

The variables for the full polynomial record both ends of a
path. We use the same cut factorizations and couplings as for
the permanent.
At node \(I\),
use \(d\) matrices \(Z_I^{(c)}\) of \(r_I\times s_I\)
indeterminates. At a merge with children \(U,V\), put
\[
 B=\begin{pmatrix}0&B_{UV}\\B_{VU}&0\end{pmatrix},\qquad
 S_I=\begin{pmatrix}S_{U,I}\\ S_{V,I}\end{pmatrix},\qquad
 T_I=\begin{pmatrix}T_{U,I}\\ T_{V,I}\end{pmatrix}.
\]
For each \(c\), use a matrix \(Z^{(c)}\) with rows from the
coordinates of \(R_U,R_V\) and columns from those of \(C_U,C_V\).
Identify \(Z_U^{(c)}\) and \(Z_V^{(c)}\) with its diagonal blocks.
Write \(\partial_{\alpha\beta}^{(c)}=\partial/\partial Z_{\alpha\beta}^{(c)}\).
Define
\begin{equation}\label{eq:cycle-operator}
 \mathcal E_I=\sum_{c=1}^d\left(
 z_c\sum_{\alpha,\beta}B_{\alpha\beta}\partial_{\alpha\beta}^{(c)}
 +\sum_{\alpha,\beta,\gamma,\delta}
 B_{\alpha\beta}Z_{\gamma\delta}^{(c)}
       \partial_{\alpha\delta}^{(c)}\partial_{\gamma\beta}^{(c)}
 \right).
\end{equation}
Let \(\eta_I\) be the substitution \(Z^{(c)}\mapsto S_I Z_I^{(c)}T_I^{\mathsf T}\)
for all \(c\). Let \(\mathcal N_I\) be the monomials in the entries
of the \(Z_I^{(c)}\) of total degree at most \(|I|\). For
\(u\in\mathcal N_U\), \(v\in\mathcal N_V\), \(w\in\mathcal N_I\), set
\begin{equation}\label{eq:cycle-entries}
 t_I^{\mathrm{cyc}}(u,v;w)=[w]\,\eta_I\!\left(
         \sum_{\ell=0}^{|I|}\frac{\mathcal E_I^\ell}{\ell!}\,uv
                                          \right).
\end{equation}
At a leaf \(\{i\}\) put
\[
 l_i^{\mathrm{cyc}}(w)=[w]\sum_{c=1}^d x_c
 \left(z_c a_{ii}+R_{\{i\}}(i,:)Z_{\{i\}}^{(c)}C_{\{i\}}(i,:)^{\mathsf T}\right).
\]

For an ordering, put \(I_i=\{1,\ldots,i\}\), \(N=\lfloor n/2\rfloor\)
and \(D_d=\binom{N+drs}{drs}\). At each cut retain the monomials
of degree at most \(N\) and pad the coefficient vectors with
zeros to length \(D_d\). For \(i\geq2\) the transition from the
monomial \(u\) at \(I_{i-1}\) to \(w\) at \(I_i\) has entry
\[
 (W_i^{\mathrm{cyc}})_{u,w}
   =\sum_{v\in\mathcal N_{\{i\}}}
       l_i^{\mathrm{cyc}}(v)t_{I_i}^{\mathrm{cyc}}(u,v;w),
\]
with children \(I_{i-1},\{i\}\). At \(i=1\) the entry from the
constant state to \(w\) is \(l_1^{\mathrm{cyc}}(w)\). Form the
strictly upper block triangular matrix \(W^{\mathrm{cyc}}\)
with block \((i-1,i)\) equal to \(W_i^{\mathrm{cyc}}\), as in
\eqref{eq:border}. Let \(e_0,e_n\) select the constant monomials
at the first and last cuts.
The endpoint variables retain the pairing of the start and end
of each path. This allows cycles to be counted when they close.

\Needspace{5\baselineskip}
\begin{theorem}\label{thm:cycle-covers}
Assume that \(\F\) has characteristic zero. Suppose the two ranks at every node of \(\mathcal T\) are bounded
by \(r,s\), respectively. Then
\begin{equation}\label{eq:cycle-tree-formula}
 F_A(\boldsymbol x,\boldsymbol z)=
 \sum_{\substack{w_I\in\mathcal N_I\ \text{for all nodes }I\\
                 w_{\{1,\ldots,n\}}=1}}
 \prod_{i=1}^n l_i^{\mathrm{cyc}}(w_{\{i\}})
 \prod_{\substack{I\ \text{with}\\\text{children }U,V}}
    t_I^{\mathrm{cyc}}(w_U,w_V;w_I).
\end{equation}
For fixed \(d,r,s\), all coefficients can be computed in
polynomial time with polynomial bit complexity for rational
entries. The existence of a suitable tree suffices for
polynomial-time computation from \(A\) alone.
If the rank bounds hold along an ordering, then
\begin{equation}\label{eq:cycle-border}
 F_A(\boldsymbol x,\boldsymbol z)=
 \det\begin{pmatrix}
        I-W^{\mathrm{cyc}}&e_n\\
        -e_0^{\mathsf T}&0
      \end{pmatrix}.
\end{equation}
The order is \((n+1)\binom{\lfloor n/2\rfloor+drs}{drs}+1\).
\end{theorem}

\subsection{}

Fix \(1\le k\le n\). Choose marks \(v_1<\cdots<v_k\) in
\([n]\). Write \(\mathcal A(v_1,\ldots,v_k)\) for the ordered
partitions into \(S_1,T_1,\ldots,S_k,T_k\) with \(v_c\in S_c\)
and \(\min(S_c\cup T_c)=v_c\). The sets \(T_c\) may be empty.
The two sets of a pair will carry permanent and determinant
weights, respectively.

Put \(N=\lfloor n/2\rfloor\) and set
\begin{equation}\label{eq:fixed-cycle-dimension}
 D^+=\sum_{h=0}^N\binom{h+r-1}{r-1}\binom{h+s-1}{s-1},
 \qquad D^-=\binom{r+s}{r},\qquad
 D_k=(k+1)(D^+D^-)^k .
\end{equation}

Ordinary polynomial variables commute. In the formulas involving
determinants of color classes we also use exterior variables.
They anticommute and have square zero. A left derivative with
respect to an exterior variable removes its occurrence in an
ordered monomial with sign \((-1)^{j-1}\) if it is the \(j\)-th
factor, and gives zero if it is absent. Different colors use
separate tensor factors.

For the determinant in the ordered case, regard \(S_c,T_c\) as
the two colors of pair \(c\). At cut \(I\) use commuting variables
\(x_I^{(c)},y_I^{(c)}\) for \(S_c\), as in \eqref{eq:tree-entries}.
Use exterior variables \(u_I^{(c)},v_I^{(c)}\) for \(T_c\),
of respective lengths \(r_I,s_I\). At a leaf \(\{i\}\), the factors for the two
choices of color are
\begin{equation}\label{eq:fixed-cycle-leaf}
 \begin{split}
 L_i^{(c,+)}&=a_{ii}+
   (R_{\{i\}}(i,:)x_{\{i\}}^{(c)})(C_{\{i\}}(i,:)y_{\{i\}}^{(c)}),\\
 L_i^{(c,-)}&=-a_{ii}-
   (R_{\{i\}}(i,:)u_{\{i\}}^{(c)})(C_{\{i\}}(i,:)v_{\{i\}}^{(c)}).
 \end{split}
\end{equation}
At a merge use the operator \(\mathcal D_I\) of
\eqref{eq:tree-entries} on each commuting pair. On each exterior
pair use the operator
\begin{equation}\label{eq:fixed-cycle-exterior}
 \sum_{(J,K)\in\{(U,V),(V,U)\}}
     \sum_{\alpha,\beta}(B_{JK})_{\alpha\beta}
       \partial_{v_{K,\beta}^{(c)}}\partial_{u_{J,\alpha}^{(c)}} .
\end{equation}
Here the derivatives on exterior variables are left derivatives. Denote the sum of
these operators over all colors by \(\mathcal D_I^{(k)}\).
Apply \(\sum_{j=0}^{|I|}(\mathcal D_I^{(k)})^j/j!\) to the
product of the two child polynomials. Then apply the changes of
coordinates in \eqref{eq:tree-entries} to both types of variables.

Along the ordering \(1,\ldots,n\), also keep the number \(\ell\)
of marks chosen so far. At the next index one may use either
color of a pair \(c\leq\ell\). If \(\ell<k\), one may choose
the index as the next mark, using \(L_i^{(\ell+1,+)}\) and
increasing \(\ell\) by one. For each choice multiply by the
corresponding factor in \eqref{eq:fixed-cycle-leaf}, merge
with the next leaf, and keep commuting degrees at most \(N\)
in the row variables and in the column variables of each
pair. At the first index use the leaf factor alone. Denote by
\(W_i^{(k)}\) the matrix of this linear map in the monomial
bases, padded with zeros to order \(D_k\). Form \(W^{(k)}\)
with block \((i-1,i)\) equal to \(W_i^{(k)}\) and all other
blocks zero.

Let \(e_0,e_n\) select the constant states with \(\ell=0,k\)
at the two ends. The marked partitions count each permutation
with \(k\) cycles once, after cancellation of the other terms.

\Needspace{5\baselineskip}
\begin{theorem}\label{thm:fixed-cycles}
Assume that \(\F\) has characteristic zero. Suppose the two ranks at every node of \(\mathcal T\) are bounded
by \(r,s\), respectively. Then
\begin{equation}\label{eq:fixed-cycle-sum}
 C_k(A)=\sum_{v_1<\cdots<v_k}\ 
       \sum_{(S_1,T_1,\ldots,S_k,T_k)\in\mathcal A(v_1,\ldots,v_k)}
       \prod_{c=1}^k
       \per A[S_c,S_c]\,(-1)^{|T_c|}\det A[T_c,T_c].
\end{equation}
For fixed \(k,r,s\), the sum can be evaluated on \(\mathcal T\)
in polynomial time, with polynomial bit complexity for rational
entries. If the rank bounds hold
along an ordering, then
\begin{equation}\label{eq:fixed-cycle-border}
 C_k(A)=\det\begin{pmatrix}
              I-W^{(k)} &e_n\\
              -e_0^{\mathsf T}&0
             \end{pmatrix}.
\end{equation}
The order is at most \(1+(n+1)D_k=O_{k,r,s}(n^{k(r+s-1)+1})\).
For \(k=1\) this gives a determinant of order
\(O_{r,s}(n^{r+s})\) for the weighted Hamiltonian-cycle sum.
\end{theorem}

For \(k=1\) the mark in \eqref{eq:fixed-cycle-sum} is \(1\).
One can instead prescribe any vertex in \(S_1\) and drop
the minimum condition. For related permanent and determinant
formulas for Hamiltonian cycles, see Staples \cite{Staples}.

\subsection{}

For a partition \(\lambda\) of \(n\), let \(\chi_\lambda\) be the
corresponding irreducible character of \(S_n\). Its immanant is
\[
 \imm_\lambda(A)=\sum_{\sigma\in S_n}
                   \chi_\lambda(\sigma)\prod_i a_{i,\sigma(i)} .
\]
The ordinary hook partitions are \((n-k,1^k)\). More generally,
the \((p,q)\)-hook is the condition \(\lambda_{p+1}\le q\), with
missing parts taken to be zero. The Durfee size is the largest
\(j\) with \(\lambda_j\ge j\).

A power-sum specialization of \(F_A\) gives all immanants
whose partitions lie in a fixed \((p,q)\)-hook.

\Needspace{5\baselineskip}
\begin{corollary}\label{cor:hook-immanants}
Under the hypotheses of Theorem~\ref{thm:cycle-covers}, for fixed
nonnegative integers \(p,q\) with \(p+q>0\), all the immanants
\(\imm_\lambda(A)\) with \(\lambda_{p+1}\leq q\) can be computed
in polynomial time. In particular this holds for all ordinary
hook immanants and for all immanants of bounded Durfee size.
For rational entries the bit complexity is polynomial.
For ordinary hooks,
\begin{equation}\label{eq:ordinary-hooks}
 \frac{F_A((1,-t),(1,-1))}{1+t}
     =\sum_{k=0}^{n-1}\imm_{(n-k,1^k)}(A)t^k .
\end{equation}
\end{corollary}

\subsection{}

Determinant methods for counting Hamiltonian cycles on graphs of
bounded treewidth were developed by Bodlaender, Cygan, Kratsch and
Nederlof \cite{BCKN}. We compare the two calculations at a single
cut. For the use of
completion matrices in counting Hamiltonian cycles on graphs
of bounded pathwidth, see Curticapean, Lindzey and Nederlof
\cite{CLN}. For \(h\geq1\)
let \(\mathcal B_h\) be the set of nonnegative \(r\times s\)
integer arrays \(M=(m_{ab})\) of sum \(h\). List \(h\) pairs
\((a,b)\) with these multiplicities, and let \(U_M,V_M\)
have corresponding rows \(e_a^{\mathsf T},e_b^{\mathsf T}\)
of lengths \(r,s\), respectively. For \(M,M'\in\mathcal B_h\)
put
\[
 A_{M,M'}=\begin{pmatrix}
      0&U_M U_{M'}^{\mathsf T}\\
      V_{M'} V_M^{\mathsf T}&0
    \end{pmatrix},\qquad
 \mathcal K_h(z)=\big(P_{A_{M,M'}}(z)\big)_{M,M'\in\mathcal B_h}.
\]
The first \(h\) vertices precede the last \(h\). Put
\(B_h=\binom{h+r-1}{r-1}\binom{h+s-1}{s-1}\).
The following proposition compares the ranks for the whole
cycle polynomial and for its one-cycle coefficient.

\Needspace{5\baselineskip}
\begin{proposition}\label{prop:cycle-cut-rank}
Assume that \(\F\) has characteristic zero. The two ranks of each \(A_{M,M'}\) at every initial interval
cut are bounded by \(r,s\), respectively. Moreover,
\begin{equation}\label{eq:cycle-cut-ranks}
 \rk_{\F(z)}\mathcal K_h(z)=\binom{h+rs-1}{rs-1},\qquad
 B_h\leq \rk_{\F}\big([z]\mathcal K_h(z)\big)\leq D^- B_h .
\end{equation}
\end{proposition}

In particular, a factorization of the full completion matrix
into a left summary depending on \(M\) and a right summary
depending on \(M'\), with a fixed bilinear pairing over
\(\F(z)\), requires at least \(\binom{h+rs-1}{rs-1}\)
coordinates. This is the number of degree \(h\) monomials
in the endpoint variables of Theorem~\ref{thm:cycle-covers}
for one color. For the one-cycle coefficient the rank has
order \(h^{r+s-2}\) for fixed \(r,s\).

\section{Applications to structured matrices}\label{sec:structured}

\subsection{}

For this subsection index rows and columns by \(1,\ldots,n\).
Let \(R(X,Y)\) be a polynomial with fixed degrees in its two
variables. The Cauchy theorem applies to the Hilbert matrix after
adding \(R(i,j)\).

\begin{corollary}\label{cor:hilbert}
For \(B_n=(1/(i+j-1)+R(i,j))_{i,j=1}^n\), formula
\eqref{eq:cauchy-mainformula} applies with
\[
 x_i=i,\qquad y_j=1-j,\qquad P(X,Y)=R(X,1-Y).
\]
It can be evaluated in polynomial time for fixed degrees of \(R\).
\end{corollary}

For example, take \(R(i,j)=ij\). Then \(P(X,Y)=X-XY\), so
\(r=s=2\) and
\[
 H=\begin{pmatrix}0&0\\1&-1\end{pmatrix},\qquad
 B_2=\begin{pmatrix}2&5/2\\5/2&13/3\end{pmatrix}.
\]
Consequently \(\per B_2=179/12\); in order three the value is
\(397087/1080\). The same formula applies for every order.

The following rank bound explains why the Cauchy construction is
needed in addition to the bounded rank-width construction.

\begin{proposition}\label{prop:cauchy-width}
Let $A$ be as in \eqref{eq:cauchy-mainmatrix}, let $n\geq2$, and put
$q_0=\min(r,s)$. For every $I\subseteq[n]$,
\[
 w_A(I)\geq2\max\{0,\min(|I|,n-|I|)-q_0\}.
\]
Consequently, every rooted binary tree on $[n]$ has width at least
\[
 2\max\{0,\lceil n/3\rceil-q_0\}.
\]
\end{proposition}

\subsection{}

For a square matrix \(B\), its undirected support graph has edge
\(ij\), \(i\ne j\), if at least one of \(b_{ij},b_{ji}\) is nonzero.
A tree decomposition uses vertex bags on a tree. Each edge lies
in a bag, and the bags containing any fixed vertex form a nonempty
connected subtree. The two crossing ranks are unchanged by
inversion. Thus a decomposition for \(B\) also gives one for
\(B^{-1}\).

\Needspace{5\baselineskip}
\begin{corollary}\label{cor:inverse}
The following permanents can be computed in polynomial time and
have the finite formulas of Theorem~\ref{thm:tree} or
Theorem~\ref{thm:path}.
\begin{enumerate}
\item The permanent of \(B^{-1}\), where \(B\) is a nonsingular
  banded matrix of fixed bandwidth. For a tridiagonal \(B\), the
  calculation from \(B\) uses \(O(n^2)\) arithmetic operations.
\item The permanent of \(B^{-1}\), where \(B\) is nonsingular and a
  tree decomposition with bags of bounded size is given for its
  undirected support graph.
\end{enumerate}
\end{corollary}

The same rank bounds give polynomial-time cycle-polynomial
formulas for these inverses. This also applies to the rational
circulants in Corollary~\ref{cor:circulant} and the banded
block-constant matrices in Corollary~\ref{cor:block}.

\subsection{}

In the diagonal formulas below, letters for offsets and transfer
orders are local to each family; they no longer denote the rank
bounds of Section~\ref{sec:results}.
For scalar and block diagonal formulas we start indices at
zero. Write
\[
 P_n=(\ind_{j-i\equiv1\pmod n}),\qquad
 J_d=(\ind_{j-i=d}),\qquad
 \mathsf J_n=(\ind_{i+j=n-1}).
\]
Here \(\ind\) denotes the indicator of the stated condition.
The matrix \(P_n\) is the cyclic shift, \(J_d\) has ones on the
ordinary diagonal of offset \(d\), and right multiplication
by \(\mathsf J_n\) reverses columns. An open band of widths
\(p,q\) has offsets \(-p,\ldots,q\); a cyclic band starting at
zero has offsets \(0,\ldots,w\) modulo the order. In a block
band these offsets refer to block rows and columns. A
Toeplitz matrix has entries \(t_{j-i}\), and a Hankel matrix
has entries \(h_{i+j}\).

For \(n>w\), write \(C_n(c_0,\ldots,c_w)=\sum_{d=0}^w c_dP_n^d\).
The two composition polynomials used below are
\begin{align}
 \cH_m(q_1,\ldots,q_D)
  &=\sum_{\substack{\alpha\in\N^D\\ \sum_{j=1}^D j\alpha_j=m}}
       \frac{|\alpha|!}{\alpha_1!\cdots\alpha_D!}
       \prod_{j=1}^D q_j^{\alpha_j},\label{eq:hpoly}\\
 \cL_m(q_1,\ldots,q_D)
  &=\sum_{\substack{\alpha\in\N^D\\ \sum_{j=1}^D j\alpha_j=m}}
       \frac{m(|\alpha|-1)!}{\alpha_1!\cdots\alpha_D!}
       \prod_{j=1}^D q_j^{\alpha_j}\qquad(m\geq1).\label{eq:lpoly}
\end{align}
We put \(\cH_m=0\) for \(m<0\), so \(\cH_0=1\). If
\(\lambda_1,\ldots,\lambda_D\) are the roots of
\(X^D-\sum_{j=1}^D q_jX^{D-j}\), counted with multiplicity, then
\(\cL_m(q_1,\ldots,q_D)=\sum_{j=1}^D\lambda_j^m\).
The coefficients of \(\cL_m\) count tilings of a labelled circle.

We use the same transfer for scalar and block diagonals. Let
\(D_0,\ldots,D_w\) be \(b\times b\) matrices and let
\(T_{b,w,k}(D_0,\ldots,D_w)\) be indexed by the \(k\)-subsets
of \(\{0,\ldots,bw-1\}\). Write \(D=[D_0\ \cdots\ D_w]\).
For states \(U,V\) its entry is
\begin{equation}\label{eq:block-transition}
 \sum_S\per D[\{0,\ldots,b-1\},S],
\end{equation}
where \(S\) runs over the \(b\)-subsets of
\(\{0,\ldots,b(w+1)-1\}\setminus U\) such that
\[
 \{0,\ldots,b-1\}\subseteq U\cup S,\qquad
 V=\{j-b:j\in U\cup S,\ j\ge b\}.
\]
For \(b=1\), we write \(T_{w,k}\) for this matrix. In this case a choice is a single
displacement \(d\notin U\), of weight \(c_d\), and
\begin{equation}\label{eq:band-transition}
 V=\{j-1:j\in U\cup\{d\},\ j>0\},\qquad 0\in U\cup\{d\}.
\end{equation}

\subsection{}

For the signed formula, take cyclic offsets \(0,s,t,s+t\)
with all four residues distinct and \(\gcd(n,s,t)=1\).
The entries along each diagonal may vary with the row.
For three offsets \(0,s,t\) with the same gcd condition,
omit the last diagonal. The allowed nonzero steps correspond
to \((1,0),(0,1),(1,1)\) in the lattice construction
\[
 \Lambda=\{(x,y)\in\Z^2:sx+ty\equiv0\pmod n\}.
\]
Put \(g_0=\gcd(n,s)\), \(h_0=n/g_0\), and
\(u_0\equiv-(s/g_0)^{-1}t\pmod{h_0}\), \(0\le u_0<h_0\).
Write \(b_1=(h_0,0)\), \(b_2=(u_0,g_0)\).
For index \(i\), let \(\rho_i\) be the representative of \(i\)
under the map \((x,y)\mapsto sx+ty\) in the rectangle
\(0\le x<h_0,\ 0\le y<g_0\).
For an edge \(i\to j\) of step \(\delta\), write
\(\rho_i+\delta-\rho_j=m_1b_1+m_2b_2\).
For \(\epsilon,\eta\in\{0,1\}\), let \(K_{\epsilon\eta}\)
keep the diagonal of \(A\) and multiply the entry on this edge
by \(-(-1)^{\epsilon m_1+\eta m_2}\).

For constant weights write \(C=aI+bP_n^s+cP_n^t+dP_n^{s+t}\).
Given \(r,z,v\ge0\) with \(r+z+v\le n\), set
\[
 R=r+v,\qquad S=z+v,\qquad M=a^{n-r-z-v}b^rc^zd^v.
\]
When only \(0,s,t\) occur take \(v=0\). If \(n\mid sR+tS\), put
\(k=\gcd(R,S,(sR+tS)/n)\), with \(\gcd(0,0,0)=0\).

For the residue formula, write \(n=bm\) and let
\(C=\sum_{s\in S}a_sP_n^s\). Take integers \(u,h,w\) and
\(0\le d_s\le w\) such that \(s-u\equiv d_s h\pmod m\).
Set \(g=\gcd(m,h)\), \(L=m/g\) and require \(w<L\). Form the blocks
and transfers
\[
 D_d=\sum_{s:d_s=d}a_sP_b^{(s-u-dh)/m},\qquad
 T_k=T_{b,w,k}(D_0,\ldots,D_w).
\]
On states, \(\tau\) replaces each position \(b\ell+j\) by
\(b\ell+(j-h/g\bmod b)\), where \(0\le j<b\).
For the scalar band \(C_n(c_0,\ldots,c_w)\), put
\(\det(I-zT_{w,k})=1-\sum_{j=1}^{\binom wk}q_{k,j}z^j\).

For independent supported entries, let \(L(A)\) be the least
number of determinants in an identity
\begin{equation}\label{eq:signing-length}
 \per A=\sum_{\nu=1}^{L(A)}c_\nu\det(\varepsilon_\nu\circ A),
 \qquad \varepsilon_{\nu,ij}\in\{1,-1\},
\end{equation}
with constant coefficients \(c_\nu\) and constant signs, where
\(\circ\) is entrywise multiplication. The additive support
admits a four-determinant formula even with independent entries.
The last part of the theorem compares it with another
four-diagonal support.

\Needspace{5\baselineskip}
\begin{theorem}\label{thm:circulant}
With the above notation, the following hold.
\begin{enumerate}
\item For the additive support, or for its three-shift version,
\begin{equation}\label{eq:4det}
 \per A=\frac{-\det K_{00}+\det K_{10}+\det K_{01}+\det K_{11}}{2}.
\end{equation}
This uses \(O(n^3)\) arithmetic operations uniformly in \(s,t\),
with polynomial bit complexity for rational entries. A contribution to the
monomial \(M\) requires \(n\mid sR+tS\). Every permutation
contributing to \(M\) has \(k\) nontrivial cycles, all with the
same primitive displacement in \(\Lambda\). In particular,
\begin{equation}\label{eq:additive-sign}
 [M]\per C=(-1)^{r+z+v-k}[M]\det C .
\end{equation}

\item For the residue support,
\begin{equation}\label{eq:cyclic-residue}
 \per C=\left(\sum_{k=0}^{bw}\ \sum_{|U|=k}
                  (T_k^{\,L})_{U,\tau U}\right)^g .
\end{equation}
For fixed \(b,w\) this uses \(O_{b,w}(\log(n+1))\)
arithmetic operations on the weights. If \(|S|\) is fixed and
\(w\ge1\), the formula can still be evaluated in polynomial time
when \(bw=O(\log(n+1))\).

\item For \(n>w\), the consecutive cyclic band satisfies
\begin{equation}\label{eq:cyclic}
 \per C_n(c_0,\ldots,c_w)
  =\sum_{k=0}^w\tr(T_{w,k}^{\,n})
  =\sum_{k=0}^w\cL_n(q_{k,1},\ldots,q_{k,\binom wk}).
\end{equation}

\item For the connected additive four-shift support, \(L(A)=4\).
For the support \(\{0,1,2,4\}\) with \(n\ge6\),
\[
 L(A)\ge\left\lceil(8/3)^{\lfloor n/5\rfloor}\right\rceil.
\]
Its permanent is given by the row-dependent trace formula
\eqref{eq:block-cyclic} with \(b=1\) and \(w=4\), using \(O(n)\)
arithmetic operations.
\end{enumerate}
\end{theorem}

The residue formula includes the unweighted linear-jump case
with \(h=1,u=0\) and \(b,w\) fixed, as in \cite{GolinLeungWang}.
For example, let \(n=2m\), \(m\ge4\), and \(w=\lfloor\log_2 m\rfloor\).
For the shifts \(0,1,m,m+w\), take \(b=2\), \(u=0\) and \(h=1\).
The only nonzero blocks are
\[
 D_0=a_0I_2+a_mP_2,\qquad D_1=a_1I_2,\qquad D_w=a_{m+w}P_2.
\]
There are \(2^{2w}\le m^2\) states in \eqref{eq:cyclic-residue}.

For the five-diagonal formulas put
\begin{equation}\label{eq:penta-par}
 \begin{split}
 q_1&=c,\quad q_2=ae+bd,\quad q_3=ad^2+b^2e,\\
 q_4&=a^2e^2+abde,\quad q_5=-a^2ce^2,\quad q_6=-a^3e^3 .
 \end{split}
\end{equation}
For the positive formulas put \(x=ae,\ y=bd,\ h=ad^2+b^2e\), and define
\begin{equation}\label{eq:firstreturn}
 F(z)=cz(1+xz^2)+x^2z^4+
         \frac{yz^2(1+xz^2)+hz^3}{1-xz^2}.
\end{equation}

For matrices with at most three cyclic offsets, or with four
offsets of additive form, write \(\mathcal D_n(A)\) for the
following expression. Translate an offset to zero by shifting
columns. For one offset use the product of its entries. For
two offsets \(0,s\) use
\begin{equation}\label{eq:two-shift}
 \prod_O\left(\prod_{i\in O}a_{ii}+\prod_{i\in O}a_{i,i+s}\right),
\end{equation}
where \(O\) runs over the orbits of addition by \(s\). For
three or four offsets \(0,s,t\) or \(0,s,t,s+t\), use the
product of \eqref{eq:4det} on the classes modulo
\(\gcd(n,s,t)\), with the divided offsets in each class.

For constant weights in the connected additive case define
\(D_{\epsilon\eta}\) in an algebraic closure of \(\F\) by
\begin{equation}\label{eq:products}
 D_{\epsilon\eta}=
 \prod_{\substack{X^{h_0}=(-1)^\epsilon\\
                  X^{u_0}Y^{g_0}=(-1)^\eta}}
                 (a-bX-cY-dXY).
\end{equation}
In the three-shift case put \(d=0\). Put
\(f_\epsilon=X^{h_0}-(-1)^\epsilon\) and let \(R_{\epsilon\eta}\)
be the remainder of \((a-bX)^{g_0}-(-1)^\eta X^{-u_0}(c+dX)^{g_0}\)
modulo \(f_\epsilon\), using the inverse of \(X\) in this quotient.
Over \(\mathbb C\), set \(\zeta=\exp(\pi\mathrm i/n)\),
\(x_\epsilon=\zeta^{\epsilon g_0}\) and
\(y_{\epsilon\eta}=\zeta^{\eta h_0-\epsilon u_0}\).
The pairs in the product are
\((x_\epsilon\zeta^{2sj},y_{\epsilon\eta}\zeta^{2tj})\)
for \(0\le j<n\).

The corollary gives explicit low-width formulas, the product
over cyclic components, and the product and resultant forms
for constant weights.

\Needspace{8\baselineskip}
\begin{corollary}\label{cor:circulant}
\leavevmode
\begin{enumerate}
\item For \(n\ge3,4,5\), respectively,
\begin{align}
 \per C_n(a,b,c)&=a^n+c^n+
     \sum_{j=0}^{\lfloor n/2\rfloor}
        \frac{n}{n-j}\binom{n-j}{j}(ac)^j b^{n-2j},
                                                        \label{eq:three}\\
 \per C_n(a,b,c,d)&=a^n+d^n+
                \cL_n(b,ac,a^2d)+\cL_n(c,bd,ad^2).
                                                        \label{eq:four}
\end{align}
\begin{align}
 \per C_n(a,b,c,d,e)={}&a^n+e^n+\cL_n(b,ac,a^2d,a^3e)\notag\\
             &+\cL_n(q_1,\ldots,q_6)+\cL_n(d,ce,be^2,ae^3).
                                                        \label{eq:five}
\end{align}
The expanded four-diagonal formula is
\begin{equation}\label{eq:four-expanded}
 \per C_n(a,b,c,d)=a^n+d^n+
 \sum_{\substack{u,v,w\in\N\\u+2v+3w=n}}\frac{n(u+v+w-1)!}{u!v!w!}
 \left(a^{v+2w}b^uc^vd^w+a^wb^vc^ud^{v+2w}\right).
\end{equation}
The middle term for five diagonals has the positive form
\begin{equation}\label{eq:positive-cyclic}
 \cL_n(q_1,\ldots,q_6)=2\ind_{2\mid n}x^{\lfloor n/2\rfloor}
       +n\sum_{j=1}^n\frac{[z^n]F(z)^j}{j}.
\end{equation}

\item For \(g=\gcd(n,h)\) and \(n/g>w\),
\begin{equation}\label{eq:progression}
 \per\left(\sum_{d=0}^w c_dP_n^{u+dh}\right)
       =\bigl(\per C_{n/g}(c_0,\ldots,c_w)\bigr)^g .
\end{equation}
For \(n=2m\), \(m>2\) and \(\gcd(m,h)=1\),
\eqref{eq:cyclic-residue} for \(aI+bP_n^h+cP_n^{2h}+dP_n^m\)
uses \(D_0=aI_2+dP_2,\ D_1=bI_2,\ D_2=cI_2\). There are \(16\)
states, and \(\tau\) switches the two labels in a block
if \(h\) is odd.

\item In the cases defining \(\mathcal D_n(A)\),
\(\per A=\mathcal D_n(A)\). In particular, for constant
weights and \(g=\gcd(n,s,t)\), the component factors are equal:
\[
 \per C=\left[\mathcal D_{n/g}
   \bigl(aI+bP_{n/g}^{s/g}+cP_{n/g}^{t/g}
                         +dP_{n/g}^{(s+t)/g}\bigr)\right]^g.
\]

\item In \eqref{eq:4det}, for constant weights one may replace
\(\det K_{\epsilon\eta}\) by \(D_{\epsilon\eta}\). Also
\begin{equation}\label{eq:resultant}
 D_{\epsilon\eta}=\operatorname{Res}_X(f_\epsilon,R_{\epsilon\eta}),
\end{equation}
This uses \(O(h_0^2\log(g_0+1))\) operations in \(\F\).
For \(s=1,\ t=q,\ 2\le q\le n-2\), with
\(z_{\epsilon j}=\exp(\pi\mathrm i(2j+\epsilon)/n)\), the product is
\[
 D_{\epsilon\eta}=\prod_{j=0}^{n-1}
   (a-bz_{\epsilon j}-(-1)^{\epsilon+\eta}c z_{\epsilon j}^{q}
                     -(-1)^{\epsilon+\eta}d z_{\epsilon j}^{q+1}).
\]
One may use \((-1)^\eta\) in place of \((-1)^{\epsilon+\eta}\)
in all four products.

\item For Laurent polynomials \(f,g\) supported in fixed
intervals and \(g(P_n)\) nonsingular, the permanent of
\(f(P_n)g(P_n)^{-1}\) is the determinant \eqref{eq:border}
with bounded cut ranks. If \(P_n\) is replaced by \(P_n^h\),
use the formula in order \(n/\gcd(n,h)\) and raise it to the
power \(\gcd(n,h)\). In particular,
\begin{equation}\label{eq:eulerian}
 \per (I-\rho P_n)^{-1}
  =\frac{\displaystyle\sum_{k=0}^{n-1}
                  \eul{n}{k}\rho^{nk}}
         {(1-\rho^n)^n},\qquad \rho^n\ne1,
\end{equation}
Here \(\eul{n}{k}\) is the Eulerian number, which equals
\(\sum_{j=0}^{k+1}(-1)^j\binom{n+1}{j}(k+1-j)^n\).
\end{enumerate}
\end{corollary}

Apart from \(a^n,d^n\), the coefficient of a monomial
\(a^{\alpha_0}b^{\alpha_1}c^{\alpha_2}d^{\alpha_3}\) of degree \(n\)
in \eqref{eq:four}, with \(m=\alpha_1+2\alpha_2+3\alpha_3\), is
\[
 \begin{cases}
 \displaystyle\frac{n(\alpha_1+\alpha_2+\alpha_3-1)!}
                  {\alpha_1!\alpha_2!\alpha_3!},&m=n,\\[6pt]
 \displaystyle\frac{n(\alpha_0+\alpha_1+\alpha_2-1)!}
                  {\alpha_0!\alpha_1!\alpha_2!},&m=2n,\\[6pt]
 0,&\text{otherwise}.
 \end{cases}
\]
For example, \([a^4b^2d^2]\per C_8(a,b,c,d)=12\).
The six- and seven-diagonal sector orders are
\((1,5,10,10,5,1)\) and \((1,6,15,20,15,6,1)\), respectively.
In general the total number of states is \(2^w\) and the largest
sector has order \(\binom{w}{\lfloor w/2\rfloor}\).

\subsection{}

Let \(G\) be a finite abelian group, and let \(S\) be a nonempty
subset of \(G\). Write \(P_s\) for the translation matrix with a
one at \((i,i+s)\) and put \(C_S=\sum_{s\in S}a_sP_s\), where the
\(a_s\) are independent. Put \(H=\langle S-S\rangle\), \(m=|H|\)
and \(g=[G:H]\). Translating \(S\) into \(H\) and ordering by the
cosets of \(H\) gives \(g\) equal blocks. We write
\(\lvert\det C_S\rvert_{\rm coeff}\) for the polynomial obtained by
taking the absolute values of its integer coefficients.

Equality with the permanent holds exactly when all permutations
with the same shift multiplicities have the same sign.

\Needspace{5\baselineskip}
\begin{theorem}\label{thm:sign-classification}
The equality \(\per C_S=\lvert\det C_S\rvert_{\rm coeff}\)
holds in precisely the following cases.
\begin{enumerate}
\item \(|S|\leq2\).
\item \(|S|=3\) and one of the following holds.
\begin{itemize}
\item \(g=1\).
\item \(S\) is an arithmetic progression, with either \(m\) even or \(g=2\).
\item After translation and an isomorphism \(H\simeq\mathbb Z_7\),
      the support is \(\{0,1,3\}\).
\end{itemize}
\item \(|S|=4\), \(g=1\), and the four elements partition into
      two pairs with the same sum.
\end{enumerate}
\end{theorem}

In particular, for connected cyclic supports, the equality holds
for every support of size at most three and for exactly the additive
supports of size four. For \(C=aI_6+bP_6+cP_6^2+dP_6^4\), for
example,
\[
 [ab^2cd^2]\per C=12,\qquad [ab^2cd^2]\det C=0.
\]
The permutations \((0,2,4,1,5,3)\) and \((0,5,3,4,2,1)\), in
one-line notation with indices starting at zero, give opposite
signs at this monomial.

\subsection{}

Consider first a cyclic block matrix \(M\) with \(m\ge2\) block
rows and columns of size \(b\ge2\). Its only supported blocks
are \(A_i\) on the diagonal and \(B_i\) from \(i\) to \(i+1\)
modulo \(m\). Each supported block has two consecutive internal
cyclic offsets. For such a block \(D\), take its offsets to be
\(d,d+1\) with \(0\le d<b\). Write \(D^{(\eta)}\) for the matrix
that multiplies an entry of offset \(s\) in row \(j\) by \(\eta\)
if \(j+s\ge b\), leaving the other entries unchanged.
Put \(E=\operatorname{diag}((-1)^j)_{j=0}^{b-1}\).
For \(\epsilon=\pm1\) let \(F_\epsilon(x)\) be the determinant
of the block array obtained by the replacements
\[
 A_i\ \longmapsto\ A_i^{((-1)^b\epsilon)},\qquad
 B_i\ \longmapsto\ E B_i^{(\epsilon)}E,
\]
with the last forward block also multiplied by \((-1)^{m-1}x\).

For an open block band with offsets \(-p,\ldots,q\) put
\(w=p+q\) and \(B=\{0,\ldots,bp-1\}\). Use
\(T_i=T_{b,w,bp}\) on the blocks of block row \(i\) in order
of increasing offset, filling out-of-range blocks with zeros.
For a cyclic block band with offsets \(0,\ldots,w\), \(w<m\),
let \(T_{i,k}\) be the transfer on \(k\)-subsets from row \(i\).
For the block-constant case, take \(r\) row classes of sizes
\(m_i\) and \(s\) column classes of sizes \(n_j\), and let the
entries of block \(ij\) be \(b_{ij}\).
The theorem gives four determinants for the first class
and finite transfer and sum formulas for the other classes.

\Needspace{5\baselineskip}
\begin{theorem}\label{thm:block}
\leavevmode
\begin{enumerate}
\item For the cyclic block matrix \(M\),
\begin{equation}\label{eq:block-four-det}
 \per M=
 \frac{F_+(1)-F_+(-1)+F_-(1)+F_-(-1)}{2}.
\end{equation}
The block size and the number of blocks may both grow.

\item The open and cyclic block bands, respectively, satisfy
\begin{equation}\label{eq:block-open}
 \per A=(T_0\cdots T_{m-1})_{B,B}.
\end{equation}
\begin{equation}\label{eq:block-cyclic}
 \per A=\sum_{k=0}^{bw}\tr(T_{0,k}\cdots T_{m-1,k}).
\end{equation}

\item For the block-constant matrix \(A\),
\begin{equation}\label{eq:constant-blocks}
 \per A=\left(\prod_i m_i!\right)\left(\prod_j n_j!\right)
 \sum_{\substack{k_{ij}\geq0\\\sum_j k_{ij}=m_i,\ \sum_i k_{ij}=n_j}}
             \prod_{i,j}\frac{b_{ij}^{k_{ij}}}{k_{ij}!}.
\end{equation}
For fixed \(r,s\), at most \((n+1)^{(r-1)(s-1)}\) arrays
occur in the sum.
\end{enumerate}
\end{theorem}

For a block Toeplitz band the transfer is the same matrix
\(T\) at every row, using the untruncated block entries.
Let \(T'\) omit state \(B\). A block circulant band has common
sector matrices \(T_k\). The corollary gives the resulting
coefficient and trace formulas.

\Needspace{5\baselineskip}
\begin{corollary}\label{cor:block}
For a block Toeplitz band,
\begin{equation}\label{eq:block-toeplitz}
 \per A=[z^m]\frac{\det(I-zT')}{\det(I-zT)}.
\end{equation}
Its order is \(\binom{b(p+q)}{bp}\), which is polynomial in \(m\)
for fixed \(b,p\), even when \(q\) is as large as \(m-1\). For a block
circulant band with common transfers \(T_k\), the permanent is
\(\sum_k\tr(T_k^m)\). The orders and the number of transfers are
fixed when \(b,w\) are fixed. A cyclic shift of the block columns
allows any common starting offset. Reversing the block columns
of a block Hankel matrix gives the block Toeplitz formula.

For a block circulant with just two cyclic block offsets,
translate one offset to zero and order the block indices along
the cycles of the remaining offset. When its supported blocks
satisfy Theorem~\ref{thm:block}(1), use \eqref{eq:block-four-det}
on each cycle of length at least two and multiply the results.

For a block-constant matrix with a fixed band of nonzero blocks
on a common ordered partition, \eqref{eq:border} gives a
determinant of polynomial order, even when the number and the
sizes of the blocks vary.
\end{corollary}

\subsection{}

For an open band with widths \(p,q\), let \(w=p+q\) and
\(B=\{0,\ldots,p-1\}\). Form \(T_i=T_{w,p}\) with weights
\(c_d=a_{i,i+d-p}\), using zero for out-of-range entries.
For \(A_n=\sum_{d=-p}^q t_dJ_d\), use \(T=T_{w,p}\) with the
constant weights \(t_{d-p}\), of order \(D=\binom{p+q}{p}\).
Let \(T'\) omit state \(B\) and write
\[
 \det(I-zT)=1-\sum_{j=1}^Dq_jz^j,\qquad
 \det(I-zT')=\sum_{j=0}^{D-1}r_jz^j .
\]

For a Toeplitz sequence with rational diagonal series, write
\begin{equation}\label{eq:toeplitz-rational}
 \sum_{\ell\ge1}t_\ell z^{\ell-1},\qquad
 \sum_{\ell\ge1}t_{-\ell}z^{\ell-1}
\end{equation}
If the numerator and denominator degrees are bounded, we use
constant generators \(t_\ell=p^{\mathsf T}B^{\ell-1}q\) and
\(t_{-\ell}=u^{\mathsf T}C^{\ell-1}v\) for \(\ell\ge1\).
For each order \(n\), let \(T_n\) be the common matrix \(W_i\)
of \eqref{eq:wi}, with \(d_i=t_0\) and the degree bound
\(\lfloor n/2\rfloor\). Let \(T_n'\) omit the constant
monomial. The matrix \(T_n\) may grow with \(n\); the fixed-band
matrix \(T\) above does not.

In characteristic \(p>0\), use instead the common matrix
\(T_p=W_i^{(p)}\) of Corollary~\ref{cor:positive-ordered}, with
the same constant generators. Its order is \(D_{p,r,s}\), and
it is independent of \(n\). Write \(T_p'\) for the matrix with
the constant row and column deleted. The \(p\) in this notation
is the characteristic, unrelated to an open bandwidth.

There is also a characteristic-zero family not requiring rational
diagonal series. For a univariate polynomial \(P\), use the
one-based indices \(1,\ldots,n\) and put
\[
 C_n^{\mathrm T}=\left(\frac1{\alpha+j-i}+P(j-i)\right)_{i,j=1}^n,
\]
assuming that every denominator is nonzero. In the Cauchy formula
take \(x_i=-i\), \(y_j=-\alpha-j\) and the bivariate polynomial
\(P(X-Y-\alpha)\). A degree bound \(d\) for \(P\) permits
\(r=s=d+1\).

For the padding formula, let \(S\) be the ordinary offsets
of a matrix \(A\) of order \(n\). Choose \(N\ge n\) with \(N=n\)
or \(0\in S\), so that \(S\) modulo \(N\) has at most three
elements or is an additive four-set. The cyclic formula
\(\mathcal D_N\) is given in Corollary~\ref{cor:circulant}.
For the direct signing, take offsets \(0,p,q,p+q\) with
\(p>0>q\) and \(p+q\ne0\). Let \(g=\gcd(p,q)\), choose
\(x,y\in\{1,-1\}\) with \(x^{q/g}y^{-p/g}=-1\), and let
\(K\) be obtained by multiplying the entries on \(p,q,p+q\)
by \(-x,-y,-xy\), respectively. The next theorem collects the open-diagonal and dense
Toeplitz formulas.

\Needspace{5\baselineskip}
\begin{theorem}\label{thm:toeplitz}
With the above notation, the following hold.
\begin{enumerate}
\item The open band and the corresponding Toeplitz band satisfy
\begin{equation}\label{eq:open-product}
 \per A=(T_0\cdots T_{n-1})_{B,B}.
\end{equation}
\begin{equation}\label{eq:open}
 \per A_n=(T^n)_{B,B}
    =[z^n]\frac{\det(I-zT')}{\det(I-zT)}
    =\sum_{j=0}^{D-1}r_j\cH_{n-j}(q_1,\ldots,q_D).
\end{equation}
For fixed diagonal weights the generating function is
\begin{equation}\label{eq:open-series}
 \sum_{n\ge0}\per A_n z^n=\frac{\det(I-zT')}{\det(I-zT)}.
\end{equation}
At fixed \(w\) the power takes \(O(\log(n+1))\) arithmetic
operations. At fixed \(p\) the order is polynomial in \(n\)
for \(q\le n-1\), and \eqref{eq:open-product} takes
\(O_p(n^{p+1})\) operations. Transposition gives the
corresponding bound with the upper bandwidth fixed. For growing
\(q\), the stated bound uses \eqref{eq:open-product}, without
expanding \(\cH_{n-j}\) into its individual summands.

\item For the rational diagonal series,
\begin{equation}\label{eq:rationaltoeplitz}
 \per A_n=[z^n]\frac{\det(I-zT_n')}{\det(I-zT_n)} .
\end{equation}
The determinant orders are polynomial in \(n\) when the degrees
in \eqref{eq:toeplitz-rational} are bounded. This is a coefficient
formula for the specified \(n\); unlike \eqref{eq:open-series},
its rational function may depend on \(n\).

\item With the stated padding,
\begin{equation}\label{eq:open-padding}
 \per A=\mathcal D_N(A\oplus I_{N-n}),
\end{equation}
Any three ordinary diagonals are covered with \(N=n\).
For four offsets \(d_1<d_2<d_3<d_4\) that include zero,
one may take \(N=d_3+d_4-d_1-d_2\) if \(N\ge n\).
Both applications use \(O(n^3)\) arithmetic operations.

\item For the additive ordinary offsets,
\begin{equation}\label{eq:additive-open}
 \per A=\det K.
\end{equation}
The cost is \(O(n^3)\) arithmetic operations uniformly
in \(p,q\).

\item In characteristic \(p>0\), rational diagonal series give
\begin{equation}\label{eq:positive-toeplitz}
 \sum_{n\geq0}\per A_n z^n
        =\frac{\det(I-zT_p')}{\det(I-zT_p)}.
\end{equation}
For fixed \(p,r,s\), a specified value \(\per A_n=(T_p^n)_{1,1}\)
can be evaluated in \(O(\log(n+1))\) field operations after
\(T_p\) has been formed.

\item In characteristic zero, \(\per C_n^{\mathrm T}\) is given
by \eqref{eq:cauchy-mainformula} with the nodes and polynomial
specified above. For fixed degree of \(P\), this is a
polynomial-time determinant formula, also in bit complexity for
rational input.
\end{enumerate}
\end{theorem}

For row weights with period \(h\) before truncation, put
\(M=T_0\cdots T_{h-1}\), \(N_r=T_0\cdots T_{r-1}\) and \(N_0=I\).
For \(n=hk+r\), \(0\le r<h\), write
\[
 \det(I-zM)=1-\sum_{j=1}^Dq_jz^j,\qquad
 e_B^{\mathsf T}\operatorname{adj}(I-zM)N_re_B
       =\sum_{j=0}^{D-1}r_{r,j}z^j ,
\]
where \(e_B\) is the vector of state \(B\).

For \(A_n=aJ_{-2}+bJ_{-1}+cI+dJ_1+eJ_2\), use \(q_j,x,y,h,F\)
from \eqref{eq:penta-par} and \eqref{eq:firstreturn}. Put
\(B(0,0)=1\), \(B(0,v)=0\) for \(v>0\), and
\(B(s,v)=\binom{s+v-1}{v}\) for \(s>0\).
In the symmetric case \(a=e=v,\ b=d=u,\ c=t\), put
\[
 (s_1,\ldots,s_5)
 =(t+v,\ u^2-tv,\ vu^2+tv^2,\ v^4-tv^3,\ -v^5).
\]

For arbitrary tridiagonal entries \(d_i,u_i,l_i\) on the
diagonal, superdiagonal and subdiagonal, put \(x_i=u_il_i\).
Let \(E\) range over the nonconsecutive subsets of
\(\{1,\ldots,n-1\}\), and let \(V(E)=\bigcup_{i\in E}\{i,i+1\}\).
For the separated diagonals \(-r,0,r\) set
\(\ell_s=\lfloor(n-s)/r\rfloor+1\) for \(1\le s\le\min(r,n)\).
For a centrosymmetric tridiagonal matrix, let \(p_j\) be the
sum in \eqref{eq:tridiag} restricted to its first \(j\) indices,
with \(p_0=1\).

For geometric diagonals take \(t_0=c\), \(t_\ell=ap^{\ell-1}\),
and \(t_{-\ell}=bq^{\ell-1}\) for \(\ell\ge1\).
Put \(\tau=pq\) and \(N=\lfloor n/2\rfloor\). The matrix \(T_n\)
in \eqref{eq:rationaltoeplitz} is indexed by \(0,\ldots,N\) and
has the tridiagonal entries
\begin{equation}\label{eq:geometric}
 \begin{split}
 (T_n)_{0,0}&=c,\qquad
 (T_n)_{h,h}=c\tau^h+h(aq+bp)\tau^{h-1}\quad (h\geq1),\\
 (T_n)_{h,h+1}&=ab\tau^h,\qquad
 (T_n)_{h+1,h}=(h+1)^2\tau^h \quad (0\leq h<N).
 \end{split}
\end{equation}
Write \(T_n(z)\) for this matrix with \(c\) replaced by \(zc\)
and the factor \((h+1)^2\) by \((h+1)(h+z)\), and let \(T_n'(z)\)
omit index zero.

In characteristic three, the geometric family
\(A_n=(t^{|i-j|})\) has basis \(1,xy,x^2y^2\). Put
\(\tau=t^2\). Its fixed matrix is
\begin{equation}\label{eq:positive-geometric-matrix}
 T_3=\begin{pmatrix}
       1&\tau&0\\1&0&\tau^2\\0&\tau&-\tau^2
     \end{pmatrix}.
\end{equation}

For periodic bands the product over one period is fixed.
The remaining formulas are the Hessenberg, pentadiagonal,
tridiagonal and geometric specializations.

\Needspace{5\baselineskip}
\begin{corollary}\label{cor:toeplitz}
\leavevmode
\begin{enumerate}
\item For row period \(h\),
\begin{equation}\label{eq:periodic}
 \per A=[z^k]\,e_B^{\mathsf T}(I-zM)^{-1}N_r e_B
       =\sum_{j=0}^{D-1}r_{r,j}\cH_{k-j}(q_1,\ldots,q_D).
\end{equation}

\item For \(A_n=aJ_{-1}+bI+\sum_{j=1}^{r-1}c_jJ_j\),
\begin{equation}\label{eq:hessenberg}
 \per A_n=\cH_n(b,ac_1,\ldots,a^{r-1}c_{r-1}).
\end{equation}
At offsets \(-1,0,1,2\) with weights \(a,b,c,d\), this is
\begin{equation}\label{eq:four-open-expanded}
 \boxed{\displaystyle
 \per T=\sum_{\substack{u,v,w\in\N\\u+2v+3w=n}}\frac{(u+v+w)!}{u!v!w!}
          b^u(ac)^v(a^2d)^w.}
\end{equation}

\item For the pentadiagonal case,
\begin{equation}\label{eq:penta-open}
 \per A_n=\cH_n(q_1,\ldots,q_6)-ae\,\cH_{n-2}(q_1,\ldots,q_6).
\end{equation}
The corresponding generating function is
\begin{equation}\label{eq:penta-gf}
 \sum_{n\ge0}\per A_n z^n=\frac{1-aez^2}{1-\sum_{j=1}^6 q_jz^j}.
\end{equation}
The positive coefficient and sum formulas are
\begin{equation}\label{eq:positive-five}
 \per A_n=[z^n]\sum_{j=0}^n F(z)^j.
\end{equation}
\begin{equation}\label{eq:positive-expanded}
 \per A_n=
 \sum_{\substack{i,j,k,\ell,u,v\in\N\\
                 i+4j+2k+3\ell+2u+2v=n\\u\leq i+k}}
  \frac{(i+j+k+\ell)!}{i!j!k!\ell!}\binom{i+k}{u}
  B(k+\ell,v)c^i y^k h^\ell x^{2j+u+v}.
\end{equation}
Every term in \eqref{eq:positive-expanded} is nonnegative when
the entries are nonnegative. For the symmetric specialization,
\begin{equation}\label{eq:symmetric}
 \per A_n=\cH_n(s_1,\ldots,s_5)-v\cH_{n-1}(s_1,\ldots,s_5),
\end{equation}

\item For the tridiagonal case and the separated-diagonal case,
respectively,
\begin{equation}\label{eq:tridiag}
 \per A=\sum_E\prod_{i\in E}x_i\prod_{j\notin V(E)}d_j.
\end{equation}
\begin{equation}\label{eq:tridiag-stride}
 \per(dI+uJ_r+lJ_{-r})=
 \prod_{s=1}^{\min(r,n)}
   \sum_{k=0}^{\lfloor\ell_s/2\rfloor}
       \binom{\ell_s-k}{k}d^{\ell_s-2k}(ul)^k .
\end{equation}
For constant \(d_i=d,\ x_i=x\), the first sum is
\(\sum_{k=0}^{\lfloor n/2\rfloor}\binom{n-k}{k}d^{n-2k}x^k\).
For a centrosymmetric tridiagonal matrix, the permanent is
\[
 p_m^2+x_m p_{m-1}^2\quad(n=2m),\qquad
 d_{m+1}p_m^2+2x_m p_m p_{m-1}\quad(n=2m+1,\ m\ge1).
\]

\item For the geometric diagonals the permanent is
\eqref{eq:rationaltoeplitz} with \eqref{eq:geometric}, and
\begin{equation}\label{eq:geometric-cycles}
 P_{A_n}(z)=[t^n]\frac{\det(I-tT_n'(z))}{\det(I-tT_n(z))},
\end{equation}
The determinant orders are at most \(\lfloor n/2\rfloor+1\).
If \(M(n)\) is the cost of multiplying polynomials of degree
at most \(n\), with \(M(n)/n\) nondecreasing, computing the
permanent or evaluating \(P_{A_n}\) at a specified scalar takes
\(O(M(n)\log(n+1))\) arithmetic operations. All coefficients of
\(P_{A_n}\) can be computed in \(O(nM(n)\log(n+1))\) operations.

\item Over any field of characteristic three, for
\(A_n=(t^{|i-j|})\) and \(\tau=t^2\),
\begin{equation}\label{eq:positive-geometric}
 \sum_{n\geq0}\per A_n z^n
 =\frac{1+\tau^2z-\tau^3z^2}
        {1+(\tau^2-1)z-(\tau+\tau^2+\tau^3)z^2}.
\end{equation}
\end{enumerate}
\end{corollary}

The bound in Theorem~\ref{thm:toeplitz}(5) counts field
operations. Over \(\F_p(t)\), the degrees of rational functions
also enter the bit complexity. At \(t=1\),
\eqref{eq:positive-geometric} becomes \(1+z-z^2\), agreeing with
\(\per J_n=n!\) in characteristic three.

For \(aJ_{-2}+bJ_{-1}+cI+dJ_u+eJ_v\), with
\(0<u<v<n\), take \(p=2,\ q=v\) in \eqref{eq:open} and zero
weights on the other diagonals. For the sign formula in
Theorem~\ref{thm:toeplitz}(4), one
example is
\[
 \per(aI+bJ_2+cJ_{-3}+dJ_{-1})
   =\det(aI+bJ_2-cJ_{-3}+dJ_{-1}).
\]
For \eqref{eq:open-padding}, in order \(50\) the offsets
\(-21,0,13,24\) give \(N=58\).

For example, the matrix \((\rho^{|i-j|})_{i,j=1}^n\) is obtained by
taking \(c=1\) and \(a=b=p=q=\rho\). Its determinant is
\((1-\rho^2)^{n-1}\), so this is a dense full-rank family for
\(\rho\ne0\) and \(\rho^2\ne1\). The two ranks across each nontrivial
initial interval cut are one. In order \(2m\), the odd-row, even-column
minor has determinant \(\rho^m(1-\rho^2)^{m-1}\), so the rank
is at least \(m\) after any diagonal change. For rational
\(0<\rho<1\), the weighted clique-width is at least \(n/6\).
For \(\rho=2\), \eqref{eq:geometric-cycles} gives
\[
 P_{A_4}(z)=z^4+108z^3+848z^2+768z .
\]

In \eqref{eq:positive-cyclic}, the second term is obtained from
\eqref{eq:positive-expanded} by multiplying each summand by
\(n/(i+j+k+\ell)\) for \(n\ge1\).
The pentadiagonal denominator and its symmetric reduction also
occur in \cite{AK1,AK2}. In the symmetric specialization both
\(A_n\) and \(A_n\mathsf J_n\) are symmetric and centrosymmetric.
At \(t=u=v=1\) the generating function is
\((1-z)/(1-2z-2z^3+z^5)\). Its coefficients for \(n=0,\ldots,10\)
are \(1,1,2,6,14,31,73,172,400,932,2177\), the perfect-matching
counts for biadjacency matrix \((\ind_{|i-j|\le2})\).
For related cycle continued fractions see Deb \cite{Deb}.

More generally, when every offset is a multiple of \(r\),
the permanent is the product over the residue classes modulo
\(r\) of the ordinary diagonal formulas with offsets divided
by \(r\). In the direct signing one may take \((x,y)=(-1,1)\)
if \(q/g\) is odd and \((1,-1)\) otherwise. For constant weights,
the formula reads
\[
 \per(aI+bJ_p+cJ_q+dJ_{p+q})
 =\det(aI-bxJ_p-cyJ_q-dxyJ_{p+q}).
\]

\subsection{}

Let \(H\) have order \(n\) and put \(A=H\mathsf J_n\).
The anti-diagonal \(i+j=n-1-d\) of \(H\) is the ordinary
diagonal of offset \(d\) in \(A\). Column reversal preserves
the permanent, so the open-diagonal formulas apply to \(H\).

\Needspace{5\baselineskip}
For the characteristic-zero reciprocal family, use one-based
indices and put
\[
 C_n^{\mathrm H}=\left(\frac1{\alpha+i+j}+P(i+j)\right)_{i,j=1}^n,
\]
where the denominators are nonzero and \(P\) is a polynomial.
In \eqref{eq:cauchy-mainformula}, take \(x_i=i\),
\(y_j=-\alpha-j\) and the bivariate polynomial
\(P(X-Y-\alpha)\).

\begin{theorem}\label{thm:hankel}
We have \(\per H=\per A\), with the following specializations.
\begin{enumerate}
\item For anti-diagonals with \(-p\le d\le q\), use
\eqref{eq:open-product} for arbitrary entries. If \(H=(h_{i+j})\),
take \(t_d=h_{n-1-d}\) in \eqref{eq:open}. The same fixed one-sided
bound applies. For rational diagonal series of this \(t_d\),
the permanent is \eqref{eq:rationaltoeplitz}.
\item For any three anti-diagonals, or for four satisfying the
condition of Theorem~\ref{thm:toeplitz}(3) in the offsets \(d\),
\[
 \per H=\mathcal D_N\bigl((H\mathsf J_n)\oplus I_{N-n}\bigr).
\]
For offsets satisfying Theorem~\ref{thm:toeplitz}(4), let \(K\)
be its signed matrix for \(A\). Then
\[
 \per H=\det K=(-1)^{n(n-1)/2}\det(K\mathsf J_n).
\]

\item In characteristic zero, \(\per C_n^{\mathrm H}\) is given
by \eqref{eq:cauchy-mainformula} with the specified data. For
fixed degree of \(P\), the arithmetic cost is polynomial, and
so is the bit complexity for rational input.

\item In characteristic \(p>0\), if the reversed matrix \(A\)
has generators of orders \(r,s\), then its permanent and that
of \(H\) are given by \eqref{eq:positive-ordered}, with cost
\(p^{O(r+s)}n\) from the generators. For a reversed Toeplitz
matrix with constant generators, this is \((T_p^n)_{1,1}\).
\end{enumerate}
\end{theorem}

Column reversal gives the separated tridiagonal, pentadiagonal
and geometric parameters of Corollary~\ref{cor:toeplitz}.

\Needspace{5\baselineskip}
\begin{corollary}\label{cor:hankel}
For \(H=(dI+uJ_r+lJ_{-r})\mathsf J_n\), its permanent is the
product in \eqref{eq:tridiag-stride}. For
\(H=(aJ_{-2}+bJ_{-1}+cI+dJ_1+eJ_2)\mathsf J_n\), it is
\[
 \cH_n(q_1,\ldots,q_6)-ae\,\cH_{n-2}(q_1,\ldots,q_6),
\]
with \(q_j\) from \eqref{eq:penta-par}. Its positive form is
\eqref{eq:positive-expanded}. For \(H_{ij}=t_{n-1-i-j}\) with
\(t_0=c,\ t_\ell=ap^{\ell-1},\ t_{-\ell}=bq^{\ell-1}\)
for \(\ell\ge1\), use \eqref{eq:rationaltoeplitz} with the
tridiagonal matrix \eqref{eq:geometric}.
\end{corollary}

For a fixed sequence \(h_j\), the reversed diagonal sequence
\(t_d=h_{n-1-d}\) can change with \(n\). The last part of the
theorem is an evaluation for each order; it does not assert a
single rational generating function for all such Hankel matrices.

\subsection{}

Let \(B\) be an \(m\times m\) \((0,1)\) matrix, and put
\[
 S=\begin{pmatrix}0&B\\ B^{\mathsf T}&0\end{pmatrix},
 \qquad C=B\oplus(\mathsf J_mB\mathsf J_m)
\]
and \(D=S\oplus(\mathsf J_{2m}S\mathsf J_{2m})\).
The proposition relates their permanents to the permanent of \(B\).

\Needspace{5\baselineskip}
\begin{proposition}\label{prop:symmetry}
The matrix \(S\) is symmetric and \(C\) is centrosymmetric,
with \(\per S=\per C=(\per B)^2\). The matrix \(D\) has
both symmetries and \(\per D=(\per B)^4\).
\end{proposition}

Thus polynomial-time evaluation in either whole class would give
polynomial-time evaluation of arbitrary \((0,1)\) permanents.
For the centrosymmetric matrix with rows \((a,b),(b,a)\), the
permanent is \(a^2+b^2\), while the product from the two
one-dimensional invariant subspaces is \((a+b)(a-b)\).

\section{Proofs}

\subsection{Proof of Theorem~\ref{thm:tree}}\label{sec:proof-tree}

We first prove the formula for a given tree. Take two copies of \(I\), one for rows and one for columns. For a
matching \(M\) between the two copies write its contribution as
\begin{equation}\label{eq:partial}
 \prod_{(i,j)\in M}a_{ij}
 \prod_{\substack{i\in I\\i\text{ unmatched as a row}}} R_I(i,:)x_I
 \prod_{\substack{j\in I\\j\text{ unmatched as a column}}} C_I(j,:)y_I.
\end{equation}
Let \(\mathcal P_I\) be the sum of these contributions. A matching of
\(a\) edges leaves \(|I|-a\) factors of each kind. Thus its monomials
belong to \(\mathcal M_I\). At a leaf \(\{i\}\) the empty matching gives
the two coordinate factors and the one-edge matching gives \(a_{ii}\).
This is the polynomial defining \(l_i\).

Let \(U,V\) be the children of \(I\). Take a matching from each child,
and choose an unmatched row \(i\) in \(U\) and an unmatched column \(j\)
in \(V\). When the first sum in \(\mathcal D_I\) differentiates their
factors, the scalar obtained is
\[
 \sum_{\alpha,\beta}R_U(i,\alpha)(B_{UV})_{\alpha\beta}
                              C_V(j,\beta)
       =\big(R_UB_{UV}C_V^{\mathsf T}\big)_{ij}=a_{ij}.
\]
The other factors are untouched. The second sum gives \(a_{ij}\) for a
row in \(V\) and a column in \(U\). If several endpoint factors contain
the same variable, the product rule still makes one choice of a factor at
a time. It therefore counts choices of endpoints with their
multiplicities. For example, with one row variable and one column
variable,
\[
 \sum_{\ell=0}^{2}\frac{(b\partial_x\partial_y)^\ell}{\ell!}
       x^2y^2=x^2y^2+4bxy+2b^2.
\]
There are four choices for one pair of factors and two bijections
when both row factors are paired. The coefficients remain integers
when the variables are replaced by longer lists.

For a fixed set of \(\ell\) disjoint cross edges, the repeated
differentiation can add them in any of \(\ell!\) orders. Conversely, once
a row or column factor has been removed, it is unavailable in the next
differentiation. Thus every term after \(\ell\) applications comes from
\(\ell\) disjoint cross edges. In characteristic zero the factor
\(1/\ell!\) counts each set once. There are at most \(|I|\) cross
edges in a matching.

The same count can be used over any field. Write \(X=(x_U,x_V)\)
and \(Y=(y_U,y_V)\) for the lists of variables, of lengths
\(\rho=r_U+r_V\) and \(\sigma=s_U+s_V\). Let \(B=(b_{ab})\) be the
block matrix with zero diagonal blocks and off-diagonal blocks
\(B_{UV},B_{VU}\). For a monomial \(X^\alpha Y^\beta\), take a
\(\rho\times\sigma\) array \(M=(m_{ab})\) of nonnegative integers.
Write \(p_a=\sum_b m_{ab}\) and \(q_b=\sum_a m_{ab}\). The contraction
inside \(\rho_I\) in \eqref{eq:tree-entries} replaces the monomial by
\begin{equation}\label{eq:integer-merge}
\begin{split}
 \sum_{\substack{M\ge0\\p\le\alpha,\ q\le\beta}}
 &\left(\prod_a\binom{\alpha_a}{\alpha_a-p_a,m_{a1},\ldots,m_{a\sigma}}\right)
  \left(\prod_b\binom{\beta_b}{\beta_b-q_b,m_{1b},\ldots,m_{\rho b}}\right)\\
 &{}\cdot\left(\prod_{a,b}m_{ab}!\,b_{ab}^{m_{ab}}\right)
     X^{\alpha-p}Y^{\beta-q}.
\end{split}
\end{equation}
Here the multinomial coefficients and factorials are integers
read in the field. To obtain the term for \(M\), first distinguish
the repeated factors in the monomial. Of the \(\alpha_a\) factors
of type \(a\), choose \(m_{ab}\) to be paired with factors of
type \(b\), for each \(b\). The first multinomial counts these
choices. The second multinomial chooses the factors on the other
side. For each \(a,b\), pair the selected factors in \(m_{ab}!\)
ways and include \(b_{ab}\) for each pair. This counts every
choice of disjoint edges once. It proves \eqref{eq:integer-merge}
and gives the same coefficients as the divided derivatives in
characteristic zero.

The change of coordinates then turns, for example, \(R_U(i,:)x_U\)
into \(R_U(i,:)S_{U,I}x_I=R_I(i,:)x_I\). The column factors change in
the same way. It follows that
\[
 \mathcal P_I=\rho_I\!\left(
          \sum_{\ell=0}^{|I|}\frac{\mathcal D_I^\ell}{\ell!}
                    \mathcal P_U\mathcal P_V\right).
\]
In positive characteristic the sum in parentheses is understood
as \eqref{eq:integer-merge} applied to each monomial.
Each matching on \(I\) occurs in this identity through its unique
restrictions to \(U,V\) and its unique set of cross edges.

If \(p_J(w)=[w]\mathcal P_J\), taking \([w]\) at the parent gives
\[
 p_I(w)=\sum_{u\in\mathcal M_U,\ v\in\mathcal M_V}
                p_U(u)p_V(v)t_I(u,v;w).
\]
Expanding this identity down the tree gives one sum over all node
monomials, with a factor \(l_i\) at each leaf and \(t_I\) at each merge.
At the root the coordinate variables are empty, so every matching with an
unmatched row or column gives zero. Taking its constant coefficient gives
\eqref{eq:tree-formula}.

For a node \(I\), the number of coordinates is at most
\[
 \sum_{h=0}^{|I|}\binom{h+r-1}{r-1}\binom{h+s-1}{s-1}
       =O_{r,s}((|I|+1)^{r+s-1}).
\]
Zero-rank factors only reduce this number. There are \(n-1\)
merges. The product at a merge uses at most \(2(r+s)\)
variables of total degree at most \(2n\). The number of ordinary
monomials is thus bounded by \(\binom{2n+2(r+s)}{2(r+s)}\). For each
monomial, the number of arrays \(M\) in \eqref{eq:integer-merge} is
at most \((n+1)^{4rs}\). The integer coefficients can be formed by
binomial coefficient tables and multiplication. Polynomial multiplication,
contraction and substitution of linear forms therefore take polynomially
many field operations for fixed \(r,s\). The rank data are obtained
by elimination on submatrices of \(A\).

It remains to find a tree from the matrix.
If \(k=0\), each singleton cut has width zero, so \(A\) is diagonal
and its permanent is the product of its diagonal entries. Assume now that \(k\geq1\).
For \(I\subseteq[n]\), put
\(f(I)=\rk A[I^c,I]\). Group column \(i\) of \(A\) with the \(i\)-th
unit column. Let \(V_I\) be the span of the groups indexed by \(I\).
The selected unit columns span the coordinate subspace on \(I\).
After quotienting out this subspace, the remaining selected columns
are the columns of \(A[I^c,I]\). It follows that
\(\dim V_I=|I|+f(I)\).

We have \(V_{I\cup J}=V_I+V_J\) and \(V_{I\cap J}\subseteq V_I\cap V_J\).
The dimension formula gives
\[
 \dim V_I+\dim V_J\ \geq\ \dim V_{I\cup J}+\dim V_{I\cap J}.
\]
Subtracting \(|I|+|J|=|I\cup J|+|I\cap J|\) proves the same inequality
for \(f\). It also holds for \(I\mapsto f(I^c)\). Thus
\(w_A(I)=f(I)+f(I^c)\) is submodular, is symmetric under taking
complements, and vanishes on the empty set. It can be evaluated
by computing two ranks.

If \(\operatorname{rw}_{\F}(A)\le k\), the branch decomposition
algorithm of Korhonen and Oum \cite{KorhonenOum} finds a tree of
width at most \(k\) using \(w_A\) as its oracle. Rooting the
decomposition at a subdivided edge gives a rooted binary tree
of the same width. The one-index tree needs no construction.
For fixed \(k\) the search takes polynomial time. Over \(\mathbb Q\)
each rank query has polynomial bit cost. Both cut ranks on the
resulting tree are bounded by \(k\).
The preceding evaluation therefore applies with
\(r=s=\max(1,k)\).

For rational inputs choose the coordinate bases from rows and columns of
the original cut blocks. Solving the coordinate equations uses minors of
order bounded in terms of \(r,s\), so the resulting rational numbers have
polynomial bit length. There are polynomially many of them. Take a common
denominator \(Q\) for these numbers and the entries of \(A\). A product
in \eqref{eq:partial} has at most \(2n\) factors, so its coefficients
have a denominator dividing \(Q^{2n}\). To bound the numerator, there are
at most \((n+1)^n\) partial choices of columns before imposing the
matching condition, and each endpoint factor has at most \(\max(r,s)\)
terms. The logarithm of the number of products being added is therefore
polynomial in \(n\). At a merge at most \(|I|\) differentiations are
made, and the extra divisors \(\ell!\) all divide \(|I|!\). The
substitution contributes at most \(2n\) additional coordinate factors to
a product. These bounds also give polynomial bit length for the
intermediate coefficients.

\subsection{Proof of Theorem~\ref{thm:positive}}

We first check that reduction modulo the $p$th powers commutes
with the operations used in the matching formula.

\begin{lemma}\label{lem:positive-ideal}
Over a field of characteristic $p>0$, the ideal generated by
the $p$th powers of the variables is preserved
by multiplication, ordinary differentiation, and linear substitution.
The integer matching contraction descends to the quotient by this
ideal and is the product of the operators \eqref{eq:positive-pair}.
\end{lemma}

\begin{proof}
Multiplication preserves every ideal. If $z$ is one of the variables,
then
\[
 \partial_z(z^pf)=p z^{p-1}f+z^p\partial_zf
                =z^p\partial_zf.
\]
Differentiation in another variable also preserves the ideal. For
linear substitution, use
\[
 (a_1z_1+\cdots+a_dz_d)^p
        =a_1^pz_1^p+\cdots+a_d^pz_d^p.
\]
This identity holds over every field of characteristic $p$; no
assumption that the field has $p$ elements is being made.

For one pair of coordinate types, match $m$ of the $a$ row factors to
$m$ of the $b$ column factors. With edge weight $c$, the integer
coefficient is
\begin{equation}\label{eq:positive-integer-pair}
 \binom{a}{m}\binom{b}{m}m!c^m
   =\binom{a}{m}(b)_m c^m,
 \qquad (b)_m=b(b-1)\cdots(b-m+1).
\end{equation}
If $m\geq p$, this coefficient is divisible by $p$ (when $m>b$
the term is already zero). If $m<p$, it equals
$(a)_m(b)_m c^m/m!$, the coefficient in \eqref{eq:positive-pair}. It follows
that the full integer contraction for this pair reduces to
\eqref{eq:positive-pair}. In particular, this statement also applies before
high powers have been discarded.

Apply these pair contractions successively for all coordinate pairs.
An endpoint used by one pair is no longer available to a later pair.
Consequently each set of disjoint matches is counted once, grouped
according to its numbers of matches of each type. This is the
integer contraction in \eqref{eq:integer-merge}. Each surviving operator is
made from ordinary derivatives and scalars, so it preserves the
ideal. Linear substitution then preserves the corresponding parent
ideal. This proves the lemma.
\end{proof}

\begin{proof}[Proof of Theorem~\ref{thm:positive}]
Before taking the quotient, associate with a node $I$ the polynomial
\begin{equation}\label{eq:positive-partial}
 P_I(x_I,y_I)=
 \sum_{\substack{S,T\subseteq I\\ |S|=|T|}}
 \per A[S,T]
 \prod_{i\in I\setminus S}\bigl(R_I(i,:)x_I\bigr)
 \prod_{j\in I\setminus T}\bigl(C_I(j,:)y_I\bigr).
\end{equation}
A term records a matching inside $I$ and the unmatched rows and
columns. Both numbers of unmatched endpoints are $|I|-|S|$.
The polynomial therefore has equal degree in $x_I$ and $y_I$.

At a leaf this is the polynomial defining $l_i$. At a merge, first
choose a matching in each child. The remaining internal edges of the
parent go between the two children. For instance, differentiating
the row factor for $i\in U$ and the column factor for $j\in V$
with the coupling $B_{UV}$ gives
\[
 R_U(i,:)B_{UV}C_V(j,:)^{\mathsf T}=a_{ij}.
\]
The integer coefficients \eqref{eq:positive-integer-pair} count all possible
choices of distinct endpoints. Each matching on $I$ has a unique
restriction to its children and a unique set of cross edges. After
the substitutions $\rho_I$, every surviving endpoint has the parent
coordinate. Hence the untruncated polynomials obey the integer
matching contraction rule.

By Lemma~\ref{lem:positive-ideal}, that rule can be applied after reduction modulo
the $p$th powers at every stage. If $\overline P_I$ is the reduced
polynomial, it follows that
\[
 \overline P_I=\rho_I\bigl(E_I(\overline P_U\overline P_V)\bigr).
\]
Expanding coefficients down the tree proves \eqref{eq:positive-tree}.
At the root there are no coordinate variables. In \eqref{eq:positive-partial},
every term with an unmatched endpoint is zero, leaving exactly
$\per A$.

The basis bound follows from \eqref{eq:positive-dimension}. At a merge there
are at most $2k$ variables before substitution and at most $k$
after substitution. A polynomial with all exponents less than $p$
in $d$ variables has at most $p^d$ terms. Multiplication, the
at most $k^2$ pair contractions, and substitution of linear forms
can each be performed in $p^{O(k)}$ operations by direct polynomial
arithmetic in these bounded lists. For substitution, precompute
each linear form to powers $0,\ldots,p-1$, reducing after each
multiplication, and then substitute monomials one at a time. There
are fewer than $2n$ nodes. Rank factorizations and the coordinate
matrices use a polynomial number of field operations and zero
tests, independent of $k$. This proves the bound when a tree is supplied. For $k=0$,
the matrix is diagonal and the assertion is immediate.
The construction from the matrix uses an approximate decomposition.
The proof of Theorem~\ref{thm:tree} shows that \(w_A\) is a
symmetric submodular integer-valued function, with \(w_A(\varnothing)=0\).
Its singleton values are at most two. The approximation theorem
of Oum and Seymour \cite{OumSeymour}, in the form
\cite[Theorem 5.3]{KorhonenOum}, finds a decomposition of width
at most \(3k+c'\), where \(c'=\max_i w_A(\{i\})\), or certifies
that the width exceeds \(k\). Its cost here is
\(2^{O(k)}n^{O(1)}\), since two rank computations give the
value oracle and submodular minimization has polynomial cost.
Under the hypothesis \(\operatorname{rw}_{\F}(A)\leq k\),
we have \(c'\leq k\), so the output width is at most \(4k\).
Subdivide an edge to root the tree. The preceding evaluation
then costs \(p^{O(4k)}n^{O(1)}=p^{O(k)}n^{O(1)}\).
If \(k\) is not given, first check the diagonal case, then
try successive positive integers until a tree is returned.
This succeeds by the actual rank-width and has the same bound.

\end{proof}

\subsection{Proof of Theorem~\ref{thm:cauchy}}

\subsubsection{The determinant quotient for repeated nodes}
For balanced multiindices $a,b$, let $C_{a,b}$ be the ordinary Cauchy matrix
whose row nodes are $x_1,\ldots,x_n$ followed by $b_\beta$ copies
of each $\eta_\beta$, and whose column nodes are $y_1,\ldots,y_n$
followed by $a_\alpha$ copies of each $\zeta_\alpha$.
This matrix has repeated rows and columns. The next identity is the
confluent form of \eqref{eq:cauchy-borchardt}.

\begin{lemma}\label{lem:cauchy-confluent}
For admissible auxiliary nodes and $(a,b)\in\calI_n$,
\[
 \det K_1(a,b)\ne0,\qquad \per C_{a,b}=R_{a,b}.
\]
\end{lemma}

\begin{proof}
First replace every repeated node by distinct indeterminates and
apply \eqref{eq:cauchy-borchardt} to the resulting Cauchy matrix.
In each determinant divide by the Vandermonde product within each
row cluster and within each column cluster. The same factors are
removed from the numerator and denominator.

We describe the resulting specialization. If $h$ rows have nodes
$u+\epsilon_1,\ldots,u+\epsilon_h$, expand their entries in the
$\epsilon_i$. Multilinearity of the determinant shows that terms
with repeated Taylor orders vanish. The term of smallest total
degree has Taylor orders $0,1,\ldots,h-1$. Dividing by
$\prod_{i<j}(\epsilon_j-\epsilon_i)$ and then setting all
$\epsilon_i=0$ therefore replaces these rows by the divided
derivatives of orders $0,1,\ldots,h-1$. Apply this argument to
every row cluster, and then to every column cluster.
The required divided derivative is
\[
 \frac1{\ell!m!}\partial_u^\ell\partial_v^m(u-v)^{-q}
 =(-1)^\ell\frac{(q+\ell+m-1)!}{(q-1)!\ell!m!}
 (u-v)^{-q-\ell-m},
\]
which is exactly \eqref{eq:cauchy-kdef}.

The Cauchy determinant formula also shows that the divided
denominator has nonzero specialization. After removing the
within-cluster Vandermonde factors, its numerator is a product of
nonzero differences between distinct row nodes and between distinct
column nodes, with their appropriate multiplicities. Its denominator
is a product of differences between row and column nodes. Every such
difference is nonzero by admissibility.

The permanent itself has no singularity when nodes within a row
list or within a column list coalesce. Specializing Borchardt's
identity after the divisions consequently gives the claimed ratio.
This argument takes place in a rational function field and by
specialization holds over every field of characteristic zero.
\end{proof}

\subsubsection{Expansion according to the blocks}
The purpose of the auxiliary rows and columns is to collect many
submatrix permanents in a polynomial number of coefficients.
The following identity explains the sum in \eqref{eq:cauchy-dformula}.

\begin{lemma}\label{lem:cauchy-border}
For admissible auxiliary nodes and every $r\times s$ matrix $L$,
\begin{equation}\label{eq:cauchy-resolvent}
 \per(C+ULV)=\calD(x,y;\eta,\zeta;L).
\end{equation}
\end{lemma}

\begin{proof}
Introduce the polynomial
\begin{equation}\label{eq:cauchy-partial}
 Q(z,w)=\sum_{\substack{I,J\subseteq[n]\\|I|=|J|}}
 \per C[I^c,J^c]
 \prod_{i\in I}\left(\sum_\alpha U_{i\alpha}z_\alpha\right)
 \prod_{j\in J}\left(\sum_\beta V_{\beta j}w_\beta\right).
\end{equation}
Its monomials have the form $z^a w^b$ with $(a,b)\in\calI_n$.
The subsets in \eqref{eq:cauchy-partial} are used to prove the identity;
we will evaluate its coefficients by determinants instead of
enumerating these subsets.

Consider a perfect matching counted by $\per C_{a,b}$. Its edges
between original rows and original columns give one of the
submatrix terms in \eqref{eq:cauchy-partial}. Each remaining original row
uses an auxiliary column; each remaining original column uses an
auxiliary row. Suppose the used auxiliary columns and rows have
type multiplicities $a'$ and $b'$.
Put $c=a-a'$ and $d=b-b'$.

The auxiliary columns are labelled even when their entries agree.
For prescribed types, the number of ways of assigning them to
the original rows is $a!/c!$. Similarly the auxiliary row
assignments contribute $b!/d!$. The remaining auxiliary rows
and columns are matched using $W$. If $M_{\alpha\beta}$ records
the number of edges from row type $\beta$ to column type
$\alpha$, the number of such labelled matchings is
\[
 \frac{c!d!}{\prod_{\alpha,\beta}M_{\alpha\beta}!}.
\]
After including the weights and summing over $M$, their
contribution is $c!d!e_{W^{\mathsf T}}(c,d)$.
The product of the three factorial factors is $a!b!$.
We have therefore proved, coefficient by coefficient,
\begin{equation}\label{eq:cauchy-blockgenerating}
 \sum_{a,b}\frac{\per C_{a,b}}{a!b!}z^a w^b
 =\exp(z^{\mathsf T}W^{\mathsf T}w)Q(z,w).
\end{equation}
The sum here ranges over all balanced multiindices; it is a formal
identity. To recover $Q$ only degrees $|a|=|b|\leq n$ are needed.

Multiply \eqref{eq:cauchy-blockgenerating} by
$\exp(-z^{\mathsf T}W^{\mathsf T}w)$ and compare these coefficients.
Using Lemma~\ref{lem:cauchy-confluent} gives
\[
 [z^a w^b]Q(z,w)=
 \sum_{\substack{c\leq a,\ d\leq b\\|c|=|d|}}
 e_{-W^{\mathsf T}}(c,d)
 \frac{R_{a-c,b-d}}{(a-c)!(b-d)!}=Q_{a,b}.
\]

It remains to join the unmatched original rows and columns using
$ULV$. First expand the permanent according to the entries taken
from $C$ and those taken from $ULV$. This is the usual formula
\[
 \per(C+ULV)=\sum_{\substack{I,J\subseteq[n]\\|I|=|J|}}
 \per C[I^c,J^c]\,\per(ULV)[I,J].
\]
For each of the latter edges expand
$(ULV)_{ij}=\sum_{\alpha,\beta}U_{i\alpha}L_{\alpha\beta}V_{\beta j}$.
Once the $U$ and $V$ indices have been fixed, with multiplicities
$a,b$, the sum over bijections between the unmatched rows and
columns is $a!b!e_L(a,b)$. This follows from the same counting
by the table $M$ as above, now with weights $L_{\alpha\beta}$.
The contributions of the $C$, $U$ and $V$ entries are exactly
the coefficient $Q_{a,b}$ in \eqref{eq:cauchy-partial}. Summing gives
\eqref{eq:cauchy-dformula}, proving \eqref{eq:cauchy-resolvent}.
\end{proof}

\subsubsection{Obtaining the polynomial entries}
Let $U(t),V(t)$ be the matrices in \eqref{eq:cauchy-uv}, with the nodes
from \eqref{eq:cauchy-transformed}, and put
\[
 A(t)=C+U(t)H(t)V(t).
\]
We next show that $A(t)$ has no pole at zero and specializes there
to the matrix in \eqref{eq:cauchy-mainmatrix}. All expansions in the argument
are formal expansions at $t=0$.

For each auxiliary column,
\[
 U(t)_{i\alpha}=\frac{t}{x_it-\xi_\alpha}
 =-\sum_{a\geq0}x_i^a\xi_\alpha^{-a-1}t^{a+1}.
\]
The first $r$ coefficient rows in this expansion are the rows of
$B_x$. Consequently, for $0\leq a<r$,
\[
 \bigl(U(t)B_x^{-1}T_r(t)\bigr)_{ia}
 =x_i^a+O(t^{r-a}).
\]
To see the cancellation explicitly, the coefficient row of
$t^{k+1}$ for $k<r$ is $x_i^k$ times row $k$ of $B_x$.
Multiplication by $B_x^{-1}$ sends it to $x_i^k$ times the
$k$th unit row. Multiplication by $T_r(t)$ then divides column
$a$ by $t^{a+1}$. Terms with $k\geq r$ have exponent at least
$r-a>0$ in that column.

Likewise,
\[
 V(t)_{\beta j}=\frac{t}{\upsilon_\beta-y_jt}
 =\sum_{b\geq0}\upsilon_\beta^{-b-1}y_j^b t^{b+1},
\]
so
\[
 \bigl(T_s(t)B_y^{-1}V(t)\bigr)_{bj}
 =y_j^b+O(t^{s-b})\qquad(0\leq b<s).
\]
Using \eqref{eq:cauchy-transformed}, we conclude that
\begin{equation}\label{eq:cauchy-specialization}
 A(t)\text{ is regular at }0,\qquad A(0)=A.
\end{equation}

We now replace this specialization by a finite evaluation formula.
Write
\[
 d_i^x(t)=\prod_\alpha(x_it-\xi_\alpha),\qquad
 d_j^y(t)=\prod_\beta(\upsilon_\beta-y_jt).
\]
The expression
\begin{equation}\label{eq:cauchy-entryclearing}
 t^{r+s-2}d_i^x(t)d_j^y(t)A(t)_{ij}
\end{equation}
is a polynomial of degree at most $2r+2s-2$.
Indeed, the Cauchy term has this degree bound immediately.
A term of the perturbation containing $H_{ab}$ has a factor
$t^{-a-b}$ before multiplication by $t^{r+s-2}$, because the
two Cauchy factors contribute $t^2$ and the two diagonal factors
contribute $t^{-a-b-2}$. The remaining denominators are one
factor of $d_i^x$ and one factor of $d_j^y$. After clearing them,
its degree is at most $2r+2s-4-a-b$, and its lowest possible
exponent is nonnegative since $a+b\leq r+s-2$.

Multiplying \eqref{eq:cauchy-entryclearing} over a permutation uses every
$d_i^x$ and every $d_j^y$ once. Thus
\begin{equation}\label{eq:cauchy-fpoly}
 F(t)=t^{m_0}q(t)\per A(t)
\end{equation}
is a polynomial of degree at most $d_0$. By
\eqref{eq:cauchy-specialization},
\begin{equation}\label{eq:cauchy-targetcoef}
 [t^{m_0}]F(t)=q(0)\per A.
\end{equation}
Lagrange interpolation at $t_0,\ldots,t_{d_0}$ gives
$[t^{m_0}]F(t)=\sum_\ell\lambda_\ell F(t_\ell)$.
At each chosen point the auxiliary nodes are admissible, so
Lemma~\ref{lem:cauchy-border} evaluates $\per A(t_\ell)$ as $\calD$.
Substitution in \eqref{eq:cauchy-targetcoef} proves \eqref{eq:cauchy-mainformula}.
In particular, the final formula contains no limit.

\subsubsection{Number of operations and bit lengths}
At one evaluation point, there are $D_n$ quotients in
\eqref{eq:cauchy-rdef}. Their matrices have order at most $2n$, so their
construction and evaluation take $O_{r,s}(D_n n^3)$ field
operations. The factorials and the powers of each node difference
needed in \eqref{eq:cauchy-kdef} can be computed in advance.

Each of the $rs$ factors in \eqref{eq:cauchy-edef} has at most $n+1$
terms. When multiplying them we retain only indices in $\calI_n$.
This takes $O_{r,s}(nD_n)$ operations. The convolution
\eqref{eq:cauchy-qcoef} takes at most $O(D_n^2)$ operations; the final sum
\eqref{eq:cauchy-dformula} takes $O(D_n)$. Therefore one evaluation of
$\calD$ takes
\[
 O_{r,s}(n^{r+s+2}+n^{2r+2s-2})
\]
operations. There are $d_0+1=O_{r,s}(n)$ evaluation points.
Choosing them, forming $H(t_\ell)$ and $q(t_\ell)$, and expanding
the products in \eqref{eq:cauchy-lagrange} add a polynomial number of
operations bounded by \eqref{eq:cauchy-timebound}. This proves the field
operation bound.

For completeness, suppose the input is rational. The chosen
integers $t_\ell$ have $O_{r,s}(\log(n+1))$ bits. Each entry
of $K_q$ is obtained from a rational node difference by a power
of order $O(n)$ and a factorial quotient of order $O(n)$.
Its numerator and denominator therefore have polynomial bit
length. Clearing the entry denominators in each determinant
and using the determinant expansion bounds the resulting
integer bit lengths by a polynomial as well. Exact
fraction-free elimination consequently evaluates these
determinants in polynomial bit complexity.

The coefficient products in \eqref{eq:cauchy-edef} involve only $rs$
different weights, exponents at most $n$, and factorials at
most $n!$. A common denominator of polynomial bit length can
be chosen for each coefficient array. The determinant
quotients have polynomial bit length and there are only
polynomially many of them. Their products and sums in
\eqref{eq:cauchy-qcoef}--\eqref{eq:cauchy-mainformula} have common denominators
of polynomial bit length. The same holds for the
Lagrange coefficients, whose factors are differences of
integers of size $O_{r,s}(n)$. Performing the rational
operations exactly proves the asserted bit bound.
\qed

\subsubsection{The matrix classes and the rank bound}

For Corollary~\ref{cor:hilbert}, put \(x_i=i\) and \(y_j=1-j\).
Then \(x_i-y_j=i+j-1\), and \(R(i,j)=R(x_i,1-y_j)\).
For Theorem~\ref{thm:toeplitz}(6), use \(x_i=-i\) and
\(y_j=-\alpha-j\). Their difference is \(\alpha+j-i\),
and \(x_i-y_j-\alpha=j-i\). For
Theorem~\ref{thm:hankel}(3), use \(x_i=i\) and
\(y_j=-\alpha-j\); the corresponding difference minus
\(\alpha\) is \(i+j\). In each case both lists consist of
distinct nodes, and their disjointness is exactly the
nonvanishing of the denominators. Substituting in
Theorem~\ref{thm:cauchy} proves the three formulas and their
degree and complexity assertions.

\begin{proof}[Proof of Proposition~\ref{prop:cauchy-width}]
Every square submatrix of $C$ is a Cauchy matrix on distinct,
disjoint nodes and has nonzero determinant. Hence each rectangular
submatrix of $C$ has the largest possible rank. The matrix
$(P(x_i,y_j))$ factors as $XHY$, where $X_{ia}=x_i^a$ and
$Y_{bj}=y_j^b$, and has rank at most $q_0$.
The inequality $\rk(B+D)\geq\rk B-\rk D$ now gives
\[
 \rk A[I,I^c]\geq\min(|I|,n-|I|)-q_0.
\]
The same argument applies in the reverse direction. Both ranks
are nonnegative, proving the first assertion.

Every rooted binary tree with $n\geq2$ leaves has a node with
between $n/3$ and $2n/3$ leaves. To verify this, start at the root
and descend to a child having more than $2n/3$ leaves for as long
as such a child exists. The last node reached has more than
$2n/3$ leaves and each child has at most $2n/3$. Its larger
child has more than $n/3$ leaves and gives the required node.
For this node both $|I|$ and $n-|I|$ are at least
$\lceil n/3\rceil$. The second assertion follows.
\end{proof}

\subsection{Proof of Theorem~\ref{thm:path}}

Let \(I_i=\{1,\ldots,i\}\). Choose a basis \(F_i\) of the row space of
\(A[I_i,I_i^c]\), with zero rows appended to give \(r\) rows. Restrict
\(F_{i-1}\) to columns \(i+1,\ldots,n\). Each restricted row belongs to
the row space of \(F_i\), so there is a matrix \(B_i\) with
\[
 F_{i-1}[:,\{i+1,\ldots,n\}]=B_iF_i.
\]
Let \(q_i\) be the column of \(F_{i-1}\) with index \(i\), and choose
\(p_i\) so that the part of row \(i\) to its right is \(p_i^{\mathsf
T}F_i\). If \(j<i\), applying the restriction identity at cuts
\(j+1,\ldots,i-1\) gives
\[
 A[j,\{i,\ldots,n\}]
   =p_j^{\mathsf T}B_{j+1}\cdots B_{i-1}F_{i-1}.
\]
Taking column \(i\) gives the upper formula in \eqref{eq:generators}.
Applying the construction to \(A^{\mathsf T}\) gives the lower formula.
In forming \(F_i\), the restricted previous basis and the new row already
span the required space. Elimination involves at most \(r+1\) rows of
length \(n-i\). On the lower side there are at most \(s+1\) such rows.
Summing these costs over the cuts gives \(O_{r,s}(n^2)\).

Consider a matching between the rows and columns of \(I_{i-1}\). We give
it the product of its entries, together with a factor
\[
 \begin{aligned}
 p_j^{\mathsf T}B_{j+1}\cdots B_{i-1}x
   &\qquad\text{for each unmatched row }j,\\
 u_k^{\mathsf T}C_{k+1}\cdots C_{i-1}y
   &\qquad\text{for each unmatched column }k.
 \end{aligned}
\]
Let \(f\) be the sum of these products. When \(D_i\) acts on the factor
for row \(j\), the result is \(p_j^{\mathsf T}B_{j+1}\cdots
B_{i-1}q_i=a_{ji}\). In the product rule, that factor is removed and all
other factors are retained. Thus a summand of \(D_i f\) uses the edge
from row \(j\) to column \(i\). In the same way \(E_i\) removes the
factor for an earlier column \(k\) and gives the entry \(a_{ik}\). On any
factor that remains, changing \(x\) to \(B_i x\) or \(y\) to \(C_i y\)
appends \(B_i\) or \(C_i\) to its product. This is the change to the next
cut.

If both new vertices are joined to earlier vertices, the two
differentiations give \(\Phi_i(D_iE_i f)\). If only column \(i\) is
joined to an earlier row, the new row remains unmatched and gives
\(L_i\Phi_i(D_i f)\). If only row \(i\) is joined to an earlier column,
the new column gives \(M_i\Phi_i(E_i f)\). In the two cases with no edge
to an earlier vertex, either the new row and column are matched to each
other, giving \(d_i\Phi_i f\), or both are left unmatched, giving
\(L_iM_i\Phi_i f\). Every matching on \(I_i\) restricts to one matching
on \(I_{i-1}\) and belongs to exactly one of these five cases. In
particular the changes of coordinates act on the retained factors before
the new factors \(L_i,M_i\) are inserted.

After \(i\) steps, if there are \(h\) unmatched rows there are also \(h\)
unmatched columns. A completion has to match these rows to distinct
columns outside \(I_i\) and these columns to distinct rows outside
\(I_i\). Hence \(h\leq i\) and \(h\leq n-i\). A term with \(h>N\) at any
cut cannot contribute to the final constant term. This justifies \(\Pi\)
and gives the operator \(\mathcal U_i\). Starting at the constant
monomial, the coefficient of a monomial \(g\) after step \(i\) is
obtained by multiplying the preceding coefficient row vector by \(W_i\).
At the last cut the constant coefficient counts precisely the perfect
matchings. This gives the product in \eqref{eq:border}, and
\eqref{eq:product-border} gives the determinant. To count the states, a
monomial with degrees \(h,h\) is a choice of a weak composition of \(h\)
into \(r\) parts and one into \(s\) parts. Summing their product over
\(0\leq h\leq N\) gives \eqref{eq:dimension}.

We apply \(\mathcal U_i\) to its input polynomial without forming the
matrix \(W_i\). The numbers of monomials of degrees \(h,h\), or of
degrees differing by one after a differentiation, are
\(O_{r,s}((h+1)^{r+s-2})\). There are \(O_{r,s}(n^{r+s-1})\) monomials in
all. Row and column elimination reduce any of the fixed-order matrices
\(B_i,C_i\) to a diagonal matrix of zeros and ones. Reversing the
invertible operations expresses its substitution as a bounded number of
permutations, scalings, elementary additions of one variable to another
with a scalar multiple, and zero substitutions. For example, a
substitution \(x_j\mapsto x_j+\lambda x_k\) replaces \(x_j^a\) by its
binomial expansion with at most \(a+1\leq n+1\) terms. The cost of a step
is therefore \(O_{r,s}(n^{r+s})\). When \(r=s=1\), each linear
substitution is just a scaling and the cost of a step is \(O(n)\).
Together with the construction of the generators, these are the required
bounds.

\subsection{Proof of Corollary~\ref{cor:positive-ordered}}

Process the indices in order. At index $i$, the new row and column
may both be matched to previous endpoints, exactly one may be so
matched, they may form the diagonal edge, or both may remain
unmatched. These are respectively the first, next two, diagonal,
and $L_iM_i$ terms of \(\mathcal U_i^{(p)}\). Each derivative
chooses a previous endpoint, and $\Phi_i$ advances the remaining
endpoint factors by $B_i$ or $C_i$. The products of generators are
exactly \eqref{eq:generators}. Thus the constant-to-constant entry
counts complete matchings. Lemma~\ref{lem:positive-ideal} permits reduction after
each step. Finally, $\det(I-W^{(p)})=1$ and
$(I-W^{(p)})^{-1}=I+W^{(p)}+\cdots+(W^{(p)})^n$. A Schur complement proves the second
identity. Each step acts on a fixed basis of size $D_{p,r,s}$,
with $p^{O(r+s)}$ work to construct and apply its matrix.

\subsection{Proof of Theorem~\ref{thm:cycle-covers}}

A partial cover of \(I\) consists of directed edges inside \(I\) such
that each vertex has at most one incoming and one outgoing edge. Its
components are directed paths and cycles. We include an isolated vertex
among the paths, with the same start and end. Give each component a color
\(c\), include \(x_c\) for every vertex of that color, and multiply by
the edge weights. A cycle also contributes \(z_c\). For a path of color
\(c\) starting at \(j\) and ending at \(i\), include the factor
\begin{equation}\label{eq:cycle-endpoints}
 R_I(i,:)Z_I^{(c)}C_I(j,:)^{\mathsf T}.
\end{equation}
Let \(\mathcal Q_I\) be the sum of the resulting weights. For a single
vertex \(i\) of color \(c\), the isolated path gives \(x_c
R_{\{i\}}(i,:)Z_{\{i\}}^{(c)} C_{\{i\}}(i,:)^{\mathsf T}\). The loop
gives \(x_cz_ca_{ii}\). Adding over \(c\) gives the coefficients
\(l_i^{\mathrm{cyc}}\).

At a merge we view the child variables as the two diagonal blocks of
\(Z^{(c)}\). Write \(r_i,c_j\) for the row vectors of the end \(i\) and
the start \(j\), padded with zeros on the coordinates of the other child.
Thus \(r_iBc_j^{\mathsf T}\) is \(a_{ij}\) when \(i,j\) belong to
different children, and is zero when they belong to the same child. This
follows from \eqref{eq:coupling} and the zero diagonal blocks of \(B\).

In the first term of \eqref{eq:cycle-operator},
\(\partial_{\alpha\beta}^{(c)}\) applied to \(r_iZ^{(c)}c_j^{\mathsf T}\)
gives \(r_i(\alpha)c_j(\beta)\). After summation the contribution is
\(z_c r_iBc_j^{\mathsf T}\). This adds the edge from \(i\) to \(j\),
removes the two endpoints, and records the newly closed cycle. In the
second term take paths from \(j\) to \(i\) and from \(j'\) to \(i'\),
both of color \(c\). Each path factor is linear in \(Z^{(c)}\). To add
the edge \(i\to j'\), the derivatives take \(r_i(\alpha)c_j(\delta)\)
from the first factor and \(r_{i'}(\gamma)c_{j'}(\beta)\) from the
second. The sums with \(B_{\alpha\beta}\) and \(Z_{\gamma\delta}^{(c)}\)
give
\[
 (r_i Bc_{j'}^{\mathsf T})(r_{i'}Z^{(c)}c_j^{\mathsf T}).
\]
The first scalar is the edge weight, and the other factor records the new
path from \(j\) to \(i'\). Taking the derivatives on the factors in the
other order accounts for the choice \(i'\to j\). The derivatives cannot
take both factors from the same path. The two paths must have the same
color because each summand uses a single \(c\).

Every partial cover on \(I\) restricts to partial covers on \(U,V\), with
the colors inherited from its components. The edges lost on restriction
are the crossing edges. Starting with these restrictions, they can be
restored in any order. At an intermediate step the end and start needed
by the next edge are still available because that edge has yet to be
added. A cycle is recorded only when its last edge has been restored, so
its \(z_c\) is included once in every ordering. Consequently \(\ell\)
crossing edges are counted \(\ell!\) times by \(\mathcal E_I^\ell\). The
coordinate change takes a remaining factor to
\[
 r_i S_I Z_I^{(c)} T_I^{\mathsf T}c_j^{\mathsf T}
      =R_I(i,:)Z_I^{(c)}C_I(j,:)^{\mathsf T}.
\]
We obtain the identity
\[
 \mathcal Q_I=\eta_I\!\left(
       \sum_{\ell=0}^{|I|}\frac{\mathcal E_I^\ell}{\ell!}
          \mathcal Q_U\mathcal Q_V\right).
\]
Taking coefficients in the endpoint monomials gives the products of
\(l_i^{\mathrm{cyc}}\) and \(t_I^{\mathrm{cyc}}\) in
\eqref{eq:cycle-tree-formula}. The root has no endpoint variables, so
only closed covers survive there. Each closed cover is a permutation, and
summing its colors gives the cycle factors in \eqref{eq:cycle-colors}.

For one color, set the cycle and vertex markers to one and substitute
\(Z_I^{(1)}=x_Iy_I^{\mathsf T}\) in the partial-cover polynomial. Each
path contributes the independent row and column factors in
\eqref{eq:partial}. The starts and ends are exactly the unmatched columns
and rows in the associated matching. This identifies the two partial
polynomials in this specialization.

The matrix of endpoint variables at a merge has at most \(2r\) rows and
\(2s\) columns for each color. Together with the markers there are at
most \(4drs+2d\) variables. Each family of markers has degree at most
\(n\), as do the endpoint variables in the product of the child
polynomials. The number of coefficients is polynomial in \(n\) for fixed
\(d,r,s\). The rational bit bound can be seen as in the permanent
calculation. In a partial cover there are at most \(n\) edges and at most
two coordinate factors per path. The sum of the edge count and twice the
path count is at most \(2n\). There are at most \((n+1)^n\) choices of
outgoing edges before the indegree restriction, at most \(d^n\)
colorings, and a bounded number of coordinate terms per factor. A common
denominator of polynomial bit length raised to a power at most \(2n\)
suffices for the state coefficients. The linear substitutions and the
factorials at a merge give polynomial bit length for the intermediate
coefficients as well.

Along an ordering, the restriction of a full permutation to \(I_i\) has,
say, \(h\) paths. The \(h\) outgoing edges from their ends must use
distinct targets outside \(I_i\), so \(h\leq\min(i,n-i)\). Each path
contributes one endpoint variable to a monomial. We can therefore project
at every cut to total degree at most \(N\) without losing any term that
can survive at the root. In at most \(drs\) endpoint variables the number
of retained monomials is at most \(\binom{N+drs}{drs}\). When \(I_i\) is
obtained by merging \(I_{i-1}\) with the leaf \(\{i\}\), the leaf
polynomial is fixed and the operation is linear in the coefficient vector
at \(I_{i-1}\). Its coefficients are the stated entries of
\(W_i^{\mathrm{cyc}}\). The first leaf starts from the constant state.
Thus the constant-to-constant entry is \(F_A\), and
\eqref{eq:product-border} gives \eqref{eq:cycle-border} with the stated
order.

If a tree with the stated rank bounds exists, then
\(\operatorname{rw}_{\F}(A)\le r+s\). The construction in the proof of
Theorem~\ref{thm:tree} finds a tree on which both ranks are at most
\(r+s\). These bounds suffice for the polynomial-time calculation
of the cycle polynomial.

\subsection{Proof of Theorem~\ref{thm:fixed-cycles}}

Fix \(v_1,\ldots,v_k\) in \eqref{eq:fixed-cycle-sum}. When we expand the
products of permanents and determinants, a permutation can occur for a
choice of the sets only if each of its cycles lies in a single \(S_c\) or
\(T_c\). A cycle of length \(a\) in \(S_c\) contributes its usual weight.
In \(T_c\) its sign in the determinant is \((-1)^{a-1}\), and it accounts
for a factor \((-1)^a\) in \((-1)^{|T_c|}\). The total extra sign for
that cycle is therefore \(-1\).

If the permutation has a cycle with no mark, choose such a cycle by its
least vertex. Move the whole cycle from \(S_c\) to \(T_c\), or from
\(T_c\) to \(S_c\), within the pair it occupies. The union \(S_c\cup
T_c\) is unchanged and the cycle contains none of the prescribed \(v_c\).
The new choice therefore satisfies the same restrictions. This operation
pairs two occurrences of the same weight with opposite signs. If every
cycle has a mark, the cycle containing \(v_c\) has to lie in \(S_c\). Two
marks cannot lie in that cycle, since different marks are assigned to
different \(S_c\)'s. The minimum condition requires \(v_c\) to be the
least vertex of its cycle. The permutation then occurs once with positive
sign. Summing over marks gives \eqref{eq:fixed-cycle-sum}. With one mark
fixed and the minimum condition dropped, the same pairing leaves exactly
the single cycles.

For the determinant factor at a merge, fix one determinant color and let
\(L\subseteq I\) be the vertices having this color. Omit the color
superscript on \(u_I,v_I\). Introduce exterior variables
\(\theta_i,\phi_i\) for \(i\in L\), anticommuting with one another and
with the cut variables. Set \(\Theta_L=\prod_{i\in L}(\theta_i\phi_i)\).
The order of these pairs is immaterial, since each pair has even degree.
With interior variables written before the cut variables, consider
\begin{equation}\label{eq:fixed-cycle-partial}
 (-1)^{|L|}[\Theta_L]\exp\!\left(
 \sum_{i,j\in L}\theta_i a_{ij}\phi_j
 +\sum_{i\in L}\theta_i R_I(i,:)u_I
 +\sum_{j\in L}(C_I(j,:)v_I)\phi_j\right).
\end{equation}
All terms in the exponent have even degree and its exponential is a
finite sum. If a term uses \(a\) internal edges, it uses \(a\) of the
\(\theta\)'s and \(a\) of the \(\phi\)'s in those edges. For it to
survive extraction of \(\Theta_L\), the other \(|L|-a\) occurrences of
each kind must come from the cut terms. In particular the \(u_I\)-degree
and the \(v_I\)-degree are equal. At a leaf the coefficient before the
factor \(-1\) is
\[
 a_{ii}+(R_{\{i\}}(i,:)u_{\{i\}})(C_{\{i\}}(i,:)v_{\{i\}}).
\]
Indeed the two factors linear in the cut variables give
\(\theta_i(R_{\{i\}}(i,:)u_{\{i\}}) (C_{\{i\}}(i,:)v_{\{i\}})\phi_i\).
Moving \(\phi_i\) past the two cut factors has positive sign. This proves
the second line of \eqref{eq:fixed-cycle-leaf}. At the root there are no
cut terms. A surviving product is \(\prod_{i\in L}(\theta_i
a_{i,\sigma(i)}\phi_{\sigma(i)})\) for a permutation of \(L\). The even
factors can be ordered by \(i\), and putting the \(\phi\)'s in their
natural places multiplies the coefficient by
\(\operatorname{sgn}(\sigma)\). The orderings of the even factors cancel
the factorial in the exponential. The resulting root value is
\((-1)^{|L|}\det A[L,L]\).

For the merge let \(H\) be the sum of the two child exponents. The two
sets of interior variables are disjoint. Taking the two interior
coefficients of \(\exp H\) therefore gives the product of the child
polynomials, including their two signs. The interior coefficients have
even degree in the interior variables, so passing cut derivatives across
them introduces no sign. We can consequently perform the differentiation
before taking the interior coefficient.

For an ordered pair of different children \(J,K\) put
\[
 A_\alpha=\sum_{i\in L\cap J}R_J(i,\alpha)\theta_i,\qquad
 B_\beta=\sum_{j\in L\cap K}C_K(j,\beta)\phi_j .
\]
By the rule for a left derivative, \(\partial_{u_{J,\alpha}}H=-A_\alpha\)
and \(\partial_{v_{K,\beta}}H=B_\beta\). As \(H\) is even, it commutes
with its derivatives, and \(\partial_\xi\exp H=(\partial_\xi H)\exp H\).
It follows that
\[
 \partial_{v_{K,\beta}}\partial_{u_{J,\alpha}}\exp H
       =\partial_{v_{K,\beta}}(-A_\alpha\exp H)
       =A_\alpha B_\beta\exp H .
\]
The last sign is positive because the \(v\)-derivative crosses the odd
factor \(A_\alpha\). Multiplying by the entries of \(B_{JK}\) and summing
over the indices produces \(\sum_{i\in L\cap J,j\in L\cap K} \theta_i
a_{ij}\phi_j\) by \eqref{eq:coupling}. This is exactly the missing part
of the internal exponent. It contains no cut variables and has even
degree, so repeating the derivative operator multiplies by successive
powers of this part. Summing with the factorial denominators inserts its
exponential. Terms of power greater than \(|I|\) vanish, since they
repeat an interior variable. The final substitutions give the coordinate
terms at \(I\), by \(R_JS_{J,I}=R_I[J,:]\) and \(C_JT_{J,I}=C_I[J,:]\).
This proves the exterior merge rule.

For each permanent color the same merge is given by the proof of
Theorem~\ref{thm:tree}, with the matchings restricted to the vertices of
that color. For fixed marks, a marked vertex \(v_c\) is assigned the
permanent color of pair \(c\). An unmarked vertex \(i\) may take either
color of a pair with \(v_c<i\). Multiplication of the leaf sums lists the
allowed assignments, and linearity of each merge sums their contributions
in the coefficient arrays. There are \(\binom nk\) choices of marks. For
fixed \(k,r,s\), the number of variables at each merge is bounded,
commuting degrees are at most \(2n\), and exterior exponents are zero or
one. This gives polynomially many coefficients and operations. Over the
rationals the factors in each state coefficient come from the input and
the cut coordinates, with at most \(2n\) per product. The number of color
assignments is at most \((2k)^n\). The same common-denominator argument
as for \eqref{eq:partial}, with the extra factorials from a single merge,
gives the rational bit bound.

In the ordered calculation we can sum over marks in the transfer itself.
A pair can first appear only at its least vertex, in its permanent color.
This is exactly the transition that increases \(\ell\), and all other
transitions use a pair already introduced. Thus an allowed partition
determines a unique sequence of these choices, with \(\ell=0\) at the
start and \(\ell=k\) at the end. At cut \(I_i\) a term that can
contribute has at most \(\min(i,n-i)\) unmatched rows of any color.
Projection to degree at most \(N\) in the commuting variables therefore
preserves the root coefficient.

For one permanent color the count of surviving monomials is at most
\(D^+\). For one determinant color, if the two exterior degrees are
\(h,h\), there are at most \(\binom rh\binom sh\) choices of variables.
Summing gives
\[
 \sum_{h=0}^{\min(r,s)}\binom rh\binom sh
 =\sum_h\binom rh\binom{s}{s-h}=\binom{r+s}{s}=D^- .
\]
The middle sum counts the choices of \(s\) elements from two sets of
sizes \(r,s\), according to how many come from the first set. Including
all \(k\) pairs and \(\ell\), the number of coordinates is at most
\((k+1)(D^+D^-)^k=D_k\). The constant-to-constant entry of the ordered
transfer, with the prescribed end values of \(\ell\), is \(C_k(A)\).
Equation \eqref{eq:product-border} gives \eqref{eq:fixed-cycle-border}.
The summand for \(h\) in \(D^+\) is \(O_{r,s}((h+1)^{r+s-2})\), so
\(D^+=O_{r,s}(n^{r+s-1})\). This gives the stated order and, at \(k=1\),
the Hamiltonian-cycle bound.

\subsection{Proof of Corollary~\ref{cor:hook-immanants}}

Take \(d=p+q\) in \eqref{eq:cycle-colors}. In the first \(p\) colors use
\(z_c=1\). In color \(p+j\) use \(z_{p+j}=-1\) and \(x_{p+j}=-y_j\). For
a cycle of length \(\ell\), its factor becomes
\[
 p_\ell(x\mid y)=\sum_{i=1}^p x_i^\ell+
                    (-1)^{\ell-1}\sum_{j=1}^q y_j^\ell.
\]
For a cycle type \(\mu\), the product of these power sums is
\(p_\mu(x\mid y)\). The Frobenius character formula expands \(p_\mu\) as
\(\sum_{\lambda\vdash n} \chi_\lambda(\mu)s_\lambda\). Applying the
specialization to this identity for each permutation, its coefficient of
\(s_\lambda(x\mid y)\) is the defining weighted character sum. Thus the
specialized \(F_A\) is
\begin{equation}\label{eq:immanant-schur}
 \sum_{\lambda\vdash n}\imm_\lambda(A)s_\lambda(x\mid y).
\end{equation}

The Schur specialization of Berele and Regev \cite{BereleRegev} gives
\[
 s_\lambda(x\mid y)=\sum_{\mu\subseteq\lambda}
                 s_\mu(x)s_{\lambda'/\mu'}(y),
\]
where primes denote conjugate partitions. A term with \(s_\mu(x)\ne0\)
has \(\mu_{p+1}=0\). If \(\lambda_{p+1}>q\), column \(p+1\) of
\(\lambda'/\mu'\) then has length \(\lambda_{p+1}\). A column of a
semistandard tableau uses distinct entries, so this skew Schur polynomial
in \(q\) variables vanishes. Hence only the partitions in the
\((p,q)\)-hook contribute to \eqref{eq:immanant-schur}.

For a partition in the hook, the largest possible \(x\)-degree is
\(\lambda_1+\cdots+\lambda_p\). It is attained only for
\(\mu=(\lambda_1,\ldots,\lambda_p)\). In the conjugate skew diagram the
remaining boxes are the conjugate of the tail \((\lambda_{p+1},\ldots)\),
up to horizontal translation. The component of this largest \(x\)-degree
is therefore
\[
 s_{(\lambda_1,\ldots,\lambda_p)}(x)
 s_{(\lambda_{p+1},\lambda_{p+2},\ldots)'}(y).
\]
The second partition has at most \(q\) parts, and both Schur factors are
nonzero. If a linear relation among the specialized Schur functions
existed, take the largest \(x\)-degree of their highest components among
the terms with nonzero coefficients. At this degree one obtains a
relation among distinct products of Schur basis elements in the two
alphabets. Such products are independent, and the two partitions
determine \(\lambda\). This rules out the relation.

There are at most \((n+1)^{p+q}\) choices of the two partitions, by
choosing their parts. For fixed \(p,q\) their polynomials have
polynomially many monomials. Their coefficients can be obtained using
Jacobi--Trudi with the complete functions having generating series
\[
 \frac{\prod_{j=1}^q(1+y_j t)}{\prod_{i=1}^p(1-x_i t)} .
\]
For instance one can evaluate each polynomial on an interpolation grid
with \(n+1\) values in every variable and recover its coefficients. The
arrays are linearly independent by the preceding argument, so solving
\eqref{eq:immanant-schur} in the monomial basis recovers all the
immanants. The numbers of rows and columns of this system are polynomial
in \(n\). The Schur coefficients count tableaux with a fixed alphabet, so
their bit lengths are polynomial. For rational \(A\) the same holds for
the other side of the system by Theorem~\ref{thm:cycle-covers}, and
rational elimination gives the asserted bit bound.

For the ordinary hook \(\lambda=(n-k,1^k)\), put \(x=1\) and \(y=t\). In
the displayed Schur sum, \(\mu\) is a single row, and the skew diagram in
the second factor must have at most one cell per column. Its first column
has \(n-k-|\mu|\) cells and every other column has at most one. Thus
\(|\mu|=n-k\) or \(n-k-1\). Their two weights are \(t^k,t^{k+1}\),
respectively. This gives \(s_{(n-k,1^k)}(1\mid t)=(1+t)t^k\). The markers
in this specialization are those of \eqref{eq:ordinary-hooks}, which
follows on dividing \eqref{eq:immanant-schur} by \(1+t\).

\subsection{Proof of Proposition~\ref{prop:cycle-cut-rank}}

Write a left vertex type as \((a,b)\) and a right type as \((a',b')\). An
edge from left to right has weight one if \(a=a'\) and zero otherwise. An
edge from right to left has weight one if \(b'=b\) and zero otherwise.
The upper and lower cut blocks, with zero rows or columns omitted, are
submatrices of \(U_MU_{M'}^{\mathsf T}\) and \(V_{M'}V_M^{\mathsf T}\),
respectively. Indeed, at any initial cut all left vertices precede all
right vertices in the ordering, and an edge inside one half has weight
zero. This proves the rank bounds \(r,s\).

Every nonzero cycle alternates between the two halves. There are \(2h\)
vertices, so the number of cycles is at most \(h\), with equality only
when all cycles have length two. Such a cycle requires both type
coordinates to agree. Thus the coefficient of \(z^h\) in
\(P_{A_{M,M'}}(z)\) vanishes unless \(M=M'\). If they are equal, the
two-cycles pair the \(m_{ab}\) left vertices of each type bijectively
with the \(m_{ab}\) right vertices of that type. This gives the diagonal
entry \(\prod_{a,b}m_{ab}!\). If the matrix \(\mathcal K_h(z)\) has size
\(d_h\), the coefficient of \(z^{h d_h}\) in its determinant is the
determinant of this diagonal matrix. It is nonzero in characteristic
zero. The rank is therefore \(d_h\), the number of arrays with \(rs\)
entries and sum \(h\). This is the first equality of
\eqref{eq:cycle-cut-ranks}.

For the one-cycle lower bound, a perfect matching in each direction is
possible only when the left and right lists have the same
first-coordinate totals and the same second-coordinate totals. The
forward matching uses every right vertex once and preserves the first
coordinate. The backward matching gives the second condition. There are
\(B_h=\binom{h+r-1}{r-1}\binom{h+s-1}{s-1}\) pairs of totals. For one
such pair take any list with these totals and write it as \((a_i,b_i)\),
\(1\le i\le h\). Use the right list \((a_i,b_{i+1})\) with
\(b_{h+1}=b_1\). Then left vertex \(i\) can go to right vertex \(i\), and
this vertex can go to left vertex \(i+1\). All these edges have weight
one and make one cycle. Using this choice for every pair of totals gives
a submatrix of \([z]\mathcal K_h(z)\) with zero entries between different
pairs and positive entries on the diagonal. This gives rank at least
\(B_h\).

For the upper bound use the one-cycle formula with the first left vertex
in the permanent color. At the cut between the halves the couplings are
the identities on the first and second coordinate spaces. The sum over
allowed colors within the left half gives a polynomial depending only on
\(M\). The analogous right polynomial depends only on \(M'\). The merge
uses these two polynomials bilinearly and is independent of the two
lists. If \(j\) left vertices have the determinant color, all \(j\) of
their rows and columns are unmatched within the half. Their exterior
degrees are \(j,j\). The remaining \(h-j\) left vertices give commuting
degrees \(h-j,h-j\). The number of left coordinates is therefore bounded
by
\[
 \sum_{j=0}^{\min(h,r,s)}
 \binom rj\binom sj
 \binom{h-j+r-1}{r-1}\binom{h-j+s-1}{s-1}
 \leq B_h\sum_j\binom rj\binom sj
 =B_hD^-.
\]
In matrix terms the pairing factors the one-cycle coefficient matrix
through this coordinate space. Its rank is at most the dimension of that
space, giving the other inequality of \eqref{eq:cycle-cut-ranks}.

\subsection{Proof of Corollary~\ref{cor:inverse}}

For complementary index sets \(I,J\) and nonsingular \(B\), we use the
nullity identity \cite{StrangNguyen}
\begin{equation}\label{eq:nullity}
 \rk B[I,J]=\rk B^{-1}[I,J].
\end{equation}
To see it here, regard a vector in the kernel of \(B[I,J]\) as a full
vector \(v\) supported on \(J\). The \(I\)-coordinates of \(Bv\) vanish,
so \(Bv\) is also supported on \(J\). Its image under \(B^{-1}\) has zero
\(I\)-coordinates, which means that its restriction to \(J\) lies in the
kernel of \(B^{-1}[I,J]\). Conversely, starting from that kernel and
multiplying by \(B^{-1}\) gives the original kernel. The two maps are
inverse. Since the matrices have the same number of columns, equality of
nullities gives \eqref{eq:nullity}. We also have
\begin{equation}\label{eq:subadd}
 \rk\big((AB)[I,J]\big)\leq \rk A[I,J]+\rk B[I,J].
\end{equation}
Indeed \((AB)[I,J]=A[I,I]B[I,J]+A[I,J]B[J,J]\), and the two terms have
ranks at most the two ranks on the right.

If \(B\) has open lower and upper widths \(p,q\), the upper cross block
at a cut after \(i\) has entries only in the last \(q\) rows of \(I_i\).
The lower block has entries only in the first \(p\) rows after the cut.
Their ranks are at most \(q,p\), and \eqref{eq:nullity} gives the same
bounds for the inverse. This proves the first assertion using
Theorem~\ref{thm:path}.

For the tridiagonal cost, write \(\theta_i=\det B[1:i,1:i]\) and
\(\varphi_i=\det B[i:n,i:n]\), with \(\theta_0=\varphi_{n+1}=1\).
Expanding a leading minor in its last row gives
\[
 \theta_i=b_{ii}\theta_{i-1}
          -b_{i,i-1}b_{i-1,i}\theta_{i-2}\qquad(i\ge2).
\]
The trailing minors satisfy the reversed formula, so all of these numbers
take \(O(n)\) operations. For \(i\leq j\), the cofactor giving
\((B^{-1})_{ij}\) deletes row \(j\) and column \(i\). Divide the
remaining rows into \(1,\ldots,i-1\), \(i,\ldots,j-1\), and
\(j+1,\ldots,n\). The corresponding column groups are \(1,\ldots,i-1\),
\(i+1,\ldots,j\), and \(j+1,\ldots,n\). The resulting matrix is block
lower triangular. Its middle block is lower triangular with diagonal
\(b_{k,k+1}\), \(i\le k<j\). This gives
\[
 (B^{-1})_{ij}
  =\frac{(-1)^{i+j}\theta_{i-1}\varphi_{j+1}}{\theta_n}
                     \prod_{k=i}^{j-1}b_{k,k+1}.
\]
For \(i=j\) the middle block is empty. The lower entries are obtained by
transposing the argument, with the lower diagonal in the product. See
also \cite{Usmani}. For fixed \(i\) one can successively multiply the
factors as \(j\) increases. Thus all the entries of \(B^{-1}\) take
\(O(n^2)\) operations, with division only by \(\theta_n\ne0\). Its
ordered cut ranks are at most one, so Theorem~\ref{thm:path} gives the
same bound for computing its permanent.

For part (2), root the tree decomposition and assign each vertex to the
highest bag containing it. This bag is unique because the bags containing
the vertex form a connected subtree. At a bag \(D\) take the vertex sets
assigned below its individual children and the singleton vertices
assigned to \(D\) itself. Combine these sets by binary unions. Applying
this construction from the leaves of the decomposition upward, omitting
empty sets and suppressing unary nodes, gives a binary tree on the
vertices.

Consider a cut arising from one of the binary unions at \(D\). A vertex
outside \(D\) has all its bags in one component of the decomposition tree
with \(D\) deleted. If two such vertices are adjacent, some bag contains
both, so they lie in the same component. Each component is entirely on
one side of the cut for vertices outside \(D\), because each child set
was included or excluded as a whole. Thus a crossing edge has an endpoint
in \(D\). In the block from \(I\) to \(J\) this means that a nonzero
entry is in a row belonging to \(D\cap I\) or a column belonging to
\(D\cap J\). Splitting off the first group of rows bounds the rank by
\(|D\cap I|+|D\cap J|=|D|\). The same holds for the reverse block. By
\eqref{eq:nullity} the inverse has the same bounds on this tree.
Theorem~\ref{thm:tree} now applies. In particular a tree support graph
has a decomposition with bags of size at most two, by using a bag for
each vertex and its parent in a rooting of that graph.

\subsection{The scalar and block band formulas}

We prove the transfer formulas in Theorem~\ref{thm:circulant}(3),
Theorem~\ref{thm:toeplitz}(1) and Theorem~\ref{thm:block}(2).

Consider the \(b\) rows in block row \(i\). Their available columns lie
in \(w+1\) consecutive blocks. Number these columns relative to the first
block, from \(0\) to \(b(w+1)-1\). The earlier rows could have used
columns in the first \(w\) blocks of this window, and we denote their set
by \(U\). They cannot have used the last block, since the rightmost
allowed block increases by one at each step. Choose \(b\) unused columns
\(S\) for the current \(b\) rows. For this fixed set the sum of the
weights of the bijections from the rows to \(S\) is \(\per
D[\{0,\ldots,b-1\},S]\). In particular, the weight of a transition
already includes all orders of assigning these columns to the \(b\) rows.

At the next step the first block has left the window. Its columns must
all lie in \(U\cup S\), since no later row can fill a missing one. Remove
those \(b\) columns from \(U\cup S\) and subtract \(b\) from the other
relative indices. The resulting set is \(V\) in
\eqref{eq:block-transition}. As \(|S|=b\) and \(S\cap U\) is empty, we
have \(|V|=|U|\). Thus the transition stays in the sector indexed by this
common size.

For an open band, the first block in the window for row \(i\) is \(i-p\).
At row zero declare the \(p\) blocks with negative indices occupied. This
gives the initial state \(B\), and prevents the choice of any column with
negative index. After the last row, the window starts at \(m-p\). State
\(B\) now means that all blocks still in the window below \(m\) are full
and those from \(m\) onward are empty. Every earlier block has already
left the window and was full when it left. Even with the untruncated
weights, no column of index at least \(m\) can be used on a path from
\(B\) to \(B\). Such a column would still lie in the window at the final
step and would appear in the terminal state outside \(B\). Hence the
walks from \(B\) to \(B\), together with the choices counted in each
transition, give exactly the matchings in \eqref{eq:open-product} and
\eqref{eq:block-open}.

For the scalar cost with \(p\) fixed, after \(i\ge p\) rows the state
specifies \(p\) used columns beyond the first \(i-p\) columns. There are
\(O_p(n^p)\) possible states. If the next column to leave the window is
unused, the current row must use it and there is at most one transition.
If that column is already used, only \(p-1\) other columns have to be
specified in the state. There are \(O_p(n^{p-1})\) such states when
\(p\ge1\), and at most \(n\) choices for the current row from each. Thus
the work at a row is \(O_p(n^p)\). For \(p=0\) only the forced choice
occurs. During the first \(p\) rows at most \(p\) actual columns are
chosen, so the transitions at each of these steps also number
\(O_p(n^p)\). Summing over the rows proves the operation bound.

For a cyclic band, each chosen block displacement has a unique
representative in \(0,\ldots,w\) because \(w<m\). Repeat the choices on
block rows indexed by \(\mathbb Z\), using these displacements and the
same internal column labels in every period. This is a matching of the
lifted columns. In fact, each target column has exactly one preimage
modulo the period, and precisely one translate of that preimage reaches
the given lifted column. The occupied sets at successive cuts give a
closed walk of length \(m\) in a fixed sector.

Conversely, start from a closed walk and repeat it in both directions. At
each transition keep the actual choice of columns counted by that
transition. After \(w\) steps none of the columns in an initial state can
remain in the window. Thus each occupied position has arisen from a
column used at an earlier step, and no column is chosen a second time.
The condition on the departing block makes each column used at least
once. Passing back to the cyclic indices gives the required matching.
Taking traces counts all initial states and yields \eqref{eq:cyclic} and
\eqref{eq:block-cyclic}. For a scalar matching an edge with displacement
\(d\) appears in the occupied sets at the \(d\) cuts between its row and
column. The sum of the displacements is therefore \(nk\) when each of the
\(n\) cuts has \(k\) occupied positions.

To derive the scalar coefficient formulas, let \(Q(z)=\sum_{j=1}^D
q_jz^j\), where \(\det(I-zT)=1-Q(z)\). The trace of \(T^m\) is the sum of
the \(m\)-th powers of its eigenvalues, counted with multiplicity.
Factoring \(1-Q(z)\) over an algebraic closure gives
\[
 -\log(1-Q(z))=\sum_{m\ge1}\frac{z^m}{m}\tr(T^m).
\]
In the formal expansion \(\sum_{\ell\ge1}Q(z)^\ell/\ell\), a term indexed
by \(\alpha\) with \(|\alpha|=\ell\) has multinomial coefficient
\(\ell!/\prod_j\alpha_j!\) before division by \(\ell\). Taking \(m[z^m]\)
gives \eqref{eq:lpoly}. Equivalently, if the \(\alpha_j\) specify numbers
of tiles of length \(j\), mark one of the \(|\alpha|\) boundaries. There
are \(m\) choices of its position and \(|\alpha|!/\prod_j\alpha_j!\)
ordered lists of tiles starting there. Dividing by \(|\alpha|\) counts
the unmarked tilings.

For the open case, the formal inverse is
\((I-zT)^{-1}=\sum_{m\ge0}z^mT^m\). The adjugate formula gives its
\(B,B\) entry as \(\det(I-zT')/\det(I-zT)\). Expanding
\((1-Q(z))^{-1}=\sum_{\ell\ge0}Q(z)^\ell\) in the same way, now without
division by \(\ell\), gives \eqref{eq:hpoly}. Multiplication by the
numerator yields \eqref{eq:open} and \eqref{eq:open-series}. The block
formula \eqref{eq:block-toeplitz} uses the same inverse entry with
\(|B|=bp\).

With a row period \(h\), the product for \(n=hk+r\) is \(M^kN_r\). Its
series in \(k\), at fixed \(r\), has \(B,B\) entry \(e_B^{\mathsf
T}(I-zM)^{-1}N_r e_B\). Expanding the denominator as above gives
\eqref{eq:periodic}. At fixed width and fixed block size the matrices
have fixed order, so repeated squaring gives the stated cost for their
powers.

\subsection{The residue case of Theorem~\ref{thm:circulant}}

Shift columns first so that \(u\) is subtracted from all the jumps. Since
\(g\) divides both \(m\) and \(h\), each new jump is a multiple of \(g\).
After grouping rows and columns by their residue classes modulo \(g\),
the matrix is a direct sum. Translation of both indices carries any class
to any other and preserves the weights, so the component permanents are
equal.

In the component with residue \(c\), write the indices in the form
\[
 c+th+jm\pmod n,\qquad 0\le t<L,\quad 0\le j<b.
\]
To check that these labels are distinct, reduce an equality first modulo
\(m\). It gives \(L\mid t-t'\), since \(L=m/g\) and \(\gcd(h/g,L)=1\).
Hence \(t=t'\), and equality modulo \(bm\) then gives \(j=j'\). There are
\(Lb=n/g\) such labels, so they exhaust the component.

In the lifted block indices \(t\in\mathbb Z\), a jump \(s-u=dh+\ell m\)
goes from \((t,j)\) to \((t+d,j+\ell)\), with the second coordinate read
modulo \(b\). This gives the blocks \(D_d\) of
Theorem~\ref{thm:circulant}(2). Since \(Lh=m(h/g)\), the label
\((t+L,j-h/g)\) represents the same column as \((t,j)\). Thus the
relative positions of occupied columns after a period are those of \(\tau
U\) if the starting state is \(U\).

As in the cyclic band proof, a matching of the component lifts to a walk,
now from \(U\) to \(\tau U\). For the converse, extend such a walk
periodically using the change \(j\mapsto j-h/g\) after every period. The
blocks \(D_d\) are circulant in their internal labels, so rotating all
internal labels preserves every transition and its weight. Each column is
therefore used once by the extended walk, by the same departing-block
condition as in the ordinary cyclic case. The result descends to the
component. Summing the entries of \(T_k^L\) indexed by \(U,\tau U\), and
then over \(k\), gives the component permanent. Taking its \(g\)-th power
proves \eqref{eq:cyclic-residue}.

There are \(2^{bw}\) states across all the sectors. For a fixed
state, form the transitions by choosing one column for each of
the \(b\) rows in the current block. Each row has at most \(|S|\)
choices. Discard choices with a repeated column or a violation
of the conditions in \eqref{eq:block-transition}, and add the
weights of the remaining choices to their corresponding entries.
There are at most \(|S|^b\) choices for each state. Repeated
squaring then takes at most \(O(2^{3bw}\log(n+1))\) operations
for all sectors together. If \(w\ge1\), then \(b\le bw\), so
both bounds are polynomial in \(n\) when \(|S|\) is fixed and
\(bw=O(\log(n+1))\).

The parameters for this bound can be found by testing the
divisors \(b\) of \(n\), the residues \(h\bmod m\), and
\(u\in S\). Indeed, in any representation, one can take the
offset with the smallest \(d_s\) as \(u\) and subtract that
smallest value from every \(d_s\). Discard a test if any congruence
\(d_sh\equiv s-u\bmod m\) has no solution. Otherwise take the
unique solution in \(0\le d_s<L\) for each \(s\) and use their
maximum for \(w\).

For \eqref{eq:progression} use \(b=1\), so that each \(D_d\) is the
scalar \(c_d\) and \(\tau\) is the identity. In the example of order
\(2m\) there are two internal labels, and the change of label is
subtraction of \(h\) modulo two. The jump \(m\) exchanges those labels
inside a block, giving the \(dP_2\) contribution to \(D_0\).

\subsection{The signed formulas of Theorem~\ref{thm:circulant}}

For the torus determinant method, see \cite{CodenottiRestaSparse,Kenyon}.
The map \((x,y)\mapsto sx+ty\bmod n\) is onto by the gcd condition. If
\((x,y)\in\Lambda\), reduction modulo \(g_0\) gives
\(ty\equiv0\pmod{g_0}\), and \(\gcd(t,g_0)=1\) gives \(g_0\mid y\). Write
\(y=g_0m\). Dividing the congruence by \(g_0\) gives \(x\equiv
u_0m\pmod{h_0}\). Thus each element of \(\Lambda\) is an integral
combination of \(b_1,b_2\). Conversely these two vectors satisfy the
congruence. To reduce any lattice point to the specified rectangle,
subtract a multiple of \(b_2\) to put \(y\) between \(0\) and \(g_0-1\),
and then subtract a multiple of \(b_1\) to put \(x\) between \(0\) and
\(h_0-1\). The congruence just proved also shows that two points of this
rectangle in the same coset must coincide. This proves the description of
\(\rho_i\).

Join each point of \(\mathbb Z^2\) to its east, north and northeast
neighbors. These segments have disjoint interiors, and their translates
give the supported off-diagonal edges on the quotient torus. The
nontrivial cycles of a contributing permutation are therefore disjoint
simple curves. The displacement of a lift of each cycle is the sum of its
steps. The steps have nonnegative coordinates and each increases at least
one coordinate, so the displacement is a nonzero element of \(\Lambda\).
In particular the curve is essential.

Cutting a torus along an essential simple curve gives a cylinder. An arc
joining the two boundary copies of a point on the curve closes after
gluing, and its intersection number with the original curve is \(1\) or
\(-1\). The homology class of the original curve therefore cannot be a
nontrivial integer multiple. For two disjoint curves the intersection
number is zero. In lattice coordinates this intersection number is the
determinant of their two coordinate vectors. Their primitive classes are
consequently equal up to sign. In our lift both vectors have nonnegative
coordinates, so they have the same sign. Thus the total displacement of
\(k\) nontrivial cycles is \(k\) times one primitive vector. In the basis
\(b_1,b_2\), if the total has coordinates \(H_1,H_2\), this says
\(k=\gcd(H_1,H_2)\), with the stated convention for the empty collection.

For the monomial \(M\), the number of east steps plus northeast steps is
\(R\), and the number of north steps plus northeast steps is \(S\). Their
total \((R,S)\) lies in \(\Lambda\), giving \(n\mid sR+tS\). An integer
\(j>0\) divides this vector in \(\Lambda\) precisely when \(j\) divides
\(R,S\) and \((R/j,S/j)\) still satisfies the defining congruence. The
last condition is that \(j\) divide \((sR+tS)/n\). This proves the
formula for \(k\) in ordinary coordinates. If the nontrivial cycle
lengths are \(m_1,\ldots,m_k\), their signs multiply to
\((-1)^{\sum_i(m_i-1)}\). As \(\sum_i m_i=r+z+v\), this is
\((-1)^{r+z+v-k}\), independent of the contributing permutation. Taking
coefficients proves \eqref{eq:additive-sign}.

Now keep arbitrary row weights. In the determinant expansion a nontrivial
cycle of length \(m\) has sign \((-1)^{m-1}\). Its \(m\) off-diagonal
entries supply \(m\) additional minus signs, leaving \(-1\). The
remaining signs in \(K_{\epsilon\eta}\) come from the lattice coordinates
of the edge differences. Around a cycle the differences \(\rho_i-\rho_j\)
telescope. Thus the sums of the edge coordinates over all cycles are
\(H_1,H_2\), and the ratio of the signed determinant term to the unsigned
permanent term is
\[
 (-1)^k(-1)^{\epsilon H_1+\eta H_2}.
\]
Since \(k=\gcd(H_1,H_2)\), the sign \((-1)^k\) is positive exactly when
both coordinates are even. The same parity rule is expressed by
\[
 (-1)^k
   =\frac{-1+(-1)^{H_1}+(-1)^{H_2}+(-1)^{H_1+H_2}}{2}.
\]
In the signed combination \eqref{eq:4det}, the multiplier is the product
of \((-1)^k\) with this last expression, and is therefore one. Fixed
points have their original weights throughout. This proves
\eqref{eq:4det} term by term. Each supported entry is signed using two
lattice coordinates, and the four determinants take \(O(n^3)\) arithmetic
operations. Fraction-free elimination gives polynomial bit complexity
over the rationals.

For constant weights, fix a solution \((X_0,Y_0)\) of the two root
equations in \eqref{eq:products} over an algebraic closure, and put
\(\chi(x,y)=X_0^xY_0^y\). Its values on \(b_1,b_2\) are
\((-1)^\epsilon,(-1)^\eta\). With \(Q_{ii}=\chi(\rho_i)\) on the
diagonal, the entry on an edge with step \(\delta\) acquires the factor
\(\chi(\rho_i)\chi(\delta)\chi(\rho_j)^{-1}\) under the following
conjugation:
\[
 K_{\epsilon\eta}=Q
   (aI-bX_0P_n^s-cY_0P_n^t-dX_0Y_0P_n^{s+t})Q^{-1}.
\]
This factor is \(\chi(\rho_i+\delta-\rho_j)\), which is exactly the
remaining sign for that edge. The diagonal is unchanged.

For each \(n\)-th root of unity \(\omega\), the vector with \(i\)-th
entry \(\omega^i\) is an eigenvector of \(P_n\) with eigenvalue
\(\omega\). The conjugate circulant therefore has the factor
\(a-bX-cY-dXY\) for \((X,Y)=(X_0\omega^s,Y_0\omega^t)\). These pairs
satisfy the two root equations since \(b_1,b_2\in\Lambda\). They are
distinct, as equality for two roots gives a ratio with \(s\)-th, \(t\)-th
and \(n\)-th powers all one. There are \(h_0g_0=n\) solutions of the two
root equations, so these are all of them. This proves
\eqref{eq:products}. For the explicit parametrization in the corollary,
\(x_\epsilon^{h_0}=\zeta^{\epsilon n}\) and
\(x_\epsilon^{u_0}y_{\epsilon\eta}^{g_0} =\zeta^{\eta n}\), so the chosen
starting pair satisfies the root equations.

For fixed \(X\), the equation for \(Y\) is \(Y^{g_0}=(-1)^\eta
X^{-u_0}\). Multiplying \(a-bX-(c+dX)Y\) over those roots, by the monic
polynomial for \(Y\), gives
\[
 (a-bX)^{g_0}-(-1)^\eta X^{-u_0}(c+dX)^{g_0}.
\]
The identity is polynomial in \(a-bX,c+dX\), including when the latter is
zero. Since \(X\) is invertible modulo \(f_\epsilon\), the value at each
root of \(f_\epsilon\) is unchanged when the expression is replaced by
its remainder. The product over those roots is the resultant in
\eqref{eq:resultant}. Modular powering uses \(O(\log(g_0+1))\)
multiplications of polynomials of degree at most \(h_0\). The classical
Euclidean resultant calculation uses \(O(h_0^2)\) arithmetic operations.
For \(s=1,t=q\), one has \(g_0=1,h_0=n,u_0=n-q\). The second root
equation then gives \(Y=(-1)^\eta X^{q-n}=(-1)^{\epsilon+\eta}X^q\).
Using the exponent \(\eta\) instead reindexes the two products with
\(\epsilon=1\), exchanging \(D_{10}\) and \(D_{11}\). Their coefficients
in \eqref{eq:4det} agree.

Finally, after translating an offset to zero, the three- and four-offset
cases of \(\mathcal D_n(A)\) separate into classes modulo
\(\gcd(n,s,t)\). Inside one class, division of the indices and offsets by
this gcd gives the connected case already proved. Taking the product
gives the formula for arbitrary row weights, and constant weights give
equal factors. For two shifts \(0,s\), in an orbit of addition by \(s\),
a row choosing the \(s\)-edge uses the next row's diagonal column. That
next row must also use the \(s\)-edge, and the choice continues around
the orbit. The only other choice is to use the diagonal throughout the
orbit. This gives the two terms in \eqref{eq:two-shift}.

\subsection{The specializations in Corollary~\ref{cor:circulant} and Corollary~\ref{cor:toeplitz}}

In sector \(k=1\), the total displacement is \(n\). Along any permutation
cycle, the sum of the displacements is a multiple of \(n\), since the
lifted walk returns to the same residue. For a nontrivial cycle this sum
is positive. Hence there is one nontrivial cycle in this sector, and its
steps go once around the circle in their positive order. The steps of
lengths \(d\) are intervals tiling the circle. A step of length \(d\) has
weight \(c_d\), and the \(d-1\) vertices skipped between its ends are
fixed, with weight \(c_0\) each. Its total contribution is
\(c_0^{d-1}c_d\).

If there are \(\alpha_d\) intervals of length \(d\), put
\(H=\sum_d\alpha_d\). A tiling with a marked boundary is specified by the
position of that boundary and the ordered list of intervals from it. Thus
there are \(nH!/\prod_d\alpha_d!\) marked tilings. Each tiling has \(H\)
boundaries that can be marked. Summing the unmarked weights gives
\(\cL_n(c_1,c_0c_2,\ldots,c_0^{w-1}c_w)\). The sectors \(0,w\) use only
steps \(0,w\), respectively, and give \(c_0^n,c_w^n\). Transposition
changes a displacement \(d\) into \(-d\), and a column shift then changes
it into \(w-d\). The total displacement \(nk\) becomes \(n(w-k)\). This
gives the formula for sector \(w-1\) by reversing the weights in sector
one.

For \(w=2\), the numbers of steps of lengths one and two are \(n-2j\) and
\(j\). In the marked-boundary count, \(H=n-j\), so the factor becomes
\(n\binom{n-j}{j}/(n-j)\). This gives \eqref{eq:three}. For \(w=3\) the
two middle sectors are the winding-one count and its reversal, giving
\eqref{eq:four}. If the numbers of steps of lengths one, two and three
are \(u,v,w\) in that count, the product of the weights is
\(a^{v+2w}b^uc^vd^w\). Reversing the coefficients gives the second
monomial in \eqref{eq:four-expanded}, with the same factorial factor. For
a prescribed monomial, the total displacement distinguishes the two
middle sectors. Reading its step counts in the appropriate sector gives
the stated individual coefficients.

For the open one-subdiagonal band, number rows and columns from one.
Suppose row one chooses column \(j+1\). If \(j\ge1\), column one can now
only be filled by row two, since rows after it cannot reach column one.
After that, column two can only be filled by row three if \(j\ge2\).
Continuing, row \(k+1\) must fill column \(k\) for \(1\le k\le j\). The
first \(j+1\) indices are thus matched among themselves and contribute
\(a^jc_j\). For \(j=0\) row one gives \(b\). Removing the first interval
leaves a band of the same type. Conversely, any composition into
intervals of the allowed lengths gives one such choice in each interval.
The product of the interval weights, summed over the compositions, gives
\eqref{eq:hessenberg}, and for lengths one, two and three its multinomial
expansion is \eqref{eq:four-open-expanded}.

For width four, the only sector not yet counted is sector two. Use the
order \(01,02,03,12,13,23\) for its two-element sets. For example, from
\(01\) the choices \(2,3,4\) give the new sets \(01,02,03\),
respectively. From \(12\) the only choice is zero, which gives \(01\).
Filling the other rows in the same way gives
\[
 T_{4,2}=
 \begin{pmatrix}
 c&d&e&0&0&0\\
 b&0&0&d&e&0\\
 0&b&0&c&0&e\\
 a&0&0&0&0&0\\
 0&a&0&0&0&0\\
 0&0&0&a&0&0
 \end{pmatrix}.
\]
The lower right three-by-three block of \(I-zT_{4,2}\) is unit lower
triangular. Eliminating the lower left block, the remaining determinant
is
\[
 \det\begin{pmatrix}
 1-cz&-dz&-ez\\
 -bz-adz^2&1-aez^2&0\\
 -acz^2-a^2ez^3&-bz&1
 \end{pmatrix}.
\]
For the numerator with state \(01\) deleted, the same elimination leaves
the lower right two-by-two minor of this last matrix. Its determinant is
\(1-aez^2\). For the denominator, put \(x=ae,\ y=bd,\ h=ad^2+b^2e\).
Expansion of the three-by-three determinant gives
\[
 (1-xz^2)\bigl(1-cz(1+xz^2)-x^2z^4\bigr)
                -yz^2(1+xz^2)-hz^3.
\]
This is \(1-\sum_{j=1}^6q_jz^j\) with the parameters in
\eqref{eq:penta-par}. It is also \((1-xz^2)(1-F(z))\) with \(F\) as in
\eqref{eq:firstreturn}. The trace expansion gives \eqref{eq:five}, and
the open quotient gives \eqref{eq:penta-open} and \eqref{eq:penta-gf}.

The factor \(1-xz^2\) cancels in the open quotient, so the coefficient of
\(z^n\) is that of \(\sum_{\ell=0}^n F(z)^\ell\). To expand it, take
\(i,j,k,\ell\) copies of the terms \(cz(1+xz^2)\), \(x^2z^4\),
\(yz^2(1+xz^2)/(1-xz^2)\), and \(hz^3/(1-xz^2)\), respectively. Their
multinomial factor is \((i+j+k+\ell)!/(i!j!k!\ell!)\). The remaining
power series are \((1+xz^2)^{i+k}\) and \((1-xz^2)^{-k-\ell}\). Taking
terms \(u,v\) from these series gives \(\binom{i+k}{u}B(k+\ell,v)\) and
the power \(z^{i+4j+2k+3\ell+2u+2v}\). This is
\eqref{eq:positive-expanded}. In the cyclic formula the logarithm of the
denominator splits into the logarithms of \(1-xz^2\) and \(1-F(z)\). The
first gives \(2x^{n/2}\) for even \(n\) and zero for odd \(n\). For the
second, the term with \(i+j+k+\ell\) factors is multiplied by
\(n/(i+j+k+\ell)\). This proves \eqref{eq:positive-cyclic}.

For the symmetric specialization \(a=e=v,b=d=u,c=t\), the numerator is
\(1-v^2z^2\). Substitution into the denominator just computed and
division by \(1+vz\) gives \(1-\sum_{j=1}^5s_jz^j\). The numerator after
division is \(1-vz\), giving \eqref{eq:symmetric} by the ordinary
composition expansion. The numerical values in the corollary come from
this quotient at \(t=u=v=1\). Each nonzero product in the permanent of
the resulting \((0,1)\) matrix represents one perfect matching of its
bipartite graph.

For a tridiagonal matrix, a nontrivial permutation cycle can only be a
transposition of neighbors. Indeed, at the least vertex in such a cycle
the preceding and following vertices must both be its immediate successor
in the path. A cycle of length greater than two would repeat that
successor. Thus the possible transpositions form a set \(E\) of disjoint
path edges. Their weights are \(x_i\), and the remaining vertices use the
diagonal. This proves \eqref{eq:tridiag}. For \(k\) transpositions,
subtract \(j-1\) from the \(j\)-th selected edge index, for \(1\le j\le
k\). This is a bijection with the \(k\)-subsets of \(\{1,\ldots,n-k\}\),
giving the binomial factor for constant weights. For offsets divisible by
\(r\), ordering by residues modulo \(r\) separates the matrix into
blocks. In the class of \(s\) the offsets \(-r,0,r\) become \(-1,0,1\) on
\(\ell_s\) indices. Multiplying their formulas gives
\eqref{eq:tridiag-stride}.

For even order \(2m\) with centrosymmetry, if edge \(m\) is unused the
two halves contribute \(p_m^2\). If it is used, it has weight \(x_m\) and
leaves two paths of length \(m-1\), giving \(x_mp_{m-1}^2\). For odd
order \(2m+1\), fixing the central vertex gives \(d_{m+1}p_m^2\).
Matching it to its left neighbor has weight \(x_m\) and leaves paths of
lengths \(m-1,m\). Reflection gives the same contribution when it is
matched to its right neighbor. This proves the two reflection formulas.

\subsection{The rational circulants of Corollary~\ref{cor:circulant}}

Consider a Laurent polynomial in \(P_n\) supported on offsets
\(-p,\ldots,q\). At an initial cut after \(i\), a forward edge from the
initial interval to its complement either moves by a positive offset and
starts in one of its last \(q\) rows, or wraps by a negative offset and
starts in one of its first \(p\) rows. Thus the upper cross block has
rank at most \(p+q\). The same argument after transposition bounds the
rank of the opposite block by \(p+q\). Applying \eqref{eq:nullity} to
\(g(P_n)\) bounds the cut ranks of its inverse. By \eqref{eq:subadd}, the
two ranks of \(f(P_n)g(P_n)^{-1}\) are bounded by the sums of the
corresponding bounds for \(f\) and \(g\). Equation \eqref{eq:border}
therefore applies with constant rank parameters.

The permutation \(P_n^h\) has \(\gcd(n,h)\) orbits, all of size
\(n/\gcd(n,h)\). Ordering each orbit along its cycle makes its
restriction the shift of that order. Both Laurent polynomials and the
inverse are block diagonal in this ordering with equal blocks.
Multiplicativity of the permanent over a direct sum gives the stated
power.

For the Eulerian formula, multiplication by \(I-\rho P_n\) gives
\[
 (I-\rho P_n)\sum_{d=0}^{n-1}\rho^dP_n^d
      =(1-\rho^n)I.
\]
For a permutation \(\sigma\), let \(d_i\in \{0,\ldots,n-1\}\) be the
residue of \(\sigma(i)-i\). Then \(d_i=\sigma(i)-i\) if \(\sigma(i)\ge
i\), and \(d_i=\sigma(i)-i+n\) otherwise. The differences \(\sigma(i)-i\)
sum to zero. Thus the term for \(\sigma\) in the inverse permanent has
weight
\[
 (1-\rho^n)^{-n}\rho^{\,n|\{i:\sigma(i)<i\}|}.
\]
The number of permutations with \(k\) such positions is the Eulerian
number \(\eul{n}{k}\) \cite{Stanley}. Summing the weights proves
\eqref{eq:eulerian}.

\subsection{The number of signed determinants}

Write the row vertices first and the column vertices second in the
skew-adjacency matrix of the bipartite support graph. Its Pfaffian is the
determinant of its upper-right block, up to one overall sign depending
only on the order. Changing the orientation of an edge changes the
corresponding sign in that block. Thus the least number in
\eqref{eq:signing-length} is the Pfaffian number of the support graph.
For any fixed collection of signings, comparing monomial coefficients
gives a rational linear system for the coefficients of the determinants.
Its consistency is unchanged by extension of the ground field in
characteristic zero. We may therefore use the real Pfaffian results for
the lower bounds.

For the connected additive support, denote the neighborhood of a row set
\(X\) by \(X+S\). For \(1\le |X|\le n-2\) we claim \(|X+S|\ge |X|+2\).
This follows at once if \(X+S\) is the whole group. If the stabilizer
\(D\) of \(X+S\) is trivial, Kneser's inequality gives \(|X+S|\ge
|X|+|S|-1=|X|+3\). If \(D\) is nontrivial and proper, connectedness means
\(S\) meets at least two \(D\)-cosets. Indeed, if it met only one, every
difference of two shifts would lie in \(D\). Kneser's inequality now
gives
\[
 |X+S|\ge |X+D|+|S+D|-|D|
         \ge |X|+|D|\ge |X|+2.
\]
Fix any two disjoint edges and remove their rows and columns. A set of
surviving rows has lost at most two neighbors. The inequality proves
Hall's condition on the remaining graph. Its matching, together with the
two prescribed edges, is a perfect matching. The connected support graph
is therefore a brace.

A Pfaffian brace is either the Heawood graph or is obtained from planar
braces by trisums along four-cycles \cite{RST}. Here we only need the
degree consequence of this description. A planar brace on \(v>4\)
vertices has minimum degree at least three. Since it is bipartite,
Euler's inequality gives at most \(2v-4\) edges. If \(v_3\) vertices have
degree three, the degree sum is at least \(4v-v_3\) and at most \(4v-8\).
Thus \(v_3\ge8\). In a trisum the vertices off the common four-cycle
retain their degrees, and the three operands meet only on that cycle.
Each operand has vertices off the cycle by the definition of a trisum. By
induction each operand has at least eight degree-three vertices, at most
four of which are on the common cycle. In particular the result has
degree-three vertices. The Heawood graph itself is cubic. Our graph is
four-regular and has more than four vertices, so it cannot be a Pfaffian
brace. Norine's theorem \cite{Norine} implies that a graph admitting a
combination of at most three signings is Pfaffian. Thus at least four
determinants are required, and \eqref{eq:4det} supplies four.

For the support \(\{0,1,2,4\}\), subtract two from the offsets by a
column shift. The available offsets are now \(-2,-1,0,2\). Divide the
ordered rows and columns into \(\lfloor n/5\rfloor\) groups of five,
leaving the same remaining indices on the two sides. In a group use row
labels \(r_0,\ldots,r_4\) and column labels \(c_0,\ldots,c_4\). Take the
\(K_{3,3}\) on rows \(r_1,r_2,r_3\) and columns \(c_1,c_2,c_3\), except
that the edges \(r_1c_2\) and \(r_2c_3\) are replaced by
\[
 r_1-c_0-r_0-c_2,\qquad r_2-c_4-r_4-c_3 .
\]
The three offsets on the first path are \(-1,0,2\), and on the second are
\(2,0,-1\). The seven retained direct edges have offsets in
\(\{-2,-1,0,2\}\). Thus every group contains a subdivision of \(K_{3,3}\)
of the required kind. After deleting the vertices of their union, the
remaining rows and columns can be matched on the diagonal. The union is
consequently conformal, and \cite[Corollary 4.3]{PantSingh} gives the
stated lower bound.

To see that there is a finite upper bound, regard each signed determinant
coefficient as a function of the independent edge signs. For a fixed
matching, averaging the signed determinants against the product of its
edge signs isolates that matching monomial, up to its permutation sign.
Every other monomial has an edge in the symmetric difference and averages
to zero. Thus the signed determinants span the matching monomials. There
are at most \(4^n\) such monomials. For evaluation of the permanent
itself, use \eqref{eq:block-cyclic} with \(b=1,w=4\). There are a fixed
number of transfer entries at each row, so the cost is \(O(n)\)
arithmetic operations.

\subsection{Proof of Theorem~\ref{thm:sign-classification}}

For a monomial in the shift variables, the permanent coefficient counts
all its contributing permutations. The determinant coefficient is the sum
of their signs. Its absolute value equals the count precisely when all
those signs agree. A common translation of the offsets permutes the
columns and changes every determinant sign by the same factor. An
automorphism of the group amounts to a simultaneous ordering of rows and
columns. We may use both operations when testing the equality.

We first consider one connected component, so \(H=G\). With one shift
there is just one permutation. With two shifts, translate them to
\(0,s\). The rows on each orbit of addition by \(s\) use either all
diagonal entries or all \(s\)-entries, as in \eqref{eq:two-shift}. The
number of orbits of the second kind is determined by the multiplicity of
\(s\). Since the orbits have the same size, the sign is determined. For
three shifts write \(S=\{0,s,t\}\), and for four shifts with equal pair
sums write \(S=\{0,s,t,s+t\}\). In both cases
\[
 \Lambda=\{(x,y)\in\mathbb Z^2:xs+yt=0\text{ in }G\}
\]
has index \(|G|\), and the quotient grid is the torus used in the proof
of \eqref{eq:additive-sign}. A monomial specifies the total lifted
displacement. Its divisibility in \(\Lambda\) gives the number of
nontrivial cycles by the primitive-vector argument. It also specifies how
many vertices are moved. These two numbers determine the sign, giving the
desired equality in these cases.

Suppose \(|S|=4\) and no disjoint pairs in \(S\) have the same sum. A
nonzero period would have an orbit size dividing four. In that case there
is a period \(b\) of order two, and \(S\) is the union of pairs
\(\{x,x+b\},\{y,y+b\}\). The cross pairs \(\{x,y+b\},\{x+b,y\}\) have the
same sum, a contradiction. Thus \(S\) is aperiodic. Also there is at most
one \(z\) with \(z-S=S\). Indeed the product of the reflections
\(u\mapsto z-u\) and \(u\mapsto z'-u\) translates by \(z-z'\).

If there is such a reflection with a nonfixed pair, choose a cyclic
ordering \(s_0,s_1,s_2,s_3\) of \(S\) in which that pair is adjacent. If
every element is fixed, translating one to zero gives a set
\(\{0,a,b,c\}\) in an elementary abelian two-group. The three nonzero
elements generate this group and are independent, since \(a+b=c\) would
give equal pair sums. The cycle \((0\,a\,a+b\,b)\) uses \(a,b\) twice
each and has negative sign. The two transpositions \((0\,a)(c\,c+b)\)
have the same multiplicities and positive sign. All remaining vertices
are fixed, so this gives the required cancellation. When there is no
reflection, any cyclic ordering will do.

With the chosen ordering put \(y_i=s_i+s_{i+1}\), using subscripts modulo
four, and \(Y=\{y_0,y_1,y_2,y_3\}\). Adjacent sums cannot be equal since
they share one summand. Equality of opposite sums would contradict the
condition on disjoint pairs. The four sums are therefore distinct, and
they satisfy \(y_0+y_2=y_1+y_3\). Match \(s_i\) to \(y_i\) in one
matching and to \(y_{i-1}\) in the other. The offsets used are
\(s_{i+1}\) in the first case and \(s_{i-1}\) in the second. Thus each
shift is used once, and the assignments to the four columns differ by a
four-cycle. To extend both matchings using the same edges on the
remaining vertices, it suffices by Hall's theorem to prove
\begin{equation}\label{eq:hall}
 |(X+S)\setminus Y|\ \geq\ |X|\qquad(X\subseteq G\setminus S).
\end{equation}

For \(X\ne\varnothing\) write \(\Delta=|X+S|-|X|\). Since at most four
columns are deleted, we only need consider \(\Delta\leq3\). If \(|X|=1\),
then \(X+S\) is a translate of \(S\) and differs from the additive set
\(Y\). Thus at least one of its elements remains. If \(X+S=G\), then
\(\Delta\ge4\) because \(X\subseteq G\setminus S\). If
\(X+S=G\setminus\{z\}\), then \(\Delta\le3\) requires \(X=G\setminus S\)
and \(\Delta=3\). Since \(z-s\notin X\) for each \(s\in S\), we have
\(z-S=S\). The chosen ordering places a pair with sum \(z\) next to each
other, so \(z\in Y\). Only three elements of \(Y\) are then deleted from
\(X+S\), and \eqref{eq:hall} follows.

Suppose next that the stabilizer \(D\) of \(X+S\) is nontrivial and
proper. Kneser's inequality \cite{Kemperman,Lev} gives
\[
 \Delta\geq |X+D|-|X|+|S+D|-|D|.
\]
The connected set \(S\) meets at least two \(D\)-cosets. If \(|D|=2\),
meeting only two cosets would fill both, making \(S\) additive. Meeting
at least three gives \(\Delta\ge4\). If \(|D|\ge4\), meeting at least two
gives the same conclusion. The remaining case with \(\Delta\le3\) has
\(|D|=3\), two cosets met by \(S\), and \(X+D=X\). A split of \(S\) into
two points in each coset would again make it additive. Indeed, the two
differences in a group of order three agree up to sign, and ordering the
pairs with equal differences gives equal cross sums. Hence, after
translation, \(S=D\cup\{t\}\). Translating \(X,S\) by the same element
translates both \(X+S\) and \(Y\) by twice that element, preserving
Hall's inequality. In \(G/D\) the image of \(X\) omits \(0,\bar t\),
since \(X+D=X\) and \(X\) avoids \(S\). Its sum with \(\{0,\bar t\}\)
therefore omits \(\bar t\). Two points of \(Y\) are sums of \(t\) with a
neighbor in \(D\) and lie in that missing coset. Thus \(|(X+S)\cap
Y|\le2<\Delta=3\).

It remains to take aperiodic \(X+S\), with \(|X|\ge2\) and \(|X+S|\le
|G|-2\). Kneser's inequality gives \(\Delta=3\), so
\(|X+S|=|X|+|S|-1\). This equality and the aperiodicity of
\(X+S\) are the two hypotheses in Kemperman's theorem
\cite[Theorem C]{Lev}. An elementary pair with
both sets of size at least two is either a pair of progressions, or lies
in subgroup cosets whose sum fills a coset or omits one point. Here a
four-point progression is additive. In the coset cases the subgroup would
have to be \(G\), since the differences of \(S\) generate \(G\), and the
size bound excludes these cases. Thus the theorem gives a proper nonzero
subgroup \(D\). Each of \(X,S\) is a union of full \(D\)-cosets and one
distinguished part contained in a further coset. In the quotient, the sum
of the two distinguished cosets has a unique representation. All other
cosets in \(X+S\) are full because at least one of their summands is
full.

The set \(S\) has a full coset, since its differences generate \(G\), and
its distinguished part cannot itself be full, since then \(X+S\) would be
periodic. Write \(4=a|D|+b\), where \(a\ge1\) is the number
of full cosets and \(1\le b<|D|\) is the size of the
distinguished part. Since \(|D|\ge2\), the only possibility is
\(|D|=3\), \(a=b=1\). Thus \(S=D\cup\{t\}\) after
translation. Write \(X_0\) for its distinguished part, and use bars
for quotient images. The exceptional coset of \(X+S\) has exactly
\(|X_0|\) points, since its only source is \(X_0+\{t\}\). Consequently
\[
 \begin{split}
 |X+S|&=3\bigl(|\bar X+\{0,\bar t\}|-1\bigr)+|X_0|,\\
 |X|&=3\bigl(|\bar X|-1\bigr)+|X_0|.
 \end{split}
\]
Their difference is three, giving \(|\bar X+\{0,\bar t\}|=|\bar X|+1\).

The quotient is cyclic with generator \(\bar t\). The set \(\bar X\)
omits zero. If it also omits \(\bar t\), then its sum with \(\{0,\bar
t\}\) omits \(\bar t\), and again two points of \(Y\) are absent from
\(X+S\). If \(\bar X\) contains \(\bar t\), this must be the
distinguished coset of \(X_0\) because \(t\) itself is excluded. In the
cyclic order along \(\bar t\), the number of new elements obtained by
adding the translate \(\bar X+\bar t\) is the number of maximal intervals
of \(\bar X\). The displayed cardinality equality gives one interval. It
begins at \(\bar t\) because zero is absent. If \(2\bar t\) also belonged
to it, the sum \(2\bar t\) would have the representations \(\bar t+\bar
t\) and \(2\bar t+0\), contrary to uniqueness for the distinguished
cosets. Thus \(\bar X=\{\bar t\}\). If the quotient has order at least
three, \(X+S\) omits \(D\), which contains the two adjacent sums of
elements of \(S\) that avoid \(t\). If the quotient has order two, then
\(|X|\ge2\) and \(X\subseteq(t+D)\setminus\{t\}\) give \(|X+S|=5=|G|-1\),
already treated. This proves \eqref{eq:hall} in every case. The common
extension keeps the two signs opposite and gives the same shift counts.

For \(|S|\ge5\), take any five elements of \(S\). If they contain a
nonadditive four-set, translate it into the subgroup generated by its
differences and use the construction just proved in one coset. In other
cosets use a fixed translation. Returning to the original offsets changes
both signs by the same factor and leaves equal multiplicities. Otherwise
every four-subset of the five is additive. For three fixed elements
\(a,b,c\), the possible completions are \(a+b-c,a+c-b,b+c-a\). The two
other elements are distinct choices from these three. Relabeling and
translating \(a\) to zero, the five elements can be written as
\(\{0,x,-x,c,c+x\}\). The two completions of \(\{0,x,-x\}\) must then be
\(2x,-2x\). If \(c=-2x\), then \(c+x=-x\) repeats a point. Hence
\(c=2x\), and \(3x=-2x\) shows \(5x=0\). The five elements are a
translate of a group of order five. On that group the permutations
\(i\mapsto2i\) and \(i\mapsto-i\) both use every shift once, since the
differences from \(i\) are \(i\) and \(-2i\). The first has a four-cycle
on the nonzero elements, and the second has two transpositions there.
Their signs are opposite, as required.

Now suppose \(g\ge2\). A cancellation in one connected block remains a
cancellation after filling the other blocks with the same matching. For
at most two shifts, all shift orbits have the same size, so the sign is
still determined by the total multiplicities. It remains to consider
repeated blocks for three shifts and for additive four-shift supports.

For three shifts on \(H\), write the set as \(\{0,s,t\}\). Let
\(\Lambda=\{(x,y):xs+yt=0\text{ in }H\}\). We show that every lattice
point \((r,z)\) with \(r,z\ge0\) and \(r+z\le m\) occurs as a total
displacement. The zero point comes from the identity. For a nonzero point
write \((r,z)=k(a,b)\) with \((a,b)\) primitive in \(\Lambda\), and write
\((a,b)=e(p,q)\), where \(p,q\) are relatively prime. Completing
\((a,b)\) to a basis of \(\Lambda\) gives a second vector whose
determinant with \((a,b)\) is \(\pm m\). Hence \(\ell(x,y)=py-qx\) takes
exactly the values \((m/e)\mathbb Z\) on \(\Lambda\).

Start in the strip \(0\le\ell<p+q\), at level zero. Take an east step if
\(\ell\ge q\) and a north step otherwise. The new level is \(\ell-q\) or
\(\ell+p\), respectively. Both rules are addition of \(p\) modulo \(p+q\)
within the strip. The levels are distinct until \(p+q\) steps have been
made, since \(\gcd(p,p+q)=1\). The resulting displacement has level zero,
is parallel to \((p,q)\), and has coordinate sum \(p+q\). It is therefore
\((p,q)\). After \(e\) periods the displacement is \((a,b)\), so the walk
closes modulo \(\Lambda\). No earlier repetition is possible. Equal
cosets within the strip must have equal levels, since the width is at
most \(m/e\), and equal levels in different periods are separated by a
multiple of \((p,q)\). The first positive multiple in \(\Lambda\) is
\(e(p,q)\). Thus the projected walk is a simple cycle. Every integer is a
possible starting level since \(p,q\) are coprime. Take \(k\) such strips
starting at \(0,p+q,\ldots,(k-1)(p+q)\). Their total width is at most
\(m/e\), by \(r+z\le m\). Levels from different strips are therefore
distinct modulo \(m/e\), and their projected cycles are disjoint. Fixing
all other vertices realizes \((r,z)\).

If \(u,v\) form a lattice basis and \(u,v,u+v\) lie in this triangle, use
displacements \(0,u+v\) in two blocks, and \(u,v\) in the same two blocks
for a second term. Both terms have the same total shift counts. Each of
\(u,v,u+v\) is primitive, so the first term has one nontrivial cycle and
the second has two. The signs are opposite.

We check when such a basis can be found. Suppose a difference \(s\) has
order \(h<m\), and translate the support to \(\{0,s,t\}\). Put \(d=m/h\).
The image of \(t\) generates \(H/\langle s\rangle\) of order \(d\). Thus
any lattice vector has second coordinate divisible by \(d\). There is a
unique \(0\le v<h\) with \(vs+dt=0\), so \((h,0),(v,d)\) form a basis.
Their sum is in the triangle if \(h+v+d\le hd\). For \(d\ge3\) this
follows from \(v\le h-1\) and \((h-1)(d-2)\ge1\). For \(d=2\) the only
failure is \(v=h-1\). Then \((h-1)s+2t=0\) gives \(2t=s\), so the support
is a progression.

Otherwise every difference generates \(H\). It follows that \(H\) is
cyclic, and its order is odd, since three residues modulo an even number
cannot all have odd pairwise differences. We can take the support as
\(\{0,1,q\}\) with \(2\le q\le(m+1)/2\), using \(q\mapsto1-q\) if needed.
The endpoints \(q=2,(m+1)/2\) give progressions. For the other values
write \(m=kq+v\), \(0\le v<q\). The vectors \((v,k),(v+q,k-1)\) belong to
\(\Lambda\), since each gives \(x+qy=m\). Their determinant is \(-m\), so
they form a basis. Their sum is in the triangle precisely when
\[
 v+q-1\le k(q-2).
\]
For \(k\ge3,q\ge4\), the left side is at most \(2q-2\), and the right
side is at least \(3q-6\). For \(k\ge3,q=3\), the only failure is
\(k=3,v=2\), namely \(m=11,q=3\). Replacing \(q\) by \((1-q)^{-1}=5\)
modulo eleven gives a basis satisfying the inequality. For \(k=2\) a
failure requires \(v=q-1\) or \(v=q-2\). If \(m=3q-1\), replace \(q\) by
\(q^{-1}=3\) and use the preceding cases. If \(m=3q-2\), replace \(q\) by
\(1-q^{-1}\), whose representative is \((m-1)/2\). Then the new remainder
is one, and the inequality holds for \(m\ge9\). These replacements are
obtained by permuting and rescaling the three offsets, so they preserve
the sign question.

The remaining support is \(\{0,1,3\}\) in \(\mathbb Z_7\). Listing the
nonzero points \(r,z\ge0\) with \(r+z\le7\) and \(r+3z\equiv0\bmod7\)
gives \((7,0),(0,7),(4,1),(1,2),(2,4)\). For these points the numbers of
nontrivial cycles are \(1,1,1,1,2\) by \eqref{eq:additive-sign}. The
numbers of moved vertices are \(7,7,5,3,6\), so all signs are positive.
This remains true for any number of equal blocks.

For a three-point progression translate its middle shift to zero. The
nonzero steps are then the two neighboring steps on a cycle of length
\(m\). A nontrivial permutation cycle is a transposition or a full
rotation in one direction. Let \(u,v\) count the blocks using the two
rotations, and let \(j\) count the transpositions. The two nonzero shift
counts are \(mu+j,mv+j\), and the sign is \((-1)^{(m-1)(u+v)+j}\). At
fixed shift counts, any changes satisfy \(\Delta u=\Delta v=-\Delta
j/m\). For even \(m\) both \(\Delta j\) and \(\Delta u+\Delta v\) are
even, so the sign is fixed. For odd \(m\) in two blocks we have \(0\le
j\le m-1\); a change of \(j\) by a multiple of \(m\) is possible only if
it is zero. With three blocks, the two rotations and an identity give the
shift counts \((m,m)\) and positive sign. The same counts come from \(m\)
transpositions spread over three blocks. For \(m=2a+1\), take \(a,a,1\)
transpositions, which fit in those blocks. Their sign is negative. This
gives cancellation for \(g\ge3\), using identity permutations in any
extra blocks.

Finally consider repeated blocks for \(\{0,s,t,s+t\}\). The triples
\(\{0,s,t\}\) and \(\{0,s,-t\}\) occur up to translation and generate
\(H\). A cancellation in two blocks for either triple also works for the
four shifts. Otherwise both triples are progressions or one is the
order-seven case. If both are progressions, zero cannot be the middle
point of \(\{0,s,t\}\), since \(s+t\ne0\). Interchanging \(s,t\) if
needed gives \(t=2s\). The other triple is then \(\{0,s,-2s\}\); its
midpoint condition gives \(4s=0\) or \(5s=0\). Since \(s\) generates
\(H\) and the four offsets are distinct, \(m=4\) or \(5\).

For \(m=4\) label \(H\) by \(0,1,2,3\). Compare the identity in one block
and \((3,2,0,1)\) in the other with the pair \((2,1,0,3),(3,1,2,0)\), in
one-line notation. The four-cycle in the first pair has shifts
\(3,1,2,2\), while the two transpositions in the second have shifts
\(2,2\) and \(3,1\) on their moved vertices. Both pairs give
\(a^4bc^2d\). The first pair has negative sign and the second positive.

In the order-seven case replace \(t\) by \(-t\) if necessary, which
translates the four-set. With a change of generator the exceptional
triple is \(\{0,1,t\}\) with \(t=3\) or \(5\). For \(t=3\), multiplying
the four-set by two and adding one gives \(\{0,1,2,3\}\). For \(t=5\),
adding two does so. It therefore remains to treat consecutive four shifts
at odd \(m=2j+1\ge5\), including \(m=5\).

In two blocks use the pairs of monomials
\[
 (a^m,\ b^2 c^{m-4}d^2),\qquad
 (a^{j+1}c^{j-1}d,\ a^j b^2 c^{j-2}d).
\]
Their products agree. The second monomial in the first pair has winding
number two, and reversal of offsets makes it a winding-one monomial with
\(m-2\) moved vertices. Since \(m\) is odd, the column shift in the
reversal has positive sign, as does the resulting cycle of odd length
\(m-2\). In the second pair both monomials have winding number one, with
\(j\) and \(j+1\) moved vertices, respectively. Their signs multiply to
\(-1\). The occurrence of all these monomials follows from
\eqref{eq:four}, or directly by the winding-one compositions. Thus there
is cancellation in two blocks, and additional blocks can be filled with
the identity. This completes the classification.

\subsection{The signed and block-constant cases of Theorem~\ref{thm:block}}

For square matrices \(X,Y\) of the same order, expand each product in
\(\per(X+Y)\) according to the rows that use an entry of \(Y\). If those
rows are \(I\) and their columns are \(J\), the remaining rows must use
the remaining columns. Hence
\[
 \per(X+Y)=\sum_{|I|=|J|}
                 \per X[I^c,J^c]\per Y[I,J].
\]
Apply this identity to the diagonal and forward parts of \(M\). Let
\(I_i\) be the rows of block \(i\) using its forward block, and \(J_i\)
the columns of block \(i\) filled by a forward block. Since \(J_i\) is
filled from \(I_{i-1}\), we have \(|J_i|=|I_{i-1}|\). The remaining rows
and columns at block \(i\) use \(A_i\), so \(b-|I_i|=b-|J_i|\). Thus
every \(|I_i|\) and \(|J_i|\) equals a common \(k\). For fixed sets, the
sum of the matching weights is
\begin{equation}\label{eq:block-expansion}
 \prod_{i=0}^{m-1}
       \per A_i[I_i^c,J_i^c]\per B_i[I_i,J_{i+1}],
\end{equation}
with indices read modulo \(m\).

We first compare a minor of a single block with its signed determinant.
Let its order be \(r\) and let \(d,d+1\) be the two internal offsets. For
a contributing matching, list its rows in increasing order. The lifted
targets are \(j+d\) or \(j+d+1\). As \(j\) increases, the targets are
nondecreasing. Equality would repeat a column, and the difference between
the first and last targets is at most \(b\). Equality in this last bound
would also repeat a column modulo \(b\). Thus the targets are strictly
increasing and span less than \(b\). If \(e\) of them are at least \(b\),
they form a suffix. On reducing modulo \(b\), each of these \(e\) targets
moves ahead of all the other \(r-e\) targets. There are \(e(r-e)\)
inversions and no others, so the term has sign
\((-1)^{e(r-e)}=(-1)^{e(r-1)}\). Multiplying each wrapped entry by
\((-1)^{r-1}\) cancels that sign. Therefore
\[
 \per D[I,J]=\det D^{((-1)^{r-1})}[I,J],
 \qquad |I|=|J|=r,
\]
with both values one when \(r=0\).

Now expand the determinant \(F_\epsilon(x)\) by the same sets \(I_i,J_i\)
of size \(k\). In a block, moving a chosen row set \(I_i\) to the end
while preserving its order takes
\[
 \sum_{j\in I_i}(b-1-j)-\binom{k}{2}
\]
interchanges. The first sum counts every row after a chosen row, and the
subtraction removes the pairs of chosen rows. The columns give the same
expression with \(J_i\). The terms depending only on \(b,k\) occur twice
in every block and cancel in parity. The total factor from these row and
column orderings is thus \((-1)^{\sum_i(\sum I_i+\sum J_i)}\). In the
forward minors, the left and right factors \(E\) supply exactly the same
sign, since the forward columns in \(B_i\) are \(J_{i+1}\). The two
factors cancel.

After this reordering, separate the diagonal groups from the forward
groups in both rows and columns. These two shuffles have the same sign
and cancel. In the forward part the column groups are cyclically shifted
by one group of size \(k\). Moving such a group past the other \(m-1\)
groups has sign \((-1)^{(m-1)k^2}=(-1)^{(m-1)k}\). This cancels the sign
from multiplying the last forward block by \((-1)^{m-1}x\), which also
supplies \(x^k\).

Choose \(\epsilon=(-1)^{k-1}\). The twist in a forward minor is then the
one for order \(k\). In a diagonal minor it is
\((-1)^b\epsilon=(-1)^{b-k-1}\), the one for order \(b-k\). The
coefficient \([x^k]F_\epsilon(x)\) is therefore the sum of
\eqref{eq:block-expansion} over the sets of that size. For odd \(k\) use
\(F_+\), and for even \(k\) use \(F_-\). The odd coefficients of \(F_+\)
sum to \((F_+(1)-F_+(-1))/2\); the even coefficients of \(F_-\) sum to
\((F_-(1)+F_-(-1))/2\). Adding them proves \eqref{eq:block-four-det}.

For a block circulant with two offsets, a column shift makes one of them
zero. Order the block indices along the cycles of addition by the other
offset. Rows from one cycle can use only columns in that cycle. The
matrix is therefore a direct sum of cyclic block matrices of the kind
just considered. For block Hankel matrices, reversing the order of the
block columns changes a block depending on \(i+j\) into one depending on
\(i-j\), with no change to its internal row and column order. The
permanent is preserved by this column permutation.

For the banded block-constant case in Corollary~\ref{cor:block}, the rows
inside a single block are identical. At a cut between blocks, only the
last \(q\) row blocks can have entries into later blocks. Their
restricted rows span a space of dimension at most \(q\). At a cut inside
a block there is at most one further row type from that block. The upper
ranks are therefore at most \(q+1\). Transposition gives at most \(p+1\)
for the reverse cut ranks. Theorem~\ref{thm:path} gives the determinant
assertion.

Finally, for \eqref{eq:constant-blocks} fix the numbers \(k_{ij}\) of
rows in class \(i\) using columns in class \(j\). Within row class \(i\),
the number of ways to assign the column classes is \(m_i!/\prod_j
k_{ij}!\). For column class \(j\), its \(n_j\) columns can be assigned to
the selected rows in \(n_j!\) ways. Every such matching has weight
\(\prod_{i,j}b_{ij}^{k_{ij}}\). Multiplying these factors and summing
over the arrays with the required row and column sums gives
\eqref{eq:constant-blocks}. The entries outside the first \(r-1\) rows
and first \(s-1\) columns of the array are determined by its margins.
This gives the stated bound on the number of arrays.

\subsection{The open-diagonal and Hankel formulas}

For \eqref{eq:open-padding}, embed \(A\) as the first diagonal block of
\(A\oplus I_{N-n}\). Any contributing permutation matches the two blocks
separately. Its restriction to the added block is the identity, with
weight one, so the permanent is \(\per A\). The allowed ordinary offsets
of \(A\), reduced modulo \(N\), give its cyclic offsets in the enlarged
matrix. The only added entries are diagonal and thus allowed when \(0\in
S\). When \(N=n\) there are no added entries. By the hypotheses on these
cyclic offsets, Corollary~\ref{cor:circulant}(3) applies. For the four
displayed offsets, the chosen \(N\) makes \(d_1+d_2\equiv d_3+d_4\pmod
N\). If the residues are distinct, subtracting \(d_1\) makes the second
offset the sum of the third and fourth offsets. If any coincide there are
at most three offsets. Since each \(d_i\) lies between \(1-n\) and
\(n-1\), the chosen \(N\) is \(O(n)\). This gives the stated cost in
these cases.

For \eqref{eq:additive-open}, at least one of \(q/g,p/g\) is odd, so the
stated choice of \(x,y\) satisfies \(x^{q/g}y^{-p/g}=-1\). Work in one
residue class modulo \(g\). On the grid, the map \((u,v)\mapsto pu+qv\)
has kernel generated by \((q/g,-p/g)\). Its image is \(g\mathbb Z\). Thus
the grid modulo that kernel represents the vertices of the residue class,
up to an initial translation, and the east, north and northeast edges are
the three nonzero offsets. The finite open support retains only those
edges with both ends in the matrix.

The quotient is a cylinder. Each nontrivial cycle of a contributing
permutation gives an embedded simple curve in it. Its lifted displacement
is nonzero, since the steps increase at least one coordinate, and is a
multiple of the kernel generator. A simple essential curve on a cylinder
has winding number \(1\) or \(-1\). One can see this by identifying the
cylinder with the punctured plane and using the index of a simple closed
curve about the missing point. Thus the product of the \(x,y\) factors on
the cycle is \((x^{q/g}y^{-p/g})^{\pm1}=-1\). The permutation sign of a
cycle of length \(a\), together with its \(a\) off-diagonal minus signs,
is another factor \((-1)^{a-1}(-1)^a=-1\). They cancel. A fixed point has
its original weight, so the determinant and permanent agree term by term
on each residue class and hence on the whole matrix.

For the rational series, write one of them as \(\sum_{j\ge0}\gamma_j
z^j=P(z)/Q(z)\) with \(Q(0)=1\). Take \(d\ge1\) at least \(\deg Q\) and
larger than \(\deg P\), and write \(Q(z)=1+q_1z+\cdots+q_dz^d\), using
zeros as needed. The coefficients in \(Q(z)\sum_j\gamma_jz^j\) give
\[
 \gamma_j=-\sum_{a=1}^d q_a\gamma_{j-a}
 \qquad (j\ge d).
\]
The vector \((\gamma_j,\ldots,\gamma_{j+d-1})^{\mathsf T}\) therefore
advances by a fixed matrix, with its last row \((-q_d,\ldots,-q_1)\) and
its other rows shifting the coordinates. Its first coordinate is
\(\gamma_j\). This gives the constant generators on both sides of
\eqref{eq:toeplitz-rational}, also for polynomial series or repeated
roots. For a cut after \(k\), for example,
\[
 t_{j-i}=(p^{\mathsf T}B^{k-i})
                 (B^{j-k-1}q)\qquad (i\le k<j).
\]
Thus the upper cut block has rank bounded by the order of \(B\). The
lower block is treated with the other constant generator. Since the
diagonal weight is also constant, \eqref{eq:border} gives
\((T_n^n)_{1,1}\). The bound on the degrees of the generators
fixes the number of coordinate variables, while their allowed
degrees grow with \(n\). The adjugate formula for the
constant-state entry of \((I-zT_n)^{-1}\) gives \eqref{eq:rationaltoeplitz}.

In characteristic \(p>0\), the same companion constructions give
constant generators. Corollary~\ref{cor:positive-ordered} uses the
fixed basis \(\mathcal B_{p,r,s}\), so its common matrix \(T_p\)
is independent of the order. Thus \(\per A_n=(T_p^n)_{1,1}\) for
every \(n\geq0\). Summing the formal geometric series gives
\[
 \sum_{n\geq0}\per A_n z^n=((I-zT_p)^{-1})_{1,1}.
\]
The cofactor identity proves \eqref{eq:positive-toeplitz}. The
denominator has constant coefficient one. Repeated squaring gives
the stated field-operation bound.

For the characteristic-three example, on \(x^hy^h\) the diagonal
coefficient of the ordered operator is \((1+2h)\tau^h\), the
coefficient at degree \(h+1,h+1\) is \(\tau^{h+1}\), and that
at degree \(h-1,h-1\) is \(h^2\tau^{h-1}\). Reading these for
\(h=0,1,2\) in characteristic three gives
\eqref{eq:positive-geometric-matrix}. The term \(x^3y^3\) is
zero in the quotient and cannot contribute at a later step by
Lemma~\ref{lem:positive-ideal}. The two determinants in
\eqref{eq:positive-toeplitz} give \eqref{eq:positive-geometric}.

For the geometric parameters choose \(B=p,C=q\),
\(D_i=\partial_x,E_i=\partial_y\), \(L_i=ax,M_i=by\). On a monomial
\(x^hy^h\), the substitution alone multiplies it by \(\tau^h\). The two
single derivatives followed by insertion of \(ax\) or \(by\) give
\(haq\tau^{h-1}\) and \(hbp\tau^{h-1}\) at degree \(h,h\), when
\(h\ge1\). The double derivative gives \(h^2\tau^{h-1}\) at degree
\(h-1,h-1\). The two other terms give \(c\tau^h\) at degree \(h,h\) and
\(ab\tau^h\) at degree \(h+1,h+1\). These are precisely the entries of
\eqref{eq:geometric}. At \(h=0\) the derivatives vanish. At the upper end
the projection removes degree \(N+1,N+1\).

For \(I-zT_n\), expanding a leading tridiagonal minor in its last row gives
a continuant relation. If \(\Delta_h\) is the leading determinant through
index \(h\), then for \(h\ge1\)
\[
 \Delta_h=(1-z(T_n)_{h,h})\Delta_{h-1}
            -z^2(T_n)_{h-1,h}(T_n)_{h,h-1}\Delta_{h-2},
 \qquad \Delta_{-1}=1,\quad \Delta_0=1-z(T_n)_{0,0}.
\]
Thus the determinants for \(T_n,T_n'\) can be formed by products of
two-by-two polynomial matrices of degree at most two. Multiply
consecutive groups of factors in a balanced tree. At a level with groups
of length \(a\), there are \(O(n/a)\) products of degree \(O(a)\). The
assumption on \(M(n)/n\) gives \(O(M(n))\) operations at each level.
Inverting the denominator series by doubling the precision and
multiplying by the numerator also lies within \(O(M(n)\log(n+1))\).

For the cycle polynomial, at the cut after vertex \(k\), give a path with
end \(i\) and start \(j\) the factors \(ap^{k-i}\) and \(bq^{k-j}\).
These are the weights of the edges from \(i\) to \(k+1\) and from \(k+1\)
to \(j\). A retained end gains a factor \(p\) at the next cut and a
retained start gains a factor \(q\). Suppose there are \(h\) old paths.
If the new vertex is isolated, it contributes \(ab\) from its new
endpoints, with \(\tau^h\) from the old ones. A loop instead contributes
\(zc\tau^h\). If an edge from an old end enters the new vertex, that end
is consumed and the new end contributes \(a\). There are \(h\) choices
and the retained endpoints gain \(p^{h-1}q^h\), giving \(haq\tau^{h-1}\).
The reverse choice gives \(hbp\tau^{h-1}\). In the two-edge case the old
end and start can lie on the same path in \(h\) ways, closing a cycle, or
on distinct paths in \(h(h-1)\) ways, joining them. In both cases \(h-1\)
paths remain and the retained endpoints gain \(\tau^{h-1}\). The
coefficient for this step is \(h(h-1+z)\tau^{h-1}\). These cases give
\(T_n(z)\). A complete cover starts and ends with no paths. Taking the
corresponding entry of \(T_n(z)^n\), then using the same adjugate identity,
proves \eqref{eq:geometric-cycles}.
For the full coefficient list, \(P_{A_n}\) has degree at most \(n\).
Evaluating it at \(0,1,\ldots,n\) and interpolating gives the stated
bound for all its coefficients.

For \((\rho^{|i-j|})\), the factorization
\(\rho^{j-i}=\rho^{k-i}\rho^{j-k}\) at a cut shows that the two cut ranks
are one for \(\rho\ne0\). For its determinant, subtract \(\rho\) times
row \(i-1\) from row \(i\), for \(i=n,\ldots,2\). In row \(i\) this zeros
columns \(j<i\), and the new diagonal entry is \(1-\rho^2\). The first
diagonal entry remains one. For the odd-row, even-column submatrix, write
its entries as \(\rho^{|2a-1-2b|}\), \(1\le a,b\le m\). Subtract
\(\rho^2\) times column \(b-1\) from column \(b\) in decreasing order of
\(b\). When \(a<b\), the original entry in column \(b\) is \(\rho^2\)
times the one in column \(b-1\). The diagonal entries after subtraction
are \(\rho\) and \(\rho(1-\rho^2)\) for \(b\ge2\). This gives the stated
minor and is unaffected by changing the full matrix diagonal.

For the weighted clique-width bound, consider the binary disjoint-union
tree of an expression with \(n\ge2\) vertices. There is a proper
subexpression with between \(n/3\) and \(2n/3\) vertices. Indeed,
starting at the root, descend through a child with more than \(2n/3\)
vertices while one exists. At the final union at least one of its two
children has at least \(n/3\) vertices and both have at most \(2n/3\).
Vertices with the same label in that child can never be separated by
subsequent relabelings. They therefore get equal weights to each fixed
vertex outside the child. For \(0<\rho<1\), a weight \(\rho^{|i-j|}\)
determines the distance to that fixed vertex \(j\), and at most two
vertices have a given distance. Thus each label accounts for at most two
of the at least \(n/3\) vertices. This gives at least \(n/6\) labels.

Finally, for an anti-diagonal support, \((H\mathsf
J_n)_{ij}=H_{i,n-1-j}\). If \(H_{ij}=h_{i+j}\), the new entry is
\(h_{n-1+i-j}=t_{j-i}\), with \(t_d=h_{n-1-d}\). The same index
calculation carries arbitrary entries on an anti-diagonal into the
corresponding ordinary diagonal. This proves the transfer and padding
formulas for \(H\). For the direct signing, \(K\) is the signed matrix
for \(H\mathsf J_n\), so \(\per H=\det K\). The reversal permutation has
\(n(n-1)/2\) inversions. Consequently \(\det(K\mathsf
J_n)=(-1)^{n(n-1)/2}\det K\), giving the signed Hankel formula. The
special cases of Corollary~\ref{cor:hankel} use the same substitution in
the three displayed Toeplitz families. In positive characteristic,
apply Corollary~\ref{cor:positive-ordered} to the reversed matrix.
This proves Theorem~\ref{thm:hankel}(4).

\subsection{Proof of Proposition~\ref{prop:symmetry}}

In a contributing permutation for \(S\), the first \(m\) rows must use
the last \(m\) columns and the last \(m\) rows must use the first \(m\)
columns. These are independent choices, so \(\per S=\per B\per B^{\mathsf
T}=(\per B)^2\). For \(C\) each set of rows uses the columns of its own
block. Both blocks have permanent \(\per B\) by row and column reversal,
which gives \(\per C=(\per B)^2\). Conjugating \(C\) by \(\mathsf
J_{2m}\) reverses both blocks and exchanges them, so \(C\) is
centrosymmetric.

Both blocks of \(D=S\oplus(\mathsf J_{2m}S\mathsf J_{2m})\) are
symmetric. Reversing the full row and column orders exchanges them after
internal reversal, so \(D\) is also centrosymmetric. Its permanent is
\((\per S)^2=(\per B)^4\). For \((0,1)\) matrices these values determine
\(\per B\) by taking the nonnegative integer square or fourth root,
respectively.

\end{document}